\documentclass[preprint,3p,10pt]{elsarticle}

\usepackage{graphicx,color} % for Inkscape
\usepackage{amssymb}
\usepackage{amsthm}
\usepackage{amsmath}
\usepackage{subfigure}
\usepackage{stmaryrd}
\usepackage{multirow}
\usepackage{pgfplots}
\pgfplotsset{compat=newest}
\usetikzlibrary{plotmarks}
\pgfplotsset{plot coordinates/math parser=false}
\newlength\figureheight 
\newlength\figurewidth
\tikzset{every picture/.append style={scale=0.8}}
\usetikzlibrary{shapes} % for symbols in figure caption

\newcommand{\Eref}[1]{Equation (\ref{#1})}
\newcommand{\fref}[1]{Figure (\ref{#1})}
\newcommand{\Tref}[1]{Table (\ref{#1})}

\newtheorem{remark}{Remark}

\newcommand{\bb}{\mathbf{b}}

\newcommand{\bn}{\boldsymbol{n}}
\newcommand{\bt}{\mathbf{t}}
\newcommand{\bu}{\mathbf{u}}
\newcommand{\bv}{\mathbf{v}}
\newcommand{\bx}{\mathbf{x}}
\newcommand{\bs}{\boldsymbol{\tau}}
\newcommand{\bsig}{\boldsymbol{\sigma}}
\newcommand{\be}{\boldsymbol{\epsilon}}
\newcommand{\bB}{\mathbf{B}}

\newcommand{\barB}{\bar{\mathbf{B}}}
\newcommand{\bM}{\boldsymbol{M}}
\newcommand{\bV}{\mathbf{V}}
\newcommand{\RR}{\mathbb{R}}
\newcommand{\0}{\mathbf{0}}

\newcommand{\An}{\mathrm{A}_{\theta}}
\newcommand{\AUn}{\mathrm{A}_{1}}
\newcommand{\Pn}{\mathrm{P}_{\theta}}
\newcommand{\PnUn}{\mathrm{P}_{1}}

\usepackage{lineno}

\usepackage{hyperref}

\journal{Computer Methods in Applied Mechanics and Engineering}

\usepackage{setspace}
\begin{document}
\begin{frontmatter}

\title{Skew-symmetric Nitsche's formulation in isogeometric analysis: \\
Dirichlet and symmetry conditions, patch coupling and frictionless contact}

\author[dut,lu1]{Qingyuan Hu}
\ead{huqingyuan@mail.dlut.edu.cn}
%\ead{qingyuanhucn@gmail.com}

\author[be]{Franz Chouly}
\ead{franz.chouly@univ-fcomte.fr}

\author[dut]{Ping Hu}
\ead{pinghu@dlut.edu.cn}

\author[dut]{Gengdong Cheng}
\ead{chenggd@dlut.edu.cn}

\author[lu1,be,cu1,uwa1]{St\'ephane P.A. Bordas\corref{cor1}}
%\ead{stephane.bordas@gmail.com}
\ead{stephane.bordas@alum.northwestern.edu}
\cortext[cor1]{Corresponding author}

%\address[dut]{Department of Engineering Mechanics, Dalian University of Technology, Dalian 116024, P.R. China}
\address[dut]{State Key Laboratory of Structural Analysis for Industrial Equipment, Dalian University of Technology, P.R. China}
\address[lu1]{Institute of Computational Engineering, Faculty of Sciences Technology and Communication, University of Luxembourg, Luxembourg}
\address[be]{Laboratoire de Math\'ematiques de Besan\c{c}on - UMR CNRS 6623, Universit\'e Bourgogne Franche-Comt\'e, 16 route de Gray, 25030 Besan\c{c}on cedex, France}
\address[cu1]{Institute of Mechanics and Advanced Materials, School of Engineering, Cardiff University, UK}
\address[uwa1]{Intelligent Systems for Medicine Laboratory, University of Western Australia, Perth, Australia}

\begin{abstract}
A simple skew-symmetric Nitsche's formulation is introduced into the framework of isogeometric analysis (IGA) to deal with various problems in small strain elasticity:
essential boundary conditions,
symmetry conditions for Kirchhoff plates,
patch coupling in statics and in modal analysis
as well as Signorini contact conditions.
For linear boundary or interface conditions, the skew-symmetric formulation is parameter-free.
For contact conditions,
it remains stable and accurate for a wide range of the stabilization parameter.
Several numerical tests are performed to illustrate its accuracy, stability and convergence performance.
We investigate particularly the effects introduced by Nitsche's coupling, including the convergence performance and condition numbers in statics as well as the extra ``outlier'' frequencies and corresponding eigenmodes in structural dynamics. 
We present the Hertz test, the block test, and a 3D self-contact example
showing that the skew-symmetric Nitsche's formulation is a suitable approach to simulate contact problems in IGA.
\end{abstract}

\end{frontmatter}

%\linenumbers

%\section*{Highlights}

%$\bullet$ A simple skew-symmetric Nitsche's formulation is introduced into IGA.

%$\bullet$ For linear boundary and interface conditions, skew-symmetric Nitsche is parameter-free.

%$\bullet$ For contact conditions, skew-symmetric Nitsche %remains stable and accurate for a wide range of values for the stabilization is robust w.r.t. the stabilization parameter. 

%$\bullet$ Robustness and accuracy of the method is studied numerically for various problems.

%%%%$\bullet$ In statics and dynamics, the coupling effects introduced by Nitsche's method is studied

%%%%%$\bullet$ The ``outlier'' frequencies introduced by Nitsche's patch-coupling for structural vibration formulations is captured

\section{Introduction}

The key concept in isogeometric analysis (IGA) \cite{hughes2005isogeometric} consists in using non-uniform rational B-splines (NURBS) as basis functions to approximate both the geometry and the unknown physical fields.
The mathematical foundations of IGA are developed in \cite{beiraodavega-2014}, and a recent overview is given in \cite{nguyen2015isogeometric}.
Contrary to classical Lagrange basis functions usually adopted in the finite element method (FEM), NURBS in IGA have the ability to exactly describe geometries: thus, no geometrical approximation error is introduced.
Moreover NURBS are widely adopted in commercial computer-aided design (CAD) packages,
and this CAD data can directly be used to construct approximations.
In boundary element method (BEM), this translates into the ability to solve directly from the field variables at the control points defining the geometry \cite{simpson2012two,simpson2013isogeometric,scott2013isogeometric,lian2013stress,beer2013isogeometric,marussig2015fast,lian2016implementation,peng2017isogeometric}.
In FEM, a 3D parameterization of the volume is still necessary \cite{xu2011parameterization,xu2013optimal},
except when solving shell-like problems \cite{kiendl2009isogeometric,benson2010isogeometric,benson2011large,echter2013hierarchic,hu2017isogeometric}.
The present paper focuses on two following issues.
One first issue in IGA is related to boundary conditions, especially essential boundary conditions.
Indeed, since NURBS are non-interpolatory,
enforcing boundary conditions and constraints cannot be done as simply as in Lagrange FEM:
they require tackling difficulties which are similar to those encountered in meshless methods \cite{nguyen2008meshless} and implicit/immersed boundary methods \cite{dunant2007architecture,harari2010analysis}.
One second issue in IGA comes from interface conditions and patch coupling: for complex geometries, patch-wise CAD modeling is necessary, and transmission conditions need to be satisfied.
The same also arises when gluing heterogeneous materials.

Various methods already exist to treat boundary or interface conditions weakly,
that have been firstly designed for instance in the FEM context. They are applicable,
or have already been applied, for IGA. The most widespread ones are the penalty method, mixed/mortar methods and Nitsche's method.
The penalty method \cite{babuska-73a,kikuchi-oden-88} is simple but not consistent.
Therefore the value of the penalty parameter has to be chosen with great care to achieve the best balance between accuracy and stability.
As a matter of fact, if the penalty parameter is chosen too small the boundary or interface conditions are imposed inaccurately,
whereas if it is chosen much larger than needed the penalized problem becomes ill-conditioned.
Mixed methods for boundary conditions \cite{babuska-73b} introduce a Lagrange multiplier,
which is an additional variable that represents the boundary stress,
and that allows to take into account weakly the essential boundary conditions in a consistent way.
This leads to a weak problem that has a saddle-point structure.
For patch-coupling,
the original mortar method \cite{bernardi-maday-patera-93,bernardi-maday-patera-94} has been reformulated later as a mixed/dual Lagrange multiplier method (see, e.g., \cite{benbelgacem-99,wohlmuth2000mortar} for FEM and \cite{brivadis-2015} for IGA).
%The original mortar method is improved by a dual Lagrange multiplier space that has better theoretical and practical properties \cite{wohlmuth2000mortar,brivadis-2015}.
Mortar methods, when carefully designed, are consistent, stable and optimally accurate (see, e.g., \cite{wohlmuth2000mortar} in the FEM context or \cite{brivadis-2015} in the IGA context). Moreover the newly introduced Lagrange multipliers have a clear meaning: they are the stresses needed to enforce the continuity of the displacements.
Mortar techniques have been applied as well with success to contact problems 
\cite{temizer2011contact,temizer2012three,de2012mortar,kim2012isogeometric,seitz-2016,antolin2017textit}.
However extra degrees of freedoms (DoFs) are introduced and an inf-sup condition must be fulfilled in order to ensure stability and optimal convergence, for which care is needed to build the dual space of Lagrange multipliers. % obtain a positive definite stiffness matrix.

Nitsche's method was originally proposed by J. Nitsche \cite{nitsche-71,stenberg-95} to impose weakly essential boundary conditions and   more recently has regained popularity to deal with interface conditions with non-conforming discretizations (see, e.g., \cite{becker-hansbo-stenberg-03,hansbo-05,annavarapu2012robust}).
Nowadays Nitsche's method has also found a number of natural applications in IGA \cite{embar2010imposing,nguyen2014nitsche,apostolatos2014nitsche,ruess2014weak,guo2015nitsche,du2015}.
%For mathematical analysis of Nitsche's method, please refer to \cite{becker-hansbo-stenberg-03,hansbo-05}.
Nitsche's formulation %turns the constraints into
makes use of an appropriate \emph{conjugate pair} such as displacement--force or rotation--moment,
in such a way that the method remains both primal (no extra DoFs) and consistent.
By the way, there is no need to fulfill an inf-sup condition.
However, standard (symmetric) Nitsche's method includes an extra term that penalizes the boundary/interface conditions and allows to recover stability and optimal accuracy. 
For this purpose, this extra term makes use of an additional numerical parameter,
the \emph{stabilization parameter}, that needs to be fixed %large enough  
above a given threshold.
For simple problems and numerical methods, such as piecewise linear or quadratic Lagrange FEM, a direct and accurate estimation of the aforementioned threshold can be effectuated (see, e.g., \cite{hansbo-larson-2002,hansbo-05} for a discussion on this topic), but for more realistic problems and less standard numerical methods, this can be harder to achieve. %an {\it a priori} knowledge of this threshold can be harder to achieve.
%Of course, an alternative technique remains 
Indeed this threshold for the stabilization will depend upon many parameters, related to the physical constants (Young's modulus) and to the discretization (polynomial order of basis functions, shape of the cells in the grid): see, e.g., \cite{annavarapu2012robust,jiang2015robust,shahbazi2005explicit}.
%Then %one common approach to estimate the stabilization parameter 
In such situations an alternative to estimate this threshold consists in solving a generalized eigenvalue problem along the target boundary/interface: see, e.g., \cite{hansbo-larson-2002,hansbo-05} for FEM and \cite{embar2010imposing,nguyen2014nitsche,apostolatos2014nitsche} for IGA.
%To address this difficult of the standard Nitsche's method lead to the inception of 
Nevertheless, the difficulties associated to this issue can be circumvented by using the penalty-free (skew-symmetric) variant of
Nitsche's method, such as in \cite{burman-12,kollmannsberger2015parameter,boiveau2016penalty,burman-hansbo-larson-2016-pfree,schillinger2016non}.

In this paper we present a simple and systematic procedure to derive,
for various boundary and interface conditions,
a family of Nitsche's formulations that have different symmetry properties and different degrees of dependency on the stabilization parameter.
This family is indexed by the \emph{Nitsche parameter} $\theta$.
%A parameter $\theta$ indexes this family.
We then focus on the variant known as the skew-symmetric Nitsche's method, that corresponds to the value $\theta=-1$.
This method can be parameter-free when dealing with linear boundary or interface conditions, and reveals to be very robust with respect to the stabilization parameter for non-linear boundary conditions such as contact.
Let us mention that in the context of standard FEM, the skew-symmetric method has been successfully applied to contact \cite{chouly2013nitsche,chouly2014adaptation,chouly2015symmetric,chouly-overview-2017}.
Furthermore in IGA there is already one contribution dealing with the skew-symmetric Nitsche's method for enforcing Dirichlet boundary conditions and patch coupling in the context of thin shell problems \cite{guo2016parameter}.
In this contribution we perform numerical experiments for different situations,
particularly we study how Nitsche's multi-patch coupling can affect the accuracy, the convergence rates, and the condition numbers.
Moreover, in modal analysis, literature \cite{cottrell2006isogeometric,cazzani2016analytical} shows that some outlier frequencies appear due to the discretization of the continuous problem.
This ``outlier'' phenomenon is also captured in multi-patch cases using the mortar method \cite{horger2017improved}.
Here we study this
issue of the ``outlier'' frequencies and corresponding eigenmodes,
in the context of Nitsche's method.
Finally to our knowledge Nitsche's method has never been applied in IGA for contact conditions,
and we show how to implement Nitsche's formulation for contact problems,
and how it performs in these cases.

The outline of this paper is as follows.
In Section \ref{section_iga} the concept and notations of IGA are introduced, 
the critical differences between Lagrange-based FEM and NURBS-based IGA are also explained.
In Section \ref{section_nitsche_franz} we introduce the Nitsche-based formulations for boundary/interface conditions,
starting from an abstract setting.
In Section \ref{section_example} various numerical tests are performed and we reach conclusions in Section \ref{section_conclusion}.
%Supplementary details are provided in Appendix.

\section{Brief introduction to isogeometric analysis}\label{section_iga}

Bivariate NURBS basis functions $R_A(\xi,\eta),(A=1,\cdots,nm)$ are often adopted in IGA to generate surfaces. 
They are constructed using appropriate weights $w_A$ and the tensor product of two sets of univariate B-spline basis functions $N_{i,p}(\xi),(i=1,\cdots,n)$ and $N_{j,q}(\eta),(j=1,\cdots,m)$,
where $p$ and $q$ are orders of the B-spline basis functions in directions $\xi$ and $\eta$ respectively.
One set of B-spline basis functions can be calculated from one given knot vector recursively \cite{cottrell2009isogeometric}.
By the help of NURBS basis functions, the desired surface is represented as the set $\overline{\Omega}$ of points
\[
\boldsymbol{x} (\xi,\eta) =\sum_{A=1}^{nm} R_A (\xi,\eta) \, \boldsymbol{x}_A,
\]
where $\boldsymbol{x}_A(x,y,z)$ denote positions of the control points. %We will denote by $h$ the mesh size associated to this parametrization of the surface $\overline{\Omega}$.
Following the ``iso'' concept, any (discrete) physical field $\bu^h$ defined on the surface (domain) $\overline{\Omega}$ is represented using the same set of NURBS basis functions as
\[
\bu^h (\xi,\eta) =\sum_{A=1}^{nm} R_A (\xi,\eta) \, \bu_A,
\]
where $\bu_A$ are the control point variables, as well as the degrees of freedom associated to $\bu^h$.
In the following we will denote %by $h$ the global mesh size and 
by $\bV^h$ the finite dimensional space of such discrete fields $\bu^h$ constructed using IGA (see, e.g., \cite{beiraodavega-2014} for the detailed construction of such a space). The notation $h$ will stand for the size of the cells associated to such a discretization.
Note that no essential boundary or interface conditions are prescribed in the definition of $\bV^h$.
According to \cite{xu2014geometry} the spline spaces used for the geometry and the physical field can be chosen and adapted independently,
which is known as the Geometry-Independent Field approximaTion (GIFT)
and brings more flexibility in the field approximation when preserving geometric exactness and tight CAD integration.
However the present research is restricted to IGA.

\begin{figure}[htbp]
\centering
\def\svgwidth{1.0\columnwidth}
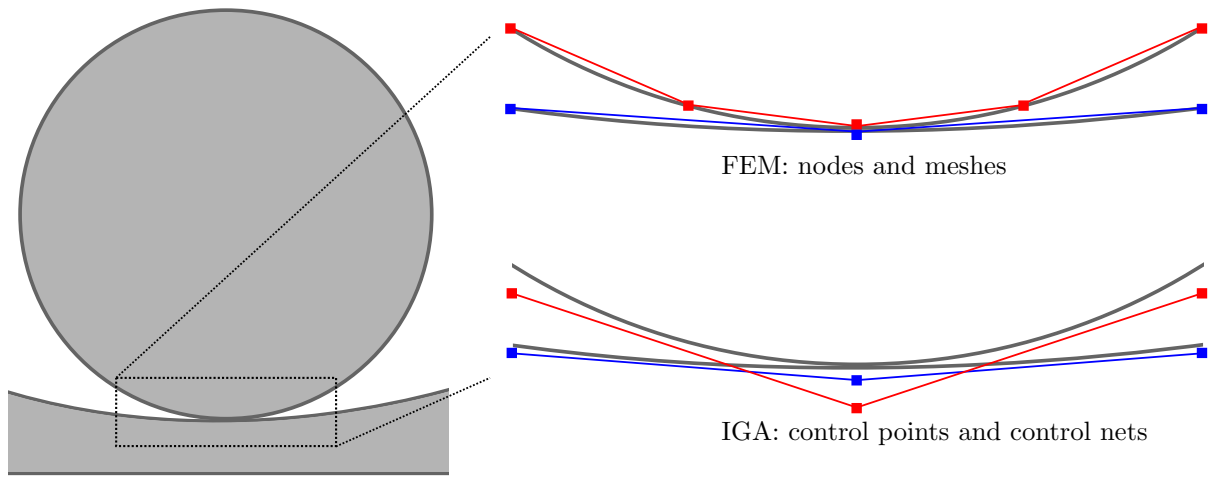
\caption{Boundary discretization: Lagrange basis function based FEM and NURBS basis functions based IGA.}
\label{lagrange_nurbs}
\end{figure}

In \fref{lagrange_nurbs} some differences between FEM and IGA are illustrated.
Consider a contact problem in the left of the figure,
and for simplicity we just consider the discretization of a portion of the boundaries that are going into contact.
For Lagrange basis function based FEM,
%since the Lagrange basis functions are interpolated at both ends,
%the FEM nodes are located inside the domain,
discrete errors are introduced by the FEM meshes.
%which does not match exactly the curved domain.
For NURBS basis functions based IGA,
the curved domain is parametrized exactly.
%for instance by using only one isogeometric element of order $p=2$ with three control points as shown in the picture.
The NURBS basis functions of order $p$
%may have a desired higher order derivatives up to continuity $C^{p-1}$.
allow to represent fields up to continuity $C^{p-1}$,
%however (most) control points are not located inside the domain,
%located inside the domain,
however (most) control points that define the boundary/interface are not interpolated,
which is owing to the corresponding non-interpolating basis functions.
%In more general cases,
This brings difficulties in directly manipulating the control variables attached to these control points when dealing with boundary and interface conditions.
In the next section we are going to introduce Nitsche's formulation to impose various boundary/interface conditions weakly.

\section{Nitsche's formulation for boundary/interface conditions}\label{section_nitsche_franz}

We first present Nitsche's method within an abstract setting,
and then show how this framework can be applied to recover various well-known Nitsche-based discretizations,
for a wide range of problems in computational mechanics.
Note that for linear boundary and interface conditions, discretized with finite elements, a general presentation can be found in, e.g., 
\cite{stenberg-95,becker-hansbo-stenberg-03,hansbo-05}.
We will consider a whole family of Nitsche's methods indexed by a real value, that we will call the \emph{Nitsche parameter} $\theta \in \mathbb{R}$, and we will pay particular attention to the skew-symmetric variant, i.e. to the case $\theta=-1$.

\subsection{Abstract setting}

Consider the domain $\Omega$ as the open set associated to the surface $\overline{\Omega}$ defined in previous Section \ref{section_iga}.
We will denote by $\Gamma$ either a portion of the boundary of $\Omega$ or an interface that subdivides $\Omega$ into two subdomains.
Our aim is to compute a field $\bu : \Omega \rightarrow \RR^d$ ($d\geq1$),
for instance a displacement field,
that is a solution to a given set of partial differential equations with prescribed boundary/interface conditions.
To simplify the presentation we consider a linear partial differential equation (but not necessarily linear boundary or interface conditions).
We will denote by $\bv$ an arbitrary test function,
that can represent the virtual displacement.
Two main ingredients are necessary to build a Nitsche-based formulation.

The first ingredient is a Green formula (inspired by Theorem 5.8 in \cite{kikuchi-oden-88}), that allows to rewrite weakly the partial differential equation satisfied by $\bu$, and
that we provide below in an abstract setting:
\begin{equation}\label{nitsche1green}
\textrm{Find }\bu \in \bV \,:\, a(\bu , \bv) - \langle \bs(\bu) , \bB (\bv) \rangle_\Gamma = L(\bv),
\quad \forall \,  \bv \in \bV,
\end{equation}
where $\bV$ is a functional space of admissible fields, $a(\cdot \,, \cdot)$ is a bilinear form (the internal work), $\langle \cdot\,, \cdot \rangle_\Gamma$ is an appropriate duality product for functions on $\Gamma$ (the boundary/interface work), and $L(\cdot)$ a linear form (the work of external loads). 
The linear operator $\bB$ is a trace-like operator: for instance $\bB (\bv)$ can be the value of $\bv$ on $\Gamma$, or of its normal component if $\bv$ is a vector field.
The dual quantity $\bs(\bu)$, where $\bs$ is a flux-like operator, is to be defined for each situation.
It is generaly related to the boundary/interface stress, if $\bu$ is a displacement (generalized stress vector in elasticity).
We can call $\bs(\bu)$ and $\bB (\bv)$ a \emph{conjugate pair}.
We suppose that both $\bs(\bu)$ and $\bB (\bv )$ can be represented at almost every point of the boundary as vectors of dimension $k$ ($1 \leq k \leq d$):
$$\bs(\bu) : \Gamma \rightarrow \RR^k, \quad \bB(\bv) : \Gamma \rightarrow \RR^k.$$

The second ingredient is a reformulation of the boundary/interface conditions as follows:
\begin{equation}\label{nitsche2proj}
 \bs(\bu) = \left [ \bs(\bu) - \gamma ( \bB ( \bu ) - \barB ) \right ]_S.
\end{equation}
In the above formula,
$\barB$ is a known, prescribed, quantity associated to the trace. The notation 
$\left [ \cdot  \right ]_S$ stands for the projection onto $S$, a closed subset of $\RR^k$ of admissible values. The set $S$ can depend on $\bu$ in some situations ($S = S(\bu)$), for instance in the case of Coulomb friction \cite{chouly-overview-2017,chouly-coulomb-18,renard-13}, but we omit this dependence to simplify our notations. For simple problems, as those studied in this paper, $S$ is generally a closed convex set, a closed convex cone, or a subspace of $\RR^k$.
Finally $\gamma$ is an arbitrary positive and one-to-one Schur operator (interface or boundary stiffness), that transforms a trace into a flux.

Nitsche-based discretizations can be obtained by following the steps we describe below in more detail and which are mathematically valid only for sufficiently smooth fields $\bu$ and $\bv$:

\begin{enumerate}
 \item Apply the following decomposition
 \[
  \bB (\bv)  = -\gamma^{-1} \left ( \theta \bs( \bv ) - \gamma \bB ( \bv ) \right ) + \theta \gamma^{-1} \bs(\bv).
 \]
 \item Insert it into \eqref{nitsche1green}, which yields
 \[
   a(\bu\, , \bv) - \theta \langle \bs(\bu) , \gamma^{-1} \bs(\bv) \rangle_\Gamma  
   +  \langle  \bs(\bu) , \gamma^{-1} ( \theta \bs( \bv ) - \gamma \bB ( \bv ) ) \rangle_\Gamma = L(\bv).
 \]
\item Inject condition \eqref{nitsche2proj} into the above formula, so as to impose it weakly
 \begin{equation} \label{nitsche3cont}
   a(\bu\, , \bv) - \theta \langle \bs(\bu) , \gamma^{-1}\bs(\bv) \rangle_\Gamma  
   + \langle [ \bs(\bu) - \gamma ( \bB (\bu) - \barB ) ]_S , \gamma^{-1} ( \theta \bs( \bv ) - \gamma \bB (\bv) ) \rangle_\Gamma = L(\bv).
 \end{equation}
\end{enumerate}

%\begin{comment} %{\textwidth}
%{\bf (Strong--weak equivalence)} Formally, the above calculations prove that \eqref{nitsche1green} and \eqref{nitsche2proj} imply \eqref{nitsche3cont}. In fact, it can (hopefully) be proved that, formally \eqref{nitsche1green}--\eqref{nitsche2proj} and \eqref{nitsche3cont} are equivalent \red{(detail)}.
%\end{comment}

The above formula may have no meaning at the continuous level. Nevertheless it becomes meaningful once all the fields are discretized. For this purpose, consider $\bV^h$, a discrete space, built from any Galerkin approximation, such as finite elements or IGA (see Section \ref{section_iga} above). For discrete fields the duality pairing $\langle \cdot, \cdot \rangle_\Gamma$ becomes simply the scalar product in $L^2(\Gamma)$ , and we will denote by $\| \cdot \|_{L^2(\Gamma)} (= \langle \cdot , \cdot \rangle^{\frac12}_\Gamma)$  the corresponding norm. Let us consider $$\gamma_h : L^2(\Gamma) \rightarrow L^2(\Gamma)$$ a discrete Schur operator, positive and one-to-one. Consider $\bu^h$ (resp. $\bv^h$) a discrete approximation to $\bu$ (resp. to $\bv$).
To simplify the notations, we introduce also the modified discrete weak form
\[
 \An ( \bu^h  , \bv^h ) :=  a(\bu^h , \bv^h) - \theta \langle \bs(\bu^h) , \gamma_h^{-1} \bs(\bv^h) \rangle_\Gamma, 
\]
and the linear operator
\[
 \Pn ( \bv^h ) := \theta \bs( \bv^h ) - \gamma_h \bB (\bv^h).
\]
Then we obtain the Nitsche-based formulation below
\begin{equation} \label{nitschefem}
\textrm{Find } \bu^h \in \bV^h \, : \, 
 \An ( \bu^h , \bv^h ) + \langle [ \PnUn (\bu^h) + \gamma_h \barB ]_S , \gamma_h^{-1} \Pn ( \bv^h ) \rangle_\Gamma = L(\bv^h), \quad \forall \, \bv^h \in \bV^h.
\end{equation}
Remark that the way the method is built ensures its consistency with respect to the partial differential equation being solved. 
% We can state already a first result of consistency:
% %\begin{proposition}
% The method \eqref{nitschefem} is consistent: for a regular enough solution $\bu$ to \eqref{nitsche1green}--\eqref{nitsche2proj}, there holds:
% \begin{equation} \label{nitscheconst}
%  \An ( \bu \, ; \bv^h ) + \langle [ \PnUn (\bu) + \gamma \barB ]_S , \gamma_h^{-1} \Pn ( \bv^h ) \rangle_\Gamma = L(\bv^h), \quad \forall \bv^h \in \bV^h.
% \end{equation}
% %\end{proposition}

% {\it Proof:} If $\bu$ is solution to \eqref{nitsche1green} and \eqref{nitsche2proj}, then it is also solution to \eqref{nitsche3cont} which makes sense for $\bv^h \in \bV^h$ and for $\bu$ regular enough, and taking $\gamma = \gamma_h$. This proves \eqref{nitscheconst}. $\Box$

An important particular case is that of linear boundary/interface conditions, which means that $S=\mathbb{R}^k$ in \eqref{nitsche2proj} and so the projection operator is merely the identity.
Then \eqref{nitschefem} reads
\[
 \An ( \bu^h ,  \bv^h ) + \langle \PnUn (\bu^h) + \gamma_h \barB , \gamma_h^{-1} \Pn ( \bv^h ) \rangle_\Gamma = L(\bv^h), 
\]
and, after re-ordering and simplifications we arrive at
\begin{eqnarray}
 & & a ( \bu^h , \bv^h) 
 - \langle \bs ( \bu^h ) , \bB (\bv^h) \rangle_\Gamma
 - \theta \langle \bs ( \bv^h ) , \bB (\bu^h)  \rangle_\Gamma 
 +  \langle \gamma_h \bB (\bu^h) , \bB (\bv^h) \rangle_\Gamma \nonumber \\
 &=& L ( \bv^h ) - \langle \theta \bs ( \bv^h ) - \gamma_h \bB (\bv^h) , \barB \rangle_\Gamma.
 \label{nitschefemlinear}
\end{eqnarray}
When $\theta=1$, we recover the well-known formulation presented for instance in \cite{stenberg-95,becker-hansbo-stenberg-03,hansbo-05}.

\begin{remark} The Nitsche parameter $\theta$ allows to select some variants of Nitsche's formulation, that yield different theoretical properties and different degrees of dependency w.r.t. the operator $\gamma_h$:

\begin{itemize}
 \item for $\theta=1$,
the standard symmetric Nitsche's method \cite{nitsche-71} is obtained. If,  $a(\cdot , \cdot)$ is symmetric, and under appropriate assumptions on $S$ (for instance if $k=1$ and $S=\RR$, $S=\RR^-$ or $S=\RR^+$), it can be derived as the first order optimality condition of the energy functional \cite{nitsche-71,chouly-overview-2017,rabii-numath-accepted}:
\[
\mathcal{J}_{\mathrm N} (\bu^h) :=
\frac12 \AUn ( \bu^h , \bu^h ) - L (\bu^h)
+ \frac12 
\langle [ \PnUn (\bu^h) + \gamma_h \barB ]_S , \gamma_h^{-1} 
[ \PnUn (\bu^h) + \gamma_h \barB ]_S  \rangle_\Gamma.
\]
Moreover a suitable choice for $\gamma_h$ is necessary in order to 
recover well-posedness and optimal accuracy (see Section \ref{sub:Schur} below);

 \item for $\theta=0$, some terms cancel out and we obtain the simple formulation
\begin{equation*}
 a ( \bu^h , \bv^h ) - \langle [ P_1(\bu^h) + \gamma_h \barB ]_S , \bB ( \bv^h ) \rangle_\Gamma = L(\bv^h),
\end{equation*}
which is close to an augmented lagrangian formulation and easier to extend to the large strain framework 
%that is more suitable for extensions to large elastic transformations 
\cite{mlika2017unbiased};

\item for $\theta=-1$,
the skew-symmetric Nitsche's method is obtained, see e.g., \cite{freund1995weakly,burman-12,boiveau2016penalty} for linear boundary conditions and \cite{chouly2015symmetric,chouly-overview-2017} for contact conditions.
Stability and optimal convergence are ensured whatever $\gamma_h$.
Note that for linear boundary/interface conditions,
we can even choose $\gamma_h= 0$,
resulting in the parameter-free formulation:
\begin{equation}\label{nitschefemlinearfree}
 a ( \bu^h , \bv^h) 
  - \langle \bs ( \bu^h ) , \bB (\bv^h) \rangle_\Gamma  
  + \langle \bs ( \bv^h ) , \bB (\bu^h)  \rangle_\Gamma 
 = L ( \bv^h ) + \langle \bs ( \bv^h ) , \barB \rangle_\Gamma.
\end{equation}
\end{itemize}

\end{remark}

\subsection{The discrete Schur operator, well-posedness and optimal accuracy}
\label{sub:Schur}

When $\theta \neq -1$, the discrete Schur operator $\gamma_h$ needs to be designed so as to preserve well-posedness of Problem \eqref{nitschefem} as well as optimal accuracy. Let us first consider a linear setting, i.e. Problem \eqref{nitschefemlinear} and $\theta=1$: the key issue in the mathematical analysis (see, e.g., \cite{thomee-97,chouly-overview-2017}) is to ensure the $\bV^h$-ellipticity of the bilinear form $\AUn (\cdot\,,\cdot)$ in the energy norm. %an appropriate discrete norm on $\bV^h$.
To this purpose, we define
\begin{equation} \label{constantTI}
C_{\mathrm{TI}} (h) := \sup_{\bv^h \in \bV^h} \frac{ %\| \bs( \bv^h ) \|^2_{-\frac12,h,\Gamma}
\| \bs(\bv^h) \|^2_{L^2(\Gamma)}
}{a(\bv^h,\bv^h)}
\end{equation}
the \emph{trace-inverse} constant associated to \eqref{nitschefemlinear}, %where $\| \cdot \|_{-\frac12,h,\Gamma}$ is a weighted discrete norm on $\Gamma$ and $C_{\mathrm{TI}} (h)$ is a constant, 
that depends on the size $h$ of the cells, but also on the other features of the discrete space $\bV^h$, such as the polynomial order of basis functions. This constant depends also upon the partial differential equation under consideration (for instance it depends on the Young's modulus in isotropic linear elasticity). Suppose that, for instance, $\gamma_h$ is such that
\begin{equation}\label{condGamma}
{\| \gamma_h^{-1} \|} {C_{\mathrm{TI}} (h)} \leq \frac12
%\| \gamma_h^{-1} \| \geq \ldots %2 C_{\mathrm{TI}} (h),
\end{equation}
with $\| \gamma_h^{-1} \| = \sup_{\tau \in L^2(\Gamma), \| \tau \|_{L^2(\Gamma)} = 1} \| \gamma_h^{-1} \tau \|_{L^2(\Gamma)}$. %an appropriate (discrete) operator norm (related to $\| \cdot \|_{-\frac12,h,\Gamma}$), 
For any $\bv^h \in \bV^h$, we can write
\begin{eqnarray*}
\AUn ( \bv^h  , \bv^h )
&\geq  &a(\bv^h , \bv^h) 
- {\| \gamma_h^{-1} \|} \| \bs( \bv^h ) \|^2_{L^2(\Gamma)}\\%\langle \bs(\bv^h) , \bs(\bv^h) \rangle_\Gamma |,\\
&\geq & \left ( 1 - {\| \gamma_h^{-1} \|} {C_{\mathrm{TI}} (h)} \right ) a(\bv^h , \bv^h) \geq \frac12 a(\bv^h , \bv^h).
\end{eqnarray*}
The same kind of argument holds for the general formulation \eqref{nitschefem} and for any value of $\theta \neq -1$ (see, e.g., \cite{chouly-overview-2017} in the case of small strain elasticity with contact).
For simple situations, as considered in this paper, where the coefficients of the partial differential equation are constant, and where the mesh is quasi-uniform, the simplest choice is to define the Schur operator globally as 
\[
\gamma_h := \gamma_0\, \mathbf{Id} 
\]
where $\gamma_0 > 0$ is a real parameter, the \emph{stabilization} parameter, and $\mathbf{Id}$ is the identity in $L^2(\Gamma)$. Then condition \eqref{condGamma} is reformulated as
\[
\gamma_0 \geq 2 C_{\mathrm{TI}} (h).
\]
As in \cite{embar2010imposing,nguyen2014nitsche,apostolatos2014nitsche} the constant $C_{\mathrm{TI}}(h)$ can be estimated as the maximum eigenvalue $\lambda^{h,\textrm{MAX}}$ associated to the problem
\begin{equation} \label{eigenvalue}
\textrm{Find } (\lambda^h,\bu^h) \in \RR \times \bV^h \, : \, 
\langle \bs(\bu^h) , \bs(\bv^h) \rangle_{\Gamma} = \lambda^h \, a ( \bu^h , \bv^h )
\quad \forall \, \bv^h \in \bV^h,
\end{equation}
and then one can choose $\gamma_0 = 2 \lambda^{h,\textrm{MAX}}$. This is actually what we do in the numerical experiments of Section \ref{section_example}. % for most of the numerical experiments (except in situations when the value of $\gamma_0$ can be inferred easily from the setting or is the object of study).

Of course, if more general situations need to be considered, such as more general (non quasi-uniform) meshes or partial differential equations with spatially variable coefficients, it is much better to consider a local, cell-wise, definition of $\gamma_h$: in this case it is chosen piecewise constant on each mesh cell $K$ and a local counterpart of \eqref{eigenvalue} can be solved to recover the value of $\gamma_h|_K$ (see, e.g., \cite{hansbo-larson-2002,hansbo-05} in the context of FEM). %TODO Question - i suppose, that you mean non-constant coefficients in space, right? Otherwise, considering the cell-wise problem would not necessarily make sense.

%For simple problems and numerical methods, such as piecewise linear or quadratic Lagrange FEM, a direct and accurate estimation of the aforementioned threshold can be effectuated (see, e.g., \cite{hansbo-larson-2002,hansbo-05} for a discussion on this topic), 

%Indeed this threshold for the stabilization will depend upon many parameters, related to the physical constants (Young's modulus) and to the discretization (polynomial order of basis functions, shape of the cells in the grid): see, e.g., \cite{annavarapu2012robust,jiang2015robust,shahbazi2005explicit}.

%In such situations an alternative to estimate this threshold consists in solving a generalized eigenvalue problem along the target boundary/interface: see, e.g., \cite{hansbo-larson-2002,hansbo-05} for FEM and \cite{embar2010imposing,nguyen2014nitsche,apostolatos2014nitsche} for IGA.
%

In case where $\theta \neq -1$, and in the context of Lagrange FEM, provided that a condition such as \eqref{condGamma} is satisfied, both stability and optimal accuracy in the energy norm can be established:
%When Nitsche formulation is combined with finite elements,
%stability and optimal accuracy can be proved: 
see, e.g., \cite{stenberg-95,thomee-97,becker-hansbo-stenberg-03,hansbo-05} for the complete mathematical analysis in the linear setting, and \cite{chouly2015symmetric,chouly-overview-2017} for unilateral contact with Tresca Friction.
Moreover standard scaling arguments show that the constant $C_{\mathrm{TI}} (h)$ scales as $\mathcal{O}(h^{-1})$, irrespectively of the polynomial order of the FEM \cite{thomee-97}.
%For this purpose,
%the stabilization parameter $\gamma$ should scale as $h^{-1}$,
%where $h$ is the size of the elements,
%and be chosen above a given threshold $Ch^{-1}>0$.
%The constant $C$ in this threshold comes from the application of a trace-inverse inequality that allows to bound locally the boundary/interface flux by the energy norm.
%This constant $C$ depends on 1) some physical constants such as the Young modulus and 2) the polynomial order of the finite element space.
In the skew-symmetric case $\theta=-1$,
the condition %$\gamma \geq Ch^{-1}$ 
\eqref{condGamma} can be relaxed,
and it suffices to take $\gamma_0 >0$, or even $\gamma_0=0$ for linear boundary/interface conditions (``penalty-free" variant).
A complete mathematical analysis for $(\theta,\gamma_0)=(-1,0)$ can be found in, e.g., \cite{burman-12} for Poisson's problem and in \cite{boiveau2016penalty} for compressible and incompressible elasticity.
The same results as above can be expected for the IGA setting though no numerical analysis has been provided to the best of our knowledge.
%See however \cite{apostolatos2014nitsche} for a heuristics to set $\gamma$ in the symmetric case $\theta=1$.

In the remaining part of this paper, we will focus on the skew-symmetric variant $\theta=-1$,
%because of its robustness properties,
but numerical tests with the symmetric variant $\theta=1$ are also performed for comparison purposes.

\subsection{Nitsche's formulation for linear boundary conditions}
\label{sub:linearbc}

We first illustrate how the above framework can be applied to deal with some linear boundary conditions, thus we consider the case where $\Gamma$ is a subset of $\partial \Omega$. The unit normal vector on $\Gamma$ pointing outward of $\Omega$ is denoted by $\bn$.
%Without losing generality, we begin with a second order problem in 2D, in which a simple elastic body occupies the domain $\Omega$.

\subsubsection{Dirichlet boundary conditions in small strain elasticity}
\label{subsub:Dirichlet}

\begin{figure}[htbp]
\centering
\def\svgwidth{0.5\columnwidth}
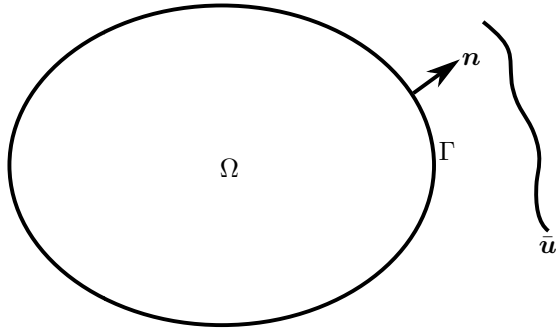
\caption{Dirichlet boundary condition: the displacement is equal to $\bar{\bu}$ on $\Gamma$.}
\label{Dirichlet_problem}
\end{figure}

Consider a linear elastic body described by a small strain constitutive model, that 
is subjected to body forces $\boldsymbol{b}$,
and surface loads $\bar{\boldsymbol{t}}$ along a Neumann boundary $\Gamma_N \subset \partial \Omega$.
%We want to find the displacement field $\bu$ associated to this problem.
The corresponding governing equations read
\begin{equation} \label{strongelast}
\begin{split}
-\nabla \cdot \boldsymbol{\sigma} (\bu) &= \boldsymbol{b} \qquad \textrm{ in $\Omega$}, \\
\boldsymbol{\sigma} (\bu) \bn &= \bar{\boldsymbol{t}} \qquad \textrm{ on $\Gamma_N$},
\end{split}
\end{equation}
where $\bu$ is the unknown displacement field, $\nabla \cdot$ is the divergence operator for vector-valued functions, and $\bsig$ is the Cauchy stress tensor. For the sake of simplicity, we choose to model the elastic behavior using Hooke's law and we denote by $E$ the Young's modulus and by $\nu$ the Poisson's ratio.
The corresponding weak form reads
\begin{equation} \label{greenelas}
a(\bu\,,\bv) - \int_\Gamma \bsig(\bu) \bn \cdot \bv \,ds = L(\bv),
\end{equation}
where
\begin{equation} \label{alelas}
 a (\bu\,,\bv) := \int_\Omega \bsig(\bu) : \be (\bv) \, d\bx, \quad 
 L(\bv) := \int_\Omega \bb \cdot \bv \, d\bx + \int_{\Gamma_N} \bar{\bt} \cdot \bv \, ds,
\end{equation}
and where $\be(\cdot)$ is the small strain tensor.
%As a result we recover a particular case of the Green formula 
The above formula \eqref{greenelas} matches with the general Green formula \eqref{nitsche1green}.
As illustrated Figure \ref{Dirichlet_problem}, we impose %in \eqref{greenelas} 
 an essential boundary condition on $\Gamma$: 
\begin{equation*}
\bu = \bar{\bu} \qquad \textrm{on $\Gamma$},
\end{equation*}
where $\bar{\bu}$ is the prescribed displacement.

% First we define the trial (resp. test) space of admissible displacements (resp. virtual displacements) as:
% \begin{equation}
% \begin{split}
% \bV &:= \left \{\bu\in (H^1(\Omega))^d \:\middle | \:\bu=\bar{\bu} \textrm{ on $\Gamma$} \right \}, \\
% \bV_0 &:= \left \{\bv \in (H^1(\Omega))^d \: \middle | \: \bv=\boldsymbol{0} \textrm{ on $\Gamma$} \right \}.
% \end{split}
% \end{equation}

With the choice
\[
 \bB (\bu) = \bu, \quad \bs(\bu) = \bsig(\bu) \bn, \quad \barB = \bar{\bu}, \quad S = \mathbb{R}^d,
\]
we obtain from \eqref{nitschefemlinear}
the following Nitsche-based formulation
\begin{eqnarray}\label{nitschefemDirichlet}
& & \textrm{Find } \bu^h \in \bV^h \, : \, \nonumber \\
& & a ( \bu^h , \bv^h)
- \int_\Gamma \left ( \bsig ( \bu^h ) \bn \right ) \cdot \bv^h \,ds
- \theta \int_\Gamma \bu^h \cdot \left (\bsig ( \bv^h )\bn \right ) \,ds
+ \int_\Gamma \gamma_h \, \bu^h \cdot \bv^h \,ds\\
& = & L ( \bv^h )
- \theta \int_\Gamma \bar{\bu} \cdot \left ( \bsig ( \bv^h ) \bn \right ) \,ds
+ \int_\Gamma \gamma_h \, \bar{\bu} \cdot \bv^h \,ds,
 \quad \forall \, \bv^h \in \bV^h. \nonumber
\end{eqnarray}
Setting $\theta=-1$ and $\gamma_h=0$, the penalty-free variant \cite{burman-12,ruess2014weak,boiveau2016penalty} is recovered.

\subsubsection{Symmetry conditions for Kirchhoff-Love plate}
\label{subsub:Kirchhoff}

\begin{figure}[htbp]
\centering
\def\svgwidth{0.6\columnwidth}
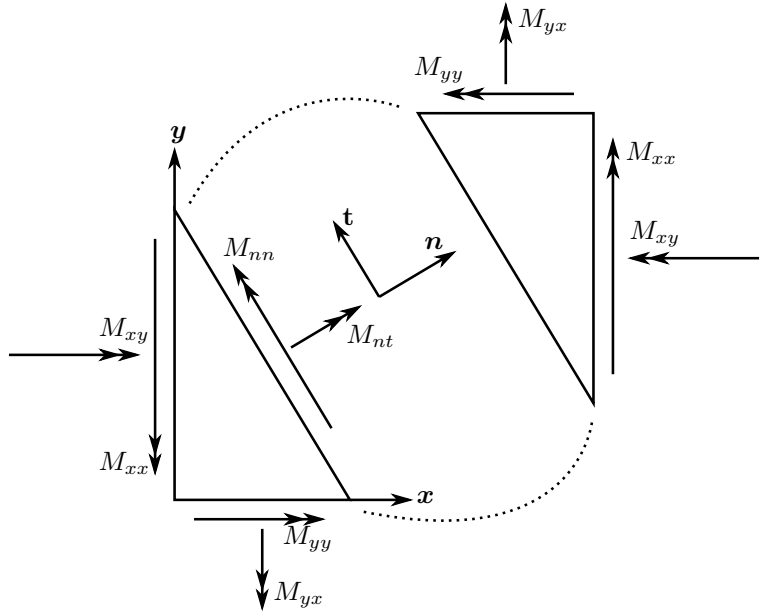
\caption{The directions of bending moments in Cartesian coordinate system $(\boldsymbol{x},\boldsymbol{y})$ and local system $(\bn,\bt)$.}
\label{bending_moment_fig}
\end{figure}

Thanks to the higher order continuity properties of NURBS basis functions,
%it experiences a renaissance 
there is a regained interest
to discretize thin-walled structures %by the 
using Kirchhoff-Love theory. %elements.
However due to the absence of rotational degrees of freedom,
additional effort is needed to apply rotational boundary conditions.
For this fourth-order problem, it is convenient to express the variables in local coordinates,
as illustrated in \fref{bending_moment_fig}.
Also, the corresponding weak form for Kirchhoff-Love plates reads
\begin{equation} \label{greenkirchhoff}
a(u,v) - \int_\Gamma M_{nn} (u) (- v_{,n}) \,ds = L(v),
\end{equation}
where $u$ is the deflection and $v$ the corresponding virtual quantity,
$\bn$ and $\bt$ indicate the outward normal direction and tangential direction respectively,
$M_{nn}(u)$ is the normal component of the moment tensor ($\bM(u) := -\, \mathcal{C} : \nabla^2 u$, where $\mathcal{C}$ is the constitutive fourth-order tensor) and $v_{,n} := (\nabla v) \cdot \bn$. In this case the bilinear and linear form read (see, e.g., \cite{harari2012embedded} for a more general formulation)
\[
 a (u,v) := - \int_\Omega \bM(u) : (\nabla^2 v) \, d\bx, \quad 
 L(v) := \int_\Omega f\,v \, d\bx,
\]
with $f$ a distributed load.
The symmetry condition on the boundary $\Gamma$ is formulated using the normal derivative of the mid-surface deflection $u$. More specifically, we impose:
\[
-u_{,n} = \bar{\theta}_t \qquad \textrm{on $\Gamma$},
\]
where $\bar{\theta}_t$ is a prescribed rotation.
Recall that Nitsche's contributions use \emph{conjugate pairs}:
in this case these are the rotation and the corresponding bending moment.
In order to form the Nitsche's contribution the rotation direction should be consistent with the direction of the corresponding bending moment (see \fref{bending_moment_fig}).

With the choice
\[
 \bB (\bu) = -u_{,n}, \quad \bs(\bu) = M_{nn}(u), \quad \barB = \bar{\theta}_t, \quad S = \mathbb{R},
\]
the Nitsche-based formulation is derived from \eqref{nitschefemlinear} (see as well \cite{embar2010imposing,harari2012embedded,schillinger2016non}):
\begin{equation}\label{bc_krichhoff}
\begin{split}
&\textrm{Find } u^h \in V^h\, :\\
&a(u^h,v^h)
-\int_{\Gamma} {M}_{nn}(u^h) \left(-v^h_{,n} \right) \,ds
-\theta\int_{\Gamma} \left(-u^h_{,n}\right) {M}_{nn}(v^h) \,ds
+ \int_{\Gamma} \gamma_h \left(-u^h_{,n}\right) \left(-v^h_{,n}\right) \,ds \\
=\:&L(v^h)
-\theta\int_{\Gamma}\bar{\theta}_t {M}_{nn}(v^h) \,ds
+ \int_{\Gamma} \gamma_h \bar{\theta}_t \left(v^h_{,n}\right) \,ds,
 \quad \forall \, v^h \in V^h.
\end{split}
\end{equation}
As previously, we recover a penalty-free method by setting $\theta=-1$, $\gamma_h=0$.

\subsection{Nitsche's formulation for interface conditions and patch coupling}
\label{sub:interface}

%\subsubsection{Interface conditions and patch coupling}

\begin{figure}[htbp]
\centering
\def\svgwidth{0.4\columnwidth}
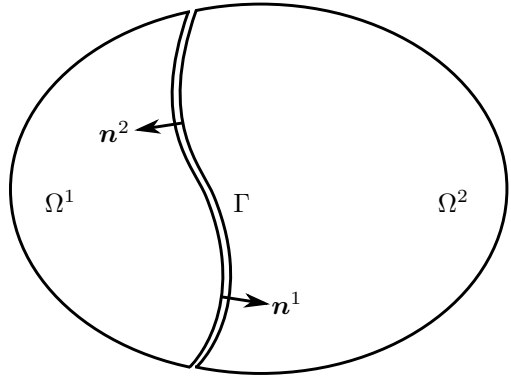
\caption{Problem with decomposed continuum domain.
Domain $\Omega$ is decomposed into two sub-domains $\Omega^1$ and $\Omega^2$.
The shared boundary is denoted by $\Gamma$ along which the outward unit normals are denoted by $\boldsymbol{n}^m$, $m=1,2$.}
\label{patch_coupling}
\end{figure}

Consider now an interface problem 
in which the domain $\Omega$ is decomposed into two sub-domains $\Omega^m$ (see \fref{patch_coupling}),
where the superscript $m=1,2$ is used to mark the partitioned domain and the corresponding variables.
The shared boundary between $\Omega^1$ and $\Omega^2$ is denoted by $\Gamma$, and
$\boldsymbol{n}^m$ is the unit normal along the interface $\Gamma$, pointing out of $\Omega^m$.
We still consider elasticity equations in small strains, and search for a displacement field $\bu = (\bu^1,\bu^2)$ solution to 
\[
\begin{split}
-\nabla \cdot \boldsymbol{\sigma} (\bu^m) &= \boldsymbol{b}^m \qquad \textrm{ in $\Omega^m$},\\
\boldsymbol{\sigma} (\bu^m) \bn^m &= \bar{\boldsymbol{t}}^m \qquad \textrm{ on $\Gamma_N^m$},\\
\bu^m &= \0 %\bar{\bu}^m 
\qquad\,\,\,\, \textrm{ on $\Gamma_D^m$},
\end{split}
\]
If Dirichlet boundary conditions on $\Gamma_D^m$ are non-homogeneous, they can be, for instance, treated as in Section \ref{subsub:Dirichlet} (but we omit this point to simplify notations).
In addition, there holds the following interface conditions %displacement and traction continuity conditions along the interface write
\[
\begin{split}
\bu^1 - \bu^2 &= \0 \qquad \textrm{on $\Gamma$}, \\
\boldsymbol{\sigma} (\bu^1) \boldsymbol{n}^1 + \boldsymbol{\sigma} (\bu^2)\boldsymbol{n}^2 &= \0 \qquad \textrm{on $\Gamma$}.
\end{split}
\]
The first equation corresponds to the continuity of the displacement along the interface, while the second one is the action-reaction principle. Note that this situation corresponds both to interface problems (when, for instance, material properties are different in $\Omega_1$ and $\Omega_2$) and to patch coupling (where the interface between subdomains is artificial, as in \cite{nguyen2014nitsche}). %The above equations are complemented by some essential and natural conditions on the external boundary $\partial \Omega$.

% By the definition of trial and test spaces containing the trial and test functions, respectively
% \begin{equation}
% \begin{split}
% S^m &= \{\bu^m(x)|\bu^m(x)\in H^1(\Omega^m),\bu^m=\bar{\bu}^m,\textrm{on $\Gamma_D^m$} \}, \\
% V^m_0 &= \{\bv^m(x)|\bv^m(x)\in H^1(\Omega^m),\bv^m=\boldsymbol{0},\textrm{on $\Gamma_D^m$} \},
% \end{split}
% \end{equation}
Let us define the jump and average operators along the interface $\Gamma$
\[
\begin{split}
\llbracket \bu \rrbracket &:= \bu^1-\bu^2, \\
\langle\boldsymbol{\sigma}(\bu) \rangle &:= \frac{1}{2}(\boldsymbol{\sigma}(\bu^1) \boldsymbol{n}^1 -\boldsymbol{\sigma}(\bu^2) \boldsymbol{n}^2).
\end{split}
\]
Let us introduce also
\begin{equation} \label{def_aL_interface}
\begin{split}
a(\bu\,,\bv) & :=  \sum_{m=1}^2 \int_{\Omega^m} \boldsymbol{\sigma}^m (\bu^m) : 
\boldsymbol{\epsilon}^m \left(\bv^m \right) \,d\bx,\\
L(\bv) & := \sum_{m=1}^2\int_{\Omega^m}
\boldsymbol{b}^m \cdot \bv^m \,d\bx +\sum_{m=1}^2 \int_{\Gamma_N^m} \boldsymbol{\bar{t}}^m 
\cdot \bv^m \, ds.
\end{split}
\end{equation}
Green formula for elasticity equations yields
\[
 a (\bu\,, \bv) - \int_\Gamma (\bsig(\bu^1) \bn^1) \cdot \bv^1 \, ds
 - \int_\Gamma (\bsig(\bu^2) \bn^2) \cdot \bv^2 \, ds = L(\bv).
\]
Remark now that the action-reaction principle implies the following identity
\[
 \bsig(\bu^1) \bn^1 = \frac12 \bsig(\bu^1) \bn^1 + \frac12 \bsig(\bu^1) \bn^1
 = \frac12 \bsig(\bu^1) \bn^1 - \frac12 \bsig(\bu^2) \bn^2 = \langle \bsig \rangle = -\bsig(\bu^2) \bn^2.
\]
This allows to impose weakly the action-reaction principle, as an essential interface condition, and we obtain
%It results that 
the appropriate Green formula in this context, as a particular form of \eqref{nitsche1green}: %, reads
\begin{equation} \label{greeninterface}
  a (\bu\, , \bv) - \int_\Gamma \langle \bsig(\bu) \rangle \llbracket \bv \rrbracket \,ds = L(\bv).
\end{equation}

With the choice
\[
 \bB (\bu) = \llbracket \bu \rrbracket, \quad \bs(\bu) = \langle \bsig(\bu) \rangle, \quad \barB = \0, \quad S = \mathbb{R}^d,
\]
formulation \eqref{nitschefemlinear} reads:
\begin{equation}\label{nitschefeminterface}
\begin{split}
&\textrm{Find } \bu^h \in \bV^h \, : \, \\
&a ( \bu^h , \bv^h) 
 - \int_\Gamma \langle \bsig(\bu^h) \rangle \llbracket \bv^h \rrbracket \,ds
 - \theta \int_\Gamma \llbracket \bu^h \rrbracket \langle \bsig(\bv^h) \rangle  \,ds
+ \int_\Gamma \gamma_h \, \llbracket \bu^h \rrbracket \llbracket \bv^h \rrbracket \,ds
 = L ( \bv^h ),\quad \forall \, \bv^h \in \bV^h.
\end{split}
\end{equation}
%Note that, here, Nitsche's method allows to recover weakly the continuity of the displacements.
Once again we recover a penalty-free formulation with $\theta=-1$, $\gamma_h=0$. 
Note as well that the same technique can be applied for patch coupling between other models, such as plates or rods, with the appropriate changes of notations \cite{nguyena2013nitsche}.

\subsection{Nitsche's formulation for frictionless contact conditions}
\label{sub:contact}

In this section we get back to the general (non-linear) formulation \eqref{nitschefem} and illustrate it in the case of frictionless contact, following \cite{chouly-overview-2017}. In \ref{subsub:biased} we first present Signorini contact and biased (master-slave) contact, in which the contact conditions are imposed on the boundary of one unique elastic body, in the same fashion as in Section \ref{sub:linearbc}. For contact between two bodies (or multi-body contact and self-contact), it is more convenient to impose contact conditions on both contact surfaces, using an unbiased formalism, as in, e.g., \cite{sauer-delorenzis-2015,rabii-numath-accepted,mlika2017unbiased}. So we
present in \ref{subsub:unbiased} an unbiased Nitsche's method, that can be derived in the same manner as presented in \ref{sub:interface} for interface conditions.

We still consider elastic bodies undergoing small strain and governed by Hooke's law, so we keep the same notations as in the previous sections, especially \ref{subsub:Dirichlet}.

\subsubsection{Biased frictionless contact conditions}
\label{subsub:biased}

\begin{figure}[htbp]
\centering
\def\svgwidth{0.45\columnwidth}
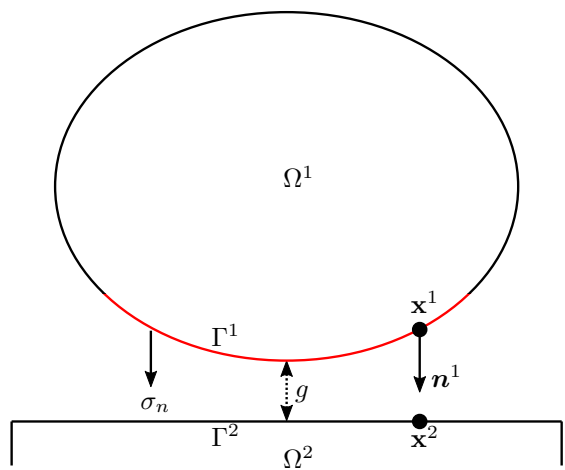
\caption{Contact problem setup, the contact slave surface is colored in red.}
\label{contact_problem}
\end{figure}

Consider a contact problem between an elastic body $\Omega := \Omega_1$ and a rigid support $\Omega_2$ as depicted in \fref{contact_problem}.
Surfaces for potential contact are denoted by $\Gamma := \Gamma_1$ and $\Gamma_2$.
To formulate the non-penetration condition, the normalized  vector is introduced
\[
\boldsymbol{n}^1 := \frac{\bx^2-\bx^1}{\left|\left|\bx^2-\bx^1 \right|\right|},
\]
where $\bx^1$ and $\bx^2$ are two mapped points on the corresponding boundary of each body $\Omega_1$ and $\Omega_2$,
for instance $\bx^2$ is the orthogonal projection of $\bx^1$ on $\Gamma_2$.
Then the gap function is defined as
\[
g := (\bx^2-\bx^1)\cdot \boldsymbol{n}^1.
\]
We deal with frictionless contact on $\Gamma$, so we impose weakly the essential condition $\bsig_{\bt} (\bu) = \0$, where $\bsig_{\bt}$ denotes the tangential stress. We start from equations \eqref{strongelast}
and we obtain the following Green formula, which is the counterpart of \eqref{greenelas} for Dirichlet boundary conditions:
\begin{equation} \label{greenelas2}
a(\bu\,,\bv) - \int_\Gamma \sigma_n(\bu) v_n \,ds = L(\bv),
\end{equation}
with the expression of $a(\cdot,\cdot)$ and $L(\cdot)$ provided in \eqref{alelas}, and 
where $\sigma_n(\bu)$ (resp. $v_n$) is the normal component of the Cauchy stress on the boundary (resp. the normal component of the virtual displacement).
The Signorini-type contact conditions (Karush-Kuhn-Tucker conditions) are expressed as (see, e.g., \cite{kikuchi-oden-88})
\begin{subequations}
\begin{align}
u_n - g & \leq 0 \qquad \textrm{on $\Gamma$} \label{contact1}, \\
\sigma_n(\bu) & \leq 0  \qquad \textrm{on $\Gamma$} \label{contact2}, \\
\sigma_n(\bu) \left( u_n - g \right) & = 0 \qquad \textrm{on $\Gamma$} \label{contact3}.
\end{align}
\end{subequations}
\Eref{contact1} is the non-penetration condition, whereas
\Eref{contact2} means the contact is non-adhesive,
and \Eref{contact3} is the complementarity condition.
%(the contact pressure can be non-zero if and only if adhesion occurs).

With the choice
\[
 \bB (\bu) = u_n, \quad \bs(\bu) = \sigma_n(\bu), \quad \barB = g, \quad S = \mathbb{R}^-,
\]
formulation \eqref{nitschefem} reads
\begin{equation} \label{nitschesigno}
\textrm{Find } \bu^h \in \bV^h \, : \, 
\An ( \bu^h , \bv^h ) 
 + \int_\Gamma \gamma_h^{-1} [ \PnUn(\bu^h) + \gamma_h g ]_{\mathbb{R}^{-}} \, \Pn ( \bv^h ) \, ds 
 = L(\bv^h), \quad \forall \, \bv^h \in \bV^h,
\end{equation}
which is exactly the formulation presented in \cite{chouly2015symmetric,chouly-overview-2017}.

\begin{remark}
\label{rem-contact-pfree}
The contact formulation is not parameter-free,
% for $\theta \neq -1$
%the value of $\gamma$ should be sufficiently large to ensure well-posedness and convergence \cite{chouly2015symmetric,fabre2014fictitious}.
%In our experience, $\gamma = E/h$ is a suitable value.
%For the skew-symmetric variant $\theta=-1$, the condition $\gamma > 0$ suffices to ensure the same theoretical properties, which is assessed numerically, see, e.g., \cite{chouly2015symmetric}. % shows that regarding the convergence rate in relative $H^1$ norm and relative $L^2$ norm,
%Note however that 
and the choice $\gamma_h = 0$ (``penalty-free'') is not permitted for contact.
% TODO_Can we explain why? Not really up to now (very interesting issue)
However see for instance \cite{burman-hansbo-larson-2016-pfree} for a first attempt at deriving a penalty-free method for Signorini contact.\end{remark}

%\begin{remark}
%When $\gamma$ is extremely large,
%the contact detection term $[ P_1(\bu^h) + \gamma g ]_{\mathbb{R}^{-}}$ is dominated by the displacement gap $u_n - g$,
%and the whole formulation behaves similarly to the penalty method.
%But one advantage of Nitsche's formulation is that when the current displacement gap is zero,
%the remaining contact normal stress term $\sigma_n$ is activated to enforce contact conditions.
%This means that penetration in the final configuration, after convergence, is not necessary, and is in practice of the same order as the discretization scheme for the bulk equation.
%In displacement-based penalty methods, no penetrations is equivalent to no reaction forces and hence no contact enforcement.
%However, in practice there could exist some slight penetrations when using Nitsche's method.
%\end{remark}

The adaptation of the above formulation \eqref{nitschesigno} %can be straightforwardly extended 
for biased (master-slave) contact between two elastic bodies (see, e.g., \cite{chouly-overview-2017}), reads:
%with the following modifications
\[
 \bB (\bu) = \llbracket u \rrbracket^{sl}_n , \quad \bs(\bu) = \sigma_n^{sl}(\bu), \quad \barB = g, \quad S = \mathbb{R}^-,
\]
where $\llbracket u \rrbracket^{sl}_n := (\bu^1(\bx^1) - \bu^2(\bx^2)) \cdot \boldsymbol{n}^1$ is the relative displacement written on the slave surface, and $\sigma_n^{sl} (\bu) (= \sigma_n(\bu^1))$ is the contact pressure on the slave surface. Also the bilinear form $a (\cdot\,,\cdot)$ and the linear form $L(\cdot)$ should incorporate the virtual work of both the master and slave elastic bodies, \emph{i.e.} they should be defined as in \eqref{def_aL_interface} (see, e.g., \cite{chouly-overview-2017} and references therein for more details).

\begin{remark}
The same methodology can be extended to other problems, involving for instance friction, dynamics and large deformations,
please refer to \cite{chouly-overview-2017} for a (non-restrictive) overview of possible extensions.
% Using the same methodology,
% extensions to other applications is feasible.
% For instance frictional contact particularly for Tresca and Coulomb friction \cite{chouly2014adaptation,mlika2017unbiased},
% and dynamic contact \cite{chouly2015dyna1,chouly2015dyna2}.
\end{remark}

\begin{remark}
\label{rem-regularity}
It is generally considered that the Sobolev regularity of contact problems is lower than $H^{\frac52} (\Omega)$ in two dimensions, due to weak singularities associated to transitions between binding and non-binding (see, e.g., \cite{moussaoui-khodja-92}). This means that, in general, FEM or IGA approximations of order higher than two do not improve the convergence rate in the energy norm, which remains limited to $\mathcal{O}(h^{\frac32})$ (see, e.g., \cite{auliac-13} in case of quadratic finite elements and, e.g., \cite{antolin2017textit} for IGA). Nevertheless, as discussed in \cite{antolin2017textit} the interest of the IGA approximation for contact is to obtain easily smooth gap functions $g$, which makes the numerical method more robust. As well, higher order approximations allow to recover smoother contact pressures. To perform better in terms of convergence, higher order approximations need to be combined with adaptive refinement, as in \cite{dimitri2014isogeometric,dorsek-melenk-2010,lee-oden-1994}.
\end{remark}

\subsubsection{Unbiased frictionless contact}
\label{subsub:unbiased}

Consider now the same situation as in the previous section \ref{subsub:biased} and depicted \fref{contact_problem}, but this time with $\Omega^1$ and $\Omega^2$ that represent both two elastic bodies in frictionless contact. Assume also, for simplicity, that there is no initial gap ($g=0$). We proceed first the same way as in Section \ref{sub:interface} and obtain the Green formula
\eqref{greeninterface}. We apply the frictionless condition, use the definition of $\langle \bsig(\bu) \rangle$, and separate the contributions on the two sides $\Gamma^1$ and $\Gamma^2$ of the interface, so the Green formula \eqref{greeninterface} can be re-written equivalently as%(considering frictionless contact) as
\begin{equation} \label{greenunbiased}
  a (\bu\,, \bv) 
  - \frac12 \int_{\Gamma^1} \sigma_n^1(\bu^1) \llbracket v \rrbracket^1_n \,ds
  - \frac12 \int_{\Gamma^2} \sigma_n^2(\bu^2) \llbracket v \rrbracket^2_n \,ds
  = L(\bv),
\end{equation}
with the notations $\llbracket v \rrbracket^1_n := (\bv^1 - \bv^2) \cdot \bn^1$ and $\llbracket v \rrbracket^2_n := (\bv^2 - \bv^1) \cdot \bn^2$.
We now apply frictionless contact conditions on the %each contact surface $\Gamma^m$, $m=1,2$, with
product set $\Gamma_1 \times \Gamma_2$:
\[
 \bB (\bu) 
 = ( \llbracket u \rrbracket_n^1 , \llbracket u \rrbracket_n^2 ), \quad \bs (\bu) = ( \sigma_n(\bu^1) , \sigma_n(\bu^2) ), \quad \barB = (0,0), \quad S = \RR^- \times \RR^-.
\]
The following unbiased formulation for contact is obtained
\begin{equation}
\textrm{Find } \bu^h \in \bV^h \, : \, 
\An(\bu^h,\bv^h) 
 + \frac{1}{2} \sum_{m=1}^2 \int_{\Gamma^m} \left( \gamma_h^{m} \right)^{-1} [ \PnUn^m(\bu^h)  ]_{\mathbb{R}^{-}} \, \Pn^m ( \bv^h ) \, ds
 = L(\bv^h), \quad \forall \, \bv^h \in \bV^h,
\end{equation}
where
\[
 \An ( \bu^h  , \bv^h ) :=  a(\bu^h , \bv^h) - \frac{\theta}{2}  
 \sum_{m=1}^2 \int_{\Gamma^m} \left ( \gamma_h^{m} \right)^{-1} \sigma_n(\bu^{h,m}) \sigma_n(\bv^{h,m}) \, ds,
\]
with $a(\cdot,\cdot)$ defined as in \eqref{def_aL_interface},
and, for $m=1,2$,
\[
 \Pn^m ( \bv^h ) := \theta \sigma_n( \bv^{h,m} ) - \gamma_h^{m} \llbracket v^h \rrbracket_n^m.
\]
This is a special case of \cite{rabii-numath-accepted} (see \cite{mlika2017unbiased,rabii-numath-accepted} for the detailed derivation in a more general setting). Remark that this formulation does not indeed differentiate between a master and a slave surface. %, and is particularly suitable for multi-body contact and self-contact .

\section{Numerical studies}\label{section_example}

To study the performance of the proposed skew-symmetric Nitsche's method,
we present some numerical tests.
We consider the IGA setting described in Section \ref{section_iga} with equal order of approximation in all directions, and carry out tests for different orders.
To avoid additional errors due to numerical integration, unless otherwise specified,
elements with $C^{p-1}$ continuity are adopted for order $p$,
and $p+1$ Gauss quadrature points are used for each element (and the same applies for other directions).
All the methods are implemented within the open source C\texttt{++} IGA library \emph{Gismo} \footnote{https://ricamsvn.ricam.oeaw.ac.at/trac/gismo/wiki/WikiStart} \cite{jlmmz2014}.

We recall that the skew-symmetric Nitsche's method corresponds to the Nitsche parameter $\theta=-1$.
For most of the numerical tests, we compare the performance of this method to the symmetric variant, that corresponds to $\theta=1$, and that will be denominated standard Nitsche's method. 
In this case the stabilization parameter is determined as $\gamma_0 = 2 \lambda^{h,\textrm{MAX}}$ where $\lambda^{h,\textrm{MAX}}$ is obtained from \eqref{eigenvalue}.
%In this case a suitable value of the stabilization parameter $\gamma_0$ needs to be chosen, as discussed in \ref{sub:Schur}. In the first part of the numerical tests, for linear boundary and interface conditions (Section \ref{sub:numlin} and Section \ref{sub:numinter}), we consider the penalty-free skew-symmetric Nitsche method ($\theta=-1$,$\gamma_0=0$) and for standard Nitsche we choose $\gamma_0 = 2 \lambda^{h,\textrm{MAX}}$ where $\lambda^{h,\textrm{MAX}}$ is obtained from \eqref{eigenvalue}. For the last part of the numerical tests, which deal with contact (Section \ref{sub:numcontact}), we can not consider a penalty-free Nitsche method anymore, and we prefer to test the robustness of the skew-symmetric variant in comparison to the symmetric standard one when the parameter $\gamma_0$ is varied within a determined range of values.

In order to evaluate the performances numerically,
the relative errors on the displacement field $\bu$ within the domain $\Omega$ is computed, in the $L^2$-norm, denoted by 
$\|\cdot\|_{L^2(\Omega)}$, and in the energy norm, denoted by 
$\|\cdot\|_{E(\Omega)} (= \sqrt{\int_{\Omega} \boldsymbol{\sigma} (\cdot) : \boldsymbol{\epsilon} (\cdot) \, d\bx})$.

\subsection{Linear boundary conditions}
\label{sub:numlin}

\subsubsection{Dirichlet boundary conditions patch test}

\begin{figure}[htbp]
\centering
\subfigure[Rectangular patch]{\def\svgwidth{0.48\columnwidth}
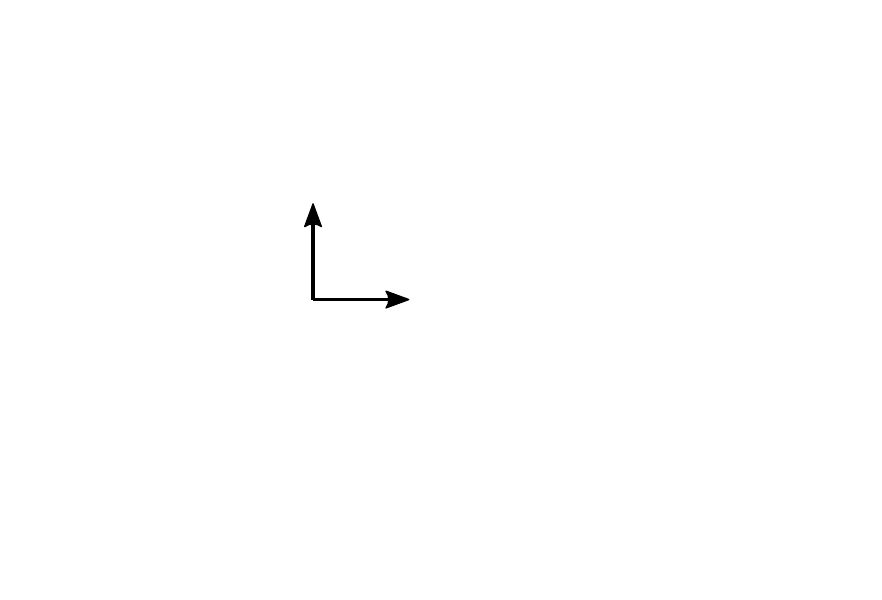}
\subfigure[Circular patch]{\def\svgwidth{0.49\columnwidth}
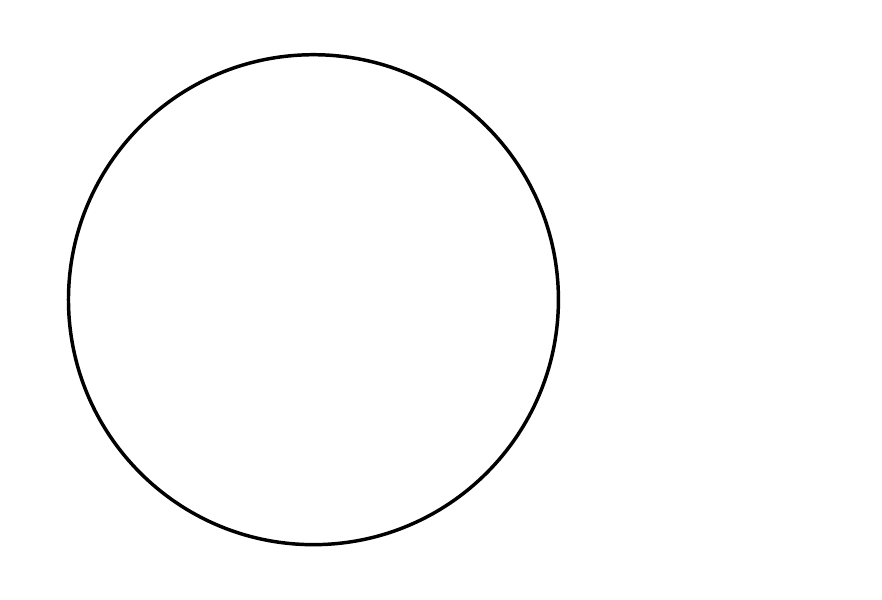}
\caption{Patch test problems.}
\label{patch_test_desc}
\end{figure}

In this section we focus on the setting presented in \ref{subsub:Dirichlet} and illustrate the skew-symmetric Nitsche's formulation is able to pass the rectangular patch tests and has optimal convergence rate in energy norm for circular shaped patch tests.
In \fref{patch_test_desc} we consider a linear elastic media in a square (resp. circular) domain $\Omega$ of length $L=20$ (resp. radius $R=10$), with a Young's modulus $E=1000$ and a Poisson's ratio $\nu=0.25$. 

The patch test is traditionally adopted to verify the consistency of a newly proposed element.
However, in the following, another type of patch test,
i.e. the B-type patch test \cite{Zienkiewicz2005329},
is used to test the effectiveness of the skew-symmetric Nitsche's method in imposing Dirichlet boundary conditions.
Firstly we set the exact solution as ${\bu}^{\textrm{ref}}$,
and impose the value $\bar{\bu} = {\bu}^{\textrm{ref}} |_{\Gamma_D}$ on the whole boundary $\Gamma_D := \partial \Omega$ in \eqref{strongelast}. No external force is imposed: $\boldsymbol{b} = \boldsymbol{0}$.
Finally problem \eqref{nitschefemlinear} is solved and the corresponding solution $\bu^h$ is compared against the exact solution ${\bu}^{\textrm{ref}}$.
In order to fulfill the equilibrium strain condition, an exact solution of displacement field up to fourth order is set up as in \cite{chen2010family}:
\[
\small
\begin{cases}
u_x(x,y) = \frac{1}{4}+x+3y-2x^2-4xy+\frac{5}{2}y^2-2x^3+x^2y-4xy^2-\frac{1}{3}y^3-\frac{7}{32}x^4-\frac{19}{24}x^3y+x^2y^2+xy^3-\frac{11}{96}y^4, \\
u_y(x,y) = 1+\frac{1}{2}x+2y-\frac{2}{3}x^2+\frac{17}{5}xy+\frac{3}{2}y^2+\frac{1}{3}x^3+12x^2y-xy^2-\frac{2}{3}y^3-\frac{11}{96}x^4+x^3y+x^2y^2-\frac{19}{24}xy^3-\frac{7}{32}y^4, \\
  \end{cases}
\]
where $u_x$ (resp. $u_y$) is the $x$-component (resp $y$-component) of 
${\bu}^{\textrm{ref}}$,
truncated to the appropriate order.
For example, if the patch test of order one is performed, then the exact solution is truncated as
\[
\begin{cases}
u_x(x,y) = \frac{1}{4}+x+3y, \\
u_y(x,y) = 1+\frac{1}{2}x+2y. \\
  \end{cases}
\]

\begin{table}[htbp]
\centering
\caption{The skew-symmetric Nitsche's formulation passes (Y) the rectangular patch tests.}
\label{patch_test_tab}

\medskip

\begin{tabular}{l cccc}
 \hline
   Patch test order & 1 & 2 & 3 & 4\\ \hline \noalign{\smallskip}
  IGA $p=q=2$ &  Y & Y & N & N \\ 
  IGA $p=q=3$ &  Y & Y & Y & N \\ 
  IGA $p=q=4$ &  Y & Y & Y & Y \\ 
    \hline
\end{tabular}
\end{table}

\begin{figure}[htbp]
\centering
\definecolor{mycolor1}{rgb}{1.00000,0.00000,1.00000}%
\begin{tikzpicture}

\begin{axis}[%
width=1.2\figurewidth,
height=1.2\figureheight,
at={(0\figurewidth,0\figureheight)},
scale only axis,
separate axis lines,
every outer x axis line/.append style={black},
every x tick label/.append style={font=\color{black}},
every x tick/.append style={black},
xmode=log,
xmin=1,
xmax=200,
xminorticks=true,
xlabel={Number of elements per side},
every outer y axis line/.append style={black},
every y tick label/.append style={font=\color{black}},
every y tick/.append style={black},
ymode=log,
ymin=1e-09,
ymax=1,
yminorticks=true,
ylabel={$||\bu^h - \bu^{\text{ref}}||_{E(\Omega)} / ||\bu^{\text{ref}}||_{E(\Omega)}$},
axis background/.style={fill=white},
xmajorgrids,
xminorgrids,
ymajorgrids,
yminorgrids,
legend style={at={(0.03,-0.22)}, anchor=south west, legend cell align=left, align=left, fill=white}
]

\addplot [color=black, line width=1.5pt, mark=square, mark options={solid, black}]
  table[row sep=crcr]{%
2	1.45045E-01\\
4   2.90857E-02\\
8   5.86244E-03\\
16  1.28367E-03\\
32  2.94155E-04\\
64  6.96464E-05\\
128 1.68593E-05\\
};
\addlegendentry{Nitsche $p=q=2$ \\ 2nd order patch test}

\addplot [color=red, line width=1.5pt, mark=o, mark options={solid, red}]
  table[row sep=crcr]{%
2	1.78698E-01\\
4   1.35493E-02\\
8   1.14534E-03\\
16  1.31567E-04\\
32  1.60735E-05\\
64  1.99559E-06\\
128  2.48895E-07\\
};
\addlegendentry{Nitsche $p=q=3$ \\ 3rd order patch test}

\addplot [color=blue, line width=1.5pt, mark=triangle, mark options={solid, blue}]
  table[row sep=crcr]{%
2	1.12701E-01\\
4   1.41507E-02\\
8   4.04110E-04\\
16  2.10613E-05\\
32  1.27057E-06\\
64  7.92761E-08\\
128 4.97321E-09\\
};
\addlegendentry{Nitsche $p=q=4$ \\ 4th order patch test}

\end{axis}

\node at (7.2,4.70) {2.08}; % black
\node [color=red] at (7.2,3.45) {3.01}; % red
\node [color=blue] at (7.2,2.40) {4.00}; % red

\end{tikzpicture}%
\caption{Circular patch test: relative errors of displacement field in energy norm. Skew-symmetric Nitsche's method is used.}
\label{Fig_patch_test_energy}
\end{figure}
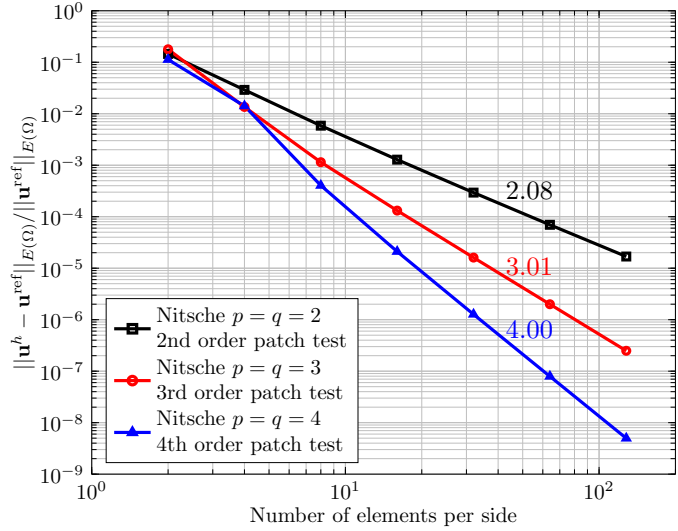

The results of the rectangular patch tests are presented in Table \ref{patch_test_tab},
showing that the skew-symmetric Nitsche's method is able to pass the appropriate patch tests of order up to $p$ for $p=2,3,4$.
The circular patch tests cannot be passed exactly,
because Nitsche's method imposes the boundary constraints weakly. % for curved boundaries.
However according to \fref{Fig_patch_test_energy} the error in the energy norm is reduced with an optimal convergence rate of order $p$,
as predicted theoretically for IGA and a conformal setting (i.e. strong imposition of Dirichlet boundary conditions) 
\cite{bazilevs2006isogeometric}.
This is also in agreement with the observed behavior of skew-symmetric Nitsche's formulation with FEM: see for instance \cite{boiveau2016penalty} where the same rates are obtained numerically for quadratic finite elements.
As regarding the convergence in $L^2(\Omega)$-norm, a sub-optimality of order $\mathcal{O}(h^{\frac12})$ is predicted by the theory, due to the lack of adjoint-consistency of skew-symmetric Nitsche's method, but this behavior seems difficult to observe in practical situations (see \cite{boiveau2016penalty}).

\subsubsection{Symmetry conditions for Kirchhoff plates}

\begin{figure}[htbp]
\centering
\def\svgwidth{0.8\columnwidth}
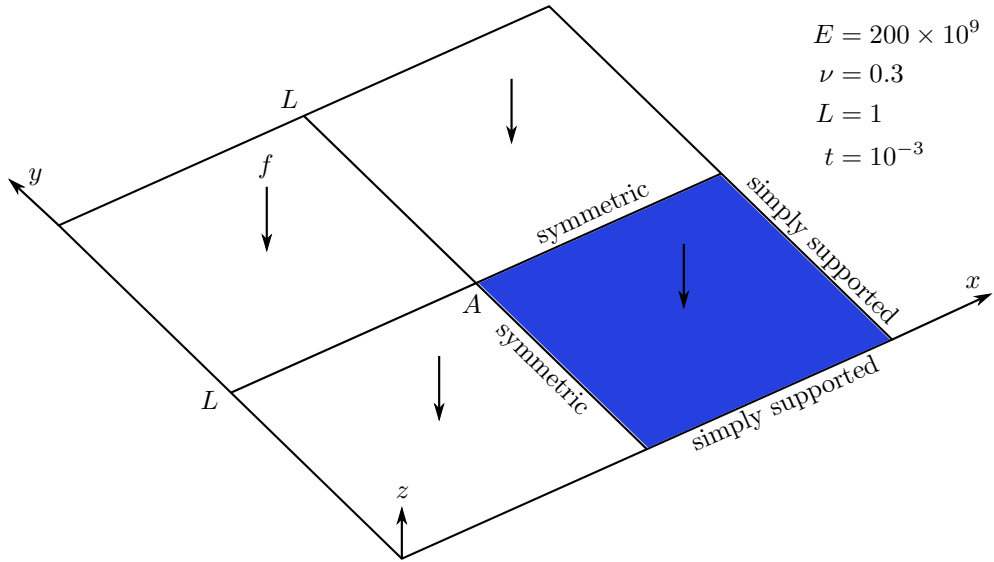
\caption{Square thin plate under distributed transverse load, only $1/4$ of the plate (blue area) is modeled.}
\label{Fig_spsp}
\end{figure}

In this section we illustrate the effectiveness of Nitsche's formulation to handle rotational boundary conditions for Kirchhoff plates,
as described in \ref{subsub:Kirchhoff}.
\fref{Fig_spsp} describes a simply supported thin square plate of thickness $t$, made of an isotropic elastic material, and subjected to a distributed transverse load $f$. \fref{Fig_spsp} provides also the values of the model parameters.
Due to the symmetry of this problem, only one quarter of the geometry is modeled,
where two of the model boundaries are simply supported,
and the other two require symmetric constraints:
\[
\bar{\theta}_t=0 \qquad \textrm{on $y=\frac{L}{2}$ and $x=\frac{L}{2}$}.
\]
%where the subscript $t$ denotes the tangential direction of the boundary.
For a sinusoidally distributed transverse load
\[
 f (x,y) = -10 \sin (\pi x) \sin (\pi y),
\]
the analytical solution of the deflection is given by \cite{reddy2006theory}
\[
u^\text{ref} (x,y)
= \frac{-10}{4 \pi^4 D} \sin(\pi x) \sin (\pi y),
\]
in which $D=\tfrac{Et^3}{12(1-\nu^2)}$ is the flexural rigidity.

\begin{figure}[htbp]
\centering
\def\svgwidth{0.7\columnwidth}
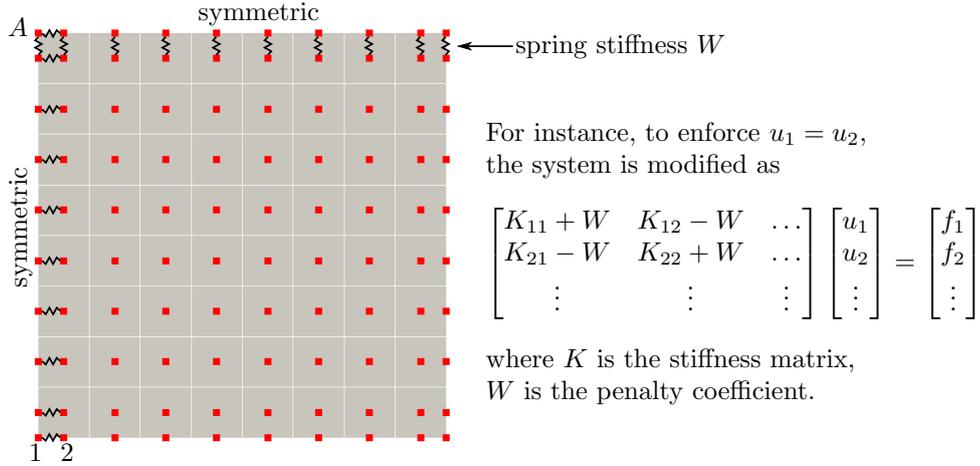
\caption{The `second row' strategy with penalty method to impose the symmetry boundary conditions. The four control points near point $A$ are penalized twice.}
\label{fig_kirchhoff_rotation_bc_penalty}
\end{figure}

\begin{figure}[htbp]
	\centering
	\subfigure[The `second row' strategy using penalty method]{% This file was created by matlab2tikz.
\definecolor{mycolor1}{rgb}{1.00000,0.00000,1.00000}%
\begin{tikzpicture}

\begin{axis}[%
width=1.2\figurewidth,
height=1.2\figureheight,
at={(0\figurewidth,0\figureheight)},
scale only axis,
separate axis lines,
every outer x axis line/.append style={black},
every x tick label/.append style={font=\color{black}},
every x tick/.append style={black},
xmode=log,
xmin=1,
xmax=200,
xminorticks=true,
xlabel={Number of elements per side},
every outer y axis line/.append style={black},
every y tick label/.append style={font=\color{black}},
every y tick/.append style={black},
ymode=log,
ymin=1e-05,
ymax=1,
yminorticks=true,
ylabel={$||\bu^h - \bu^{\text{ref}}||_{E(\Omega)} / ||\bu^{\text{ref}}||_{E(\Omega)}$},
axis background/.style={fill=white},
xmajorgrids,
xminorgrids,
ymajorgrids,
yminorgrids,
legend style={at={(0.525,0.665)}, anchor=south west, legend cell align=left, align=left, fill=white}
]

\addplot [color=black, line width=1.5pt, mark=square, mark options={solid, black}, mark size=3pt]
  table[row sep=crcr]{%
2	1.62522E-01\\
4   8.04213E-02\\
8   4.01120E-02\\
16  2.00439E-02\\
32  1.00204E-02\\
64  5.01002E-03\\
128 2.50499E-03\\
};
\addlegendentry{penalty $=1\times 10^{10}$, $p=q=2$}

\addplot [color=red, line width=1.5pt, mark=o, mark options={solid, red}, mark size=3pt]
  table[row sep=crcr]{%
2	1.73356E-02\\
4   4.14228E-03\\
8   1.02120E-03\\
16  2.54339E-04\\
32  6.36294E-05\\
64  1.75669E-05\\
128 1.56427E-05\\
};
\addlegendentry{penalty $=1\times 10^{10}$, $p=q=3$}

\addplot [color=blue, line width=1.5pt, mark=triangle, mark options={solid, blue}, mark size=3pt]
  table[row sep=crcr]{%
2	1.62522E-01\\
4   8.04213E-02\\
8   4.01120E-02\\
16  2.00439E-02\\
32  1.00209E-02\\
64  5.03290E-03\\
128 2.90526E-03\\
};
\addlegendentry{penalty $=1\times 10^{13}$, $p=q=2$}

\addplot [color=mycolor1, line width=1.5pt, mark=triangle, mark options={solid, rotate=180, mycolor1}, mark size=3pt]
  table[row sep=crcr]{%
2	1.73356E-02\\
4   4.14231E-03\\
8   1.02130E-03\\
16  2.63747E-04\\
32  4.48019E-04\\
64  1.08808E-03\\
128  2.43514E-03\\
};
\addlegendentry{penalty $=1\times 10^{13}$, $p=q=3$}

\draw (10,0.01) -- (20,0.005);
\draw (2,0.005) -- (4,0.00125);

\end{axis}

\node at (4.08,4.62) {1};
\node at (1.09,4.12) {2};

\node at (7.2,4.64) {1.00}; % black
\node [color=red] at (4.68,3.04) {2.01}; % red

\end{tikzpicture}%
}
    \subfigure[Nitsche's method]{% This file was created by matlab2tikz.
\definecolor{mycolor1}{rgb}{1.00000,0.00000,1.00000}%
\begin{tikzpicture}

\begin{axis}[%
width=1.2\figurewidth,
height=1.2\figureheight,
at={(0\figurewidth,0\figureheight)},
scale only axis,
separate axis lines,
every outer x axis line/.append style={black},
every x tick label/.append style={font=\color{black}},
every x tick/.append style={black},
xmode=log,
xmin=1,
xmax=200,
xminorticks=true,
xlabel={Number of elements per side},
every outer y axis line/.append style={black},
every y tick label/.append style={font=\color{black}},
every y tick/.append style={black},
ymode=log,
ymin=1e-06,
ymax=1,
yminorticks=true,
ylabel={$||\bu^h - \bu^{\text{ref}}||_{E(\Omega)} / ||\bu^{\text{ref}}||_{E(\Omega)}$},
axis background/.style={fill=white},
xmajorgrids,
xminorgrids,
ymajorgrids,
yminorgrids,
legend style={at={(0.42,0.692)}, anchor=south west, legend cell align=left, align=left, fill=white}
]

\addplot [color=black, line width=1.5pt, mark=square, mark options={solid, black}, mark size=3pt]
  table[row sep=crcr]{%
2	1.68036E-01\\
4   8.09234E-02\\
8   4.01491E-02\\
16  2.00464E-02\\
32  1.00206E-02\\
64  5.01003E-03\\
128 2.50499E-03\\
};
\addlegendentry{skew-symmetric Nitsche, $p=q=2$}

\addplot [color=red, line width=1.5pt, mark=o, mark options={solid, red}, mark size=3pt]
  table[row sep=crcr]{%
2	2.61545E-02\\
4   5.36347E-03\\
8   1.18681E-03\\
16  2.76044E-04\\
32  6.63051E-05\\
64  1.62288E-05\\
128 4.01313E-06\\
};
\addlegendentry{skew-symmetric Nitsche, $p=q=3$}

\addplot [color=blue, line width=1.5pt, mark=triangle, mark options={solid, blue}, mark size=3pt]
  table[row sep=crcr]{%
2	1.62520E-01\\
4   8.04212E-02\\
8   4.01120E-02\\
16  2.00439E-02\\
32  1.00204E-02\\
64  5.01002E-03\\
128 2.50499E-03\\
};
\addlegendentry{standard Nitsche, $p=q=2$}

\addplot [color=mycolor1, line width=1.5pt, mark=triangle, mark options={solid, rotate=180, mycolor1}, mark size=3pt]
  table[row sep=crcr]{%
2	1.73355E-02\\
4   4.14227E-03\\
8   1.02120E-03\\
16  2.54329E-04\\
32  6.35193E-05\\
64  1.58883E-05\\
128 1.00188E-05\\
};
\addlegendentry{standard Nitsche, $p=q=3$}

\draw (10,0.01) -- (20,0.005);
\draw (2,0.005) -- (4,0.00125);

\end{axis}

\node at (4.08,5.1) {1};
\node at (1.09,4.72) {2};

\node at (7.2,5.14) {1.00}; % black
\node [color=red] at (6.0,2.94) {2.06}; % red

\end{tikzpicture}%
}
	\caption{Kirchhoff plate: relative errors of deflection field in energy norm. Symmetric rotational boundary conditions are imposed by the `second row' strategy using penalty method (a), and Nitsche's method (b).}
	\label{kirchoff_result}
\end{figure}
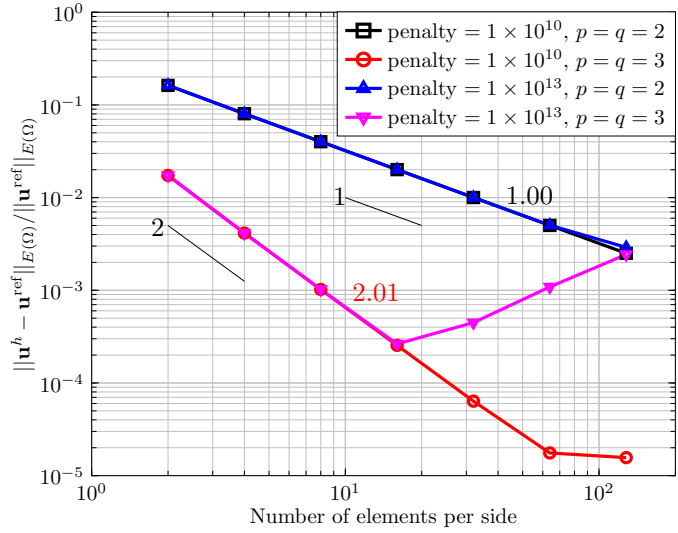
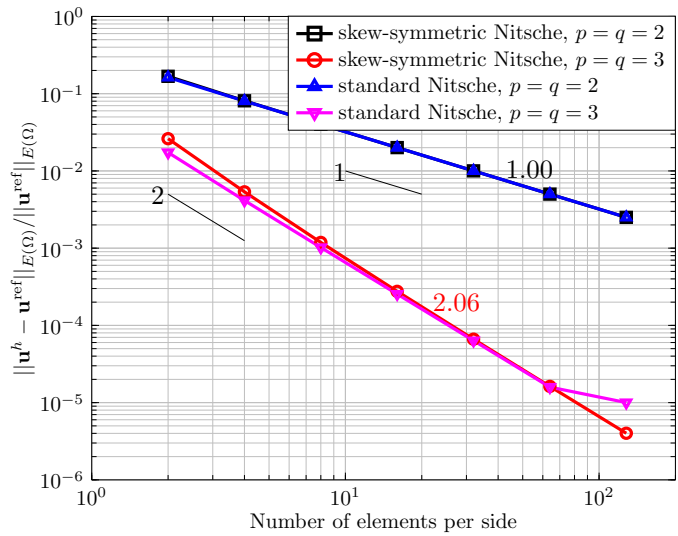

We firstly implement the `second row' strategy \cite{kiendl2009isogeometric} by penalty method for comparison.
As illustrated in \fref{fig_kirchhoff_rotation_bc_penalty},
the idea of the `second row' strategy consists in enforcing the displacements of the control points along the symmetric boundary and the neighboring row to be equal.
This is achieved by adding a penalty coefficient $W$ in the stiffness matrix \footnote{http://www.colorado.edu/engineering/CAS/courses.d/IFEM.d/IFEM.Ch09.d/IFEM.Ch09.pdf}.
However in this way,
the convergence results depend significantly on the penalty parameters,
and the constraints of the four control points are penalized twice (see \fref{fig_kirchhoff_rotation_bc_penalty}) near the corner point $A$,
which makes the corner deflection even more sensitive to the penalty coefficient.
The results using different penalty parameters are shown in \fref{kirchoff_result} (a).
It is concluded that a suitable value of the penalty parameter $W$ should be chosen carefully for different meshes and orders.

The results obtained by Nitsche's method are shown in \fref{kirchoff_result} (b).
The stabilization parameter $\gamma_0$ for the standard (symmetric) Nitsche's formulation is acquired by solving the generalized eigenvalue problem \eqref{eigenvalue}.
Remember that the Kirchhoff problem results in a fourth order system with respect to the deflection $u$,
the strain $\boldsymbol{\epsilon}$ consists of second order derivatives of $u$,
thus the "energy norm" in this situation is equivalent to the $H^2$ semi-norm on the deflection,
and the optimal convergence rate is expected to be $p-1$ in the energy norm for approximation order $p$ \cite{embar2010imposing}.
As indicated in \fref{kirchoff_result} (b),
for relative errors in energy norm the skew-symmetric Nitsche's method and the standard one are similar,
and both standard and skew-symmetric Nitsche's formulations converge optimally
with the expected orders associated to Kirchhoff plate theory.
%Note that for $p=q=3$ a reduction of the slope is observed for the most refined solutions.
%probably due to a limitation from machine precision.

\subsection{Linear interface conditions and patch coupling}
\label{sub:numinter}

\subsubsection{Patch coupling effects: statics}\label{pach_coupling_statics}

\begin{figure}[htbp]
\centering
\def\svgwidth{0.9\columnwidth}
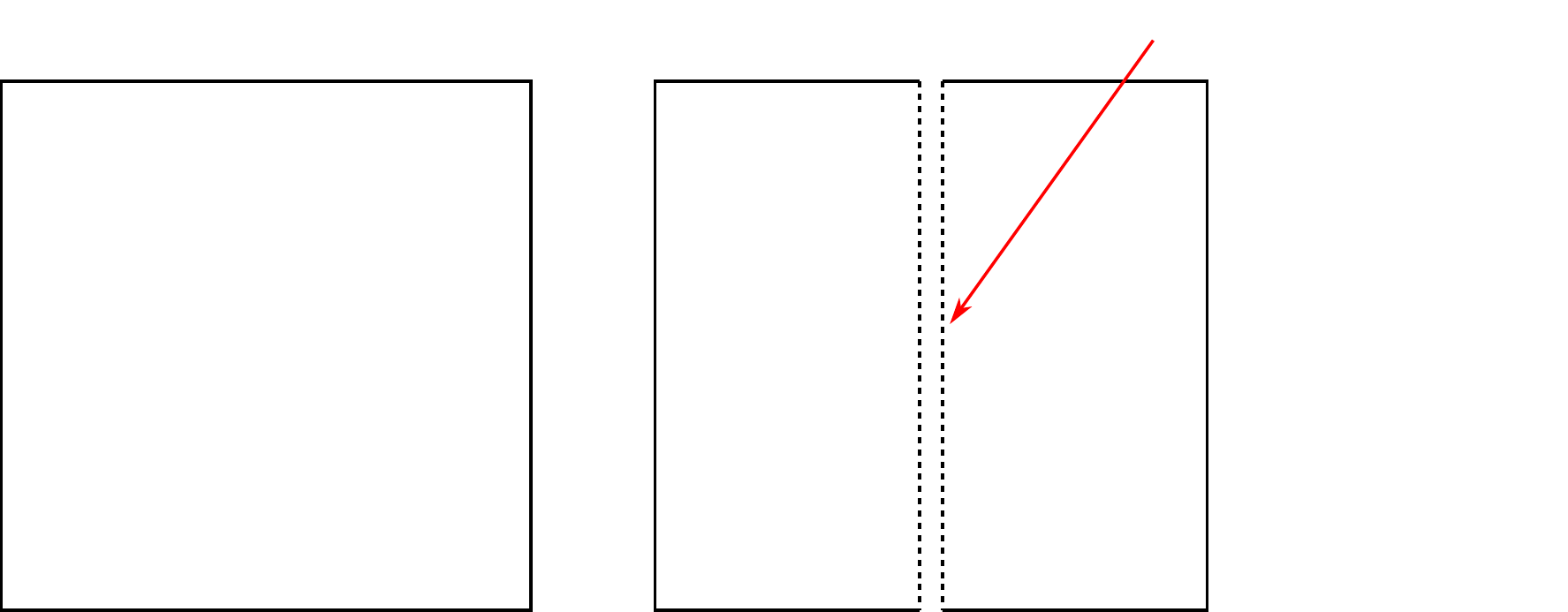
\caption{Plate model with four simply supported edges. On the left, one patch model is adopted as the control group. On the right, the plate is artificially broken into two conforming patches and the interface is coupled by Nitsche's method.}
\label{fig_sq_plate_couple}
\end{figure}

In this section,
we study whether additional effects are introduced into the accuracy, convergence performance and condition numbers when a patch coupling in statics is performed by Nitsche's method.
The problem setup is shown in \fref{fig_sq_plate_couple}:
in the left figure we present a plate model with thickness $t=10^{-1}$, Young's modulus $E=200\times 10^9$ and Poisson's ratio $\nu=0.3$.
The plate is subjected to uniform pressure $f$ with four edges being simply supported.
In the right figure the domain $\Omega$ is artificially broken into two identical patches $\Omega_1$ and $\Omega_2$. % with the setups being unchanged,
It corresponds then to the setting described in \ref{sub:interface} with the value of the parameters provided in \fref{fig_sq_plate_couple},
%The only difThe equations that are solved by 
%The corresponding equations are approximated 
and with an approximation using degenerate Reissner-Mindlin elements, % in linear elasticity,
%for %these elements 
in which only the mid-surface of the plate has to be modeled \cite{adam2015improved}.
The deflection field approximated with (conforming) IGA using one patch of 1,024 elements of order $p(=q)=5$ is adopted as the reference.
%Its value at the parameter center $(\xi=0.5,\eta=0.5)$ is denoted by $w_A$,
%and we found $w_A=-2.52083\times 10^{-5}$.

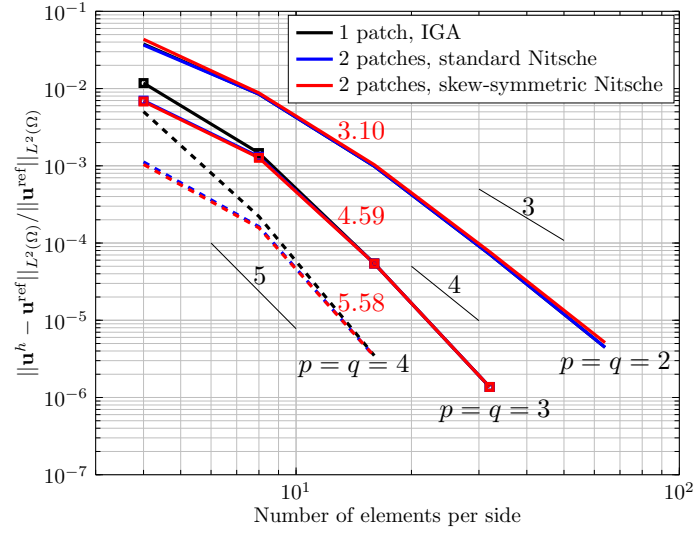
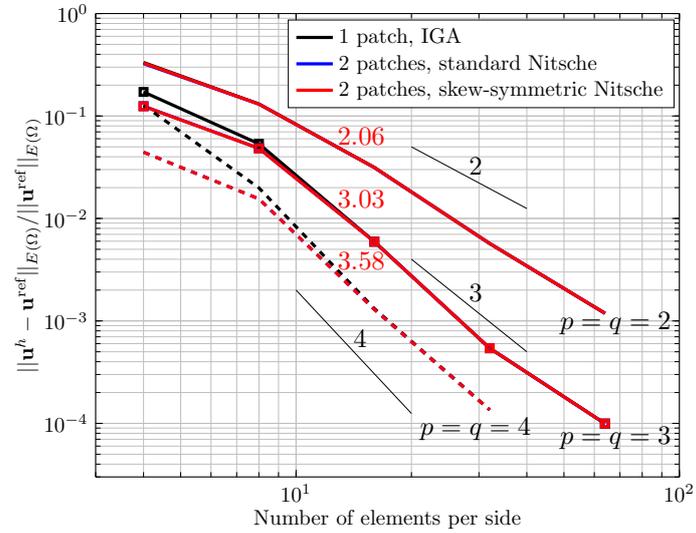
\begin{figure}[htbp]
\centering
	\subfigure[$L^2$ norm]{% This file was created by matlab2tikz.
\begin{tikzpicture}

\begin{axis}[%
width=1.2\figurewidth,
height=1.2\figureheight,
at={(0\figurewidth,0\figureheight)},
scale only axis,
separate axis lines,
every outer x axis line/.append style={black},
every x tick label/.append style={font=\color{black}},
every x tick/.append style={black},
xmode=log,
xmin=3,
xmax=100,
xminorticks=true,
xlabel={Number of elements per side},
every outer y axis line/.append style={black},
every y tick label/.append style={font=\color{black}},
every y tick/.append style={black},
ymode=log,
ymin=1e-07,
ymax=0.1,
yminorticks=true,
ylabel={$||\bu^h - \bu^{\text{ref}}||_{L^2(\Omega)} / ||\bu^{\text{ref}}||_{L^2(\Omega)}$},
axis background/.style={fill=white},
xmajorgrids,
xminorgrids,
ymajorgrids,
yminorgrids,
legend style={at={(0.49,0.71)}, anchor=south west, legend cell align=left, align=left, draw=black}
]

\addplot [color=black, line width=1.5pt]
  table[row sep=crcr]{%
4	3.77247E-02\\
8	8.46946E-03\\
16	9.87279E-04\\
32	7.15617E-05\\
64  4.45933E-06\\
};
\addlegendentry{1 patch, IGA}

\addplot [color=blue, line width=1.5pt]
  table[row sep=crcr]{%
4	3.66179E-02\\
8	8.46386E-03\\
16	9.87232E-04\\
32	7.15604E-05\\
64  4.45927E-06\\
};
\addlegendentry{2 patches, standard Nitsche}

\addplot [color=red, line width=1.5pt]
  table[row sep=crcr]{%
4	4.34769E-02\\
8	8.77962E-03\\
16	1.03266E-03\\
32	7.70307E-05\\
64  5.13640E-06\\
};
\addlegendentry{2 patches, skew-symmetric Nitsche}

%p=3
\addplot [color=black, line width=1.5pt, mark=square, mark options={solid, black}, forget plot]
  table[row sep=crcr]{%
4	1.17879E-02\\
8	1.46328E-03\\
16	5.44206E-05\\
32	1.36886E-06\\
%64  6.09465E-07\\
};
\addplot [color=blue, line width=1.5pt, mark=square, mark options={solid, blue}, forget plot]
  table[row sep=crcr]{%
4	6.98290E-03\\
8	1.31031E-03\\
16	5.43115E-05\\
32	1.36882E-06\\
%64  6.09465E-07\\
};
\addplot [color=red, line width=1.5pt, mark=square, mark options={solid, red}, forget plot]
  table[row sep=crcr]{%
4	6.80577E-03\\
8	1.26882E-03\\
16	5.42033E-05\\
32	1.36866E-06\\
%64  6.09465E-07\\
};

%p=4
\addplot [color=black, dashed, line width=1.5pt, forget plot]
  table[row sep=crcr]{%
4	5.02609E-03\\
8	2.20672E-04\\
16	3.48869E-06\\
%32	6.05882E-07\\
};
\addplot [color=blue, dashed, line width=1.5pt, forget plot]
  table[row sep=crcr]{%
4	1.12235E-03\\
8	1.62813E-04\\
16	3.40454E-06\\
%32	6.05870E-07\\
};
\addplot [color=red, dashed, line width=1.5pt, forget plot]
  table[row sep=crcr]{%
4	1.04086E-03\\
8	1.57217E-04\\
16	3.42334E-06\\
%32	6.05855E-07\\
};

\draw (30,0.0005) -- (50,0.108e-3);
\draw (20,0.00005) -- (30,0.9876543206e-5);
\draw (6,1e-4) -- (10,0.7776e-5);

\end{axis}

\node at (7.12,4.45) {3};
\node at (5.9,3.15) {4};
\node at (2.7,3.3) {5};

\node at (8.55,1.85) {$p=q=2$};
\node at (6.57,1.05) {$p=q=3$};
\node at (4.25,1.8) {$p=q=4$};

\node [color=red] at (4.37,5.70) {3.10}; %blue

\node [color=red] at (4.37,4.30) {4.59}; %blue

\node [color=red] at (4.37,2.86) {5.58}; %blue

\end{tikzpicture}%
}
    \subfigure[Energy norm]{% This file was created by matlab2tikz.
\begin{tikzpicture}

\begin{axis}[%
width=1.2\figurewidth,
height=1.2\figureheight,
at={(0\figurewidth,0\figureheight)},
scale only axis,
separate axis lines,
every outer x axis line/.append style={black},
every x tick label/.append style={font=\color{black}},
every x tick/.append style={black},
xmode=log,
xmin=3,
xmax=100,
xminorticks=true,
xlabel={Number of elements per side},
every outer y axis line/.append style={black},
every y tick label/.append style={font=\color{black}},
every y tick/.append style={black},
ymode=log,
ymin=3e-05,
ymax=1,
yminorticks=true,
ylabel={$||\bu^h - \bu^{\text{ref}}||_{E(\Omega)} / ||\bu^{\text{ref}}||_{E(\Omega)}$},
axis background/.style={fill=white},
xmajorgrids,
xminorgrids,
ymajorgrids,
yminorgrids,
legend style={at={(0.49,0.745)}, anchor=south west, legend cell align=left, align=left, draw=black}
]+

\addplot [color=black, line width=1.5pt]
  table[row sep=crcr]{%
4	3.33182E-01\\
8   1.30936E-01\\
16  3.14715E-02\\
32  5.69680E-03\\
64  1.18585E-03\\
};
\addlegendentry{1 patch, IGA}

\addplot [color=blue, line width=1.5pt]
  table[row sep=crcr]{%
4	3.24409E-01\\
8   1.30883E-01\\
16  3.14707E-02\\
32  5.69677E-03\\
64  1.18585E-03\\
};
\addlegendentry{2 patches, standard Nitsche}

\addplot [color=red, line width=1.5pt]
  table[row sep=crcr]{%
4   3.28501E-01\\
8   1.30730E-01\\
16  3.14618E-02\\
32  5.69632E-03\\
64  1.18583E-03\\
};
\addlegendentry{2 patches, skew-symmetric Nitsche}

%p=3
\addplot [color=black, line width=1.5pt, mark=square, mark options={solid, black}, forget plot]
  table[row sep=crcr]{%
4	1.71784E-01\\
8   5.35148E-02\\
16  5.94928E-03\\
32  5.41732E-04\\
64  9.94080E-05\\
};
\addplot [color=blue, line width=1.5pt, mark=square, mark options={solid, blue}, forget plot]
  table[row sep=crcr]{%
4	1.24839E-01\\
8   4.83287E-02\\
16  5.93189E-03\\
32  5.41692E-04\\
64  9.94074E-05\\
};
\addplot [color=red, line width=1.5pt, mark=square, mark options={solid, red}, forget plot]
  table[row sep=crcr]{%
4	1.24266E-01\\
8   4.83236E-02\\
16  5.93191E-03\\
32  5.41692E-04\\
64  9.94074E-05\\
};

%p=4
\addplot [color=black, dashed, line width=1.5pt, forget plot]
  table[row sep=crcr]{%
4	1.27498E-01\\
8   1.98843E-02\\
16  1.32157E-03\\
32  1.35655E-04\\
%64  9.37561E-05\\
};
\addplot [color=blue, dashed, line width=1.5pt, forget plot]
  table[row sep=crcr]{%
4	4.43602E-02\\
8   1.55320E-02\\
16  1.29654E-03\\
32  1.35518E-04\\
%64  9.37561E-05\\
};
\addplot [color=red, dashed, line width=1.5pt, forget plot]
  table[row sep=crcr]{%
4	4.43176E-02\\
8   1.55298E-02\\
16  1.29655E-03\\
32  1.35518E-04\\
%64  9.37561E-05\\
};

\draw (20,0.05) -- (40,0.0125);
\draw (20,0.004) -- (40,0.0005);
\draw (10,0.002) -- (20,0.000125);

\end{axis}

\node at (6.25,5.15) {2};
\node at (6.25,3.05) {3};
\node at (4.35,2.3) {4};

\node at (8.55,2.55) {$p=q=2$};
\node at (8.55,0.65) {$p=q=3$};
\node at (6.25,0.79) {$p=q=4$};

\node [color=red] at (4.37,5.65) {2.06}; %blue

\node [color=red] at (4.37,4.60) {3.03}; %blue

\node [color=red] at (4.37,3.58) {3.58}; %blue

\end{tikzpicture}%
}
    %In the energy norm only the strain component $\epsilon_{13}=\theta_y+dw/dx$ and corresponding stress $\sigma_{13}$ is used.
	\caption{Patch coupling of two plates: relative errors of displacement field in $L^2$ norm (a) and energy norm (b).}
	\label{couple_influence_converge}
\end{figure}

We compare symmetric and skew-symmetric variants of Nitsche's method,
and for the symmetric variant,
we still compute $\gamma_0$ by solving the generalized eigenvalue problem \eqref{eigenvalue}.
The convergence performance is plotted in \fref{couple_influence_converge}.
In this test the meshes of the left patch and the right patch are equal (for instance, for conforming patch the mesh is $8\times 8$, then for two patch coupling the meshes are $4\times 8$ and $4\times 8$ for the left patch and right patch respectively),
in order to evaluate the coupling influence on the approximation of the displacement.
By artificially breaking the patch into two patches and couple them by Nitsche's method,
the obtained errors is to some degrees different from the one patch case.
%For higher order elements ($p=q=3$ and $p=q=4$) the skew-symmetric formulation performs more accurately than the standard formulation.
Skew-symmetric and standard formulations perform very similarly, and
as we refine the mesh,  both %the standard and the skew-symmetric Nitsche's formulations 
converge with nearly optimal rates.

The condition numbers of the obtained stiffness matrix are given in \Tref{couple_influence_condition_number}.
From a general viewpoint,
h-refinement of the mesh,
which means that using more control points,
increases the corresponding condition number.
It is noticed that the condition number obtained from Nitsche's coupling is larger than one patch IGA.
This is inferred to be related to the coupling effects.
Specifically, the skew-symmetric Nitsche's formulation slightly increases the condition number compared to conforming IGA because there are more control points along the coupled interface,
moreover the condition numbers are almost independent of the mesh size $h$ and basis functions orders $p$ and $q$.
%The stabilization parameter for the standard Nitsche formulation and it is usually a large value,
%which is acquired by solving the generalized eigenvalue problem along the boundary %\cite{embar2010imposing,apostolatos2014nitsche}.
The standard Nitsche's formulation increases the condition number significantly,
because of the large value of the stabilization parameter $\gamma_0$.
%which is inferred to be attributed to the adopted stabilization parameter, which has a large value.

\begin{table}[htbp]
\centering
\caption{Patch coupling of two plates: condition numbers ($\times 10^{10}$) of the obtained stiffness matrix. The standard Nitsche's formulation increases the condition number significantly.} 
\label{couple_influence_condition_number}

\medskip

\begin{tabular}{c c cccc}
 \hline
 Method &  Number of elements per side & $p=q=2$ & $p=q=3$ & $p=q=4$ & $p=q=5$ \\ \hline \noalign{\smallskip}
Conforming &  $4$ &  0.926 & 0.884 & 0.858 & 0.850 \\ 
IGA &  $8$ &  0.958 & 0.979 & 1.103 & 1.238 \\
&  $16$ &  0.959 & 1.002 & 1.161 & 1.347 \\ 
 \hline
Standard&  $4$ &  8.942 & 15.033 & 21.711 & 28.011 \\ 
Nitsche&  $8$ &  11.963 & 18.104 & 23.693 & 29.033 \\
&  $16$ &  13.301 & 21.680 & 29.204 & 35.938 \\ 
 \hline
Skew-symmetric&  $4$ &  1.919 & 1.340 & 1.285 & 1.256 \\ 
Nitsche&  $8$ &  2.118 & 2.055 & 2.039 & 2.035 \\
&  $16$ &  2.174 & 2.487 & 2.723 & 2.887 \\ 
    \hline
\end{tabular}
\end{table}

\begin{figure}[htbp]
	\centering
	\subfigure[Conforming IGA]{\includegraphics[width=0.32\textwidth]{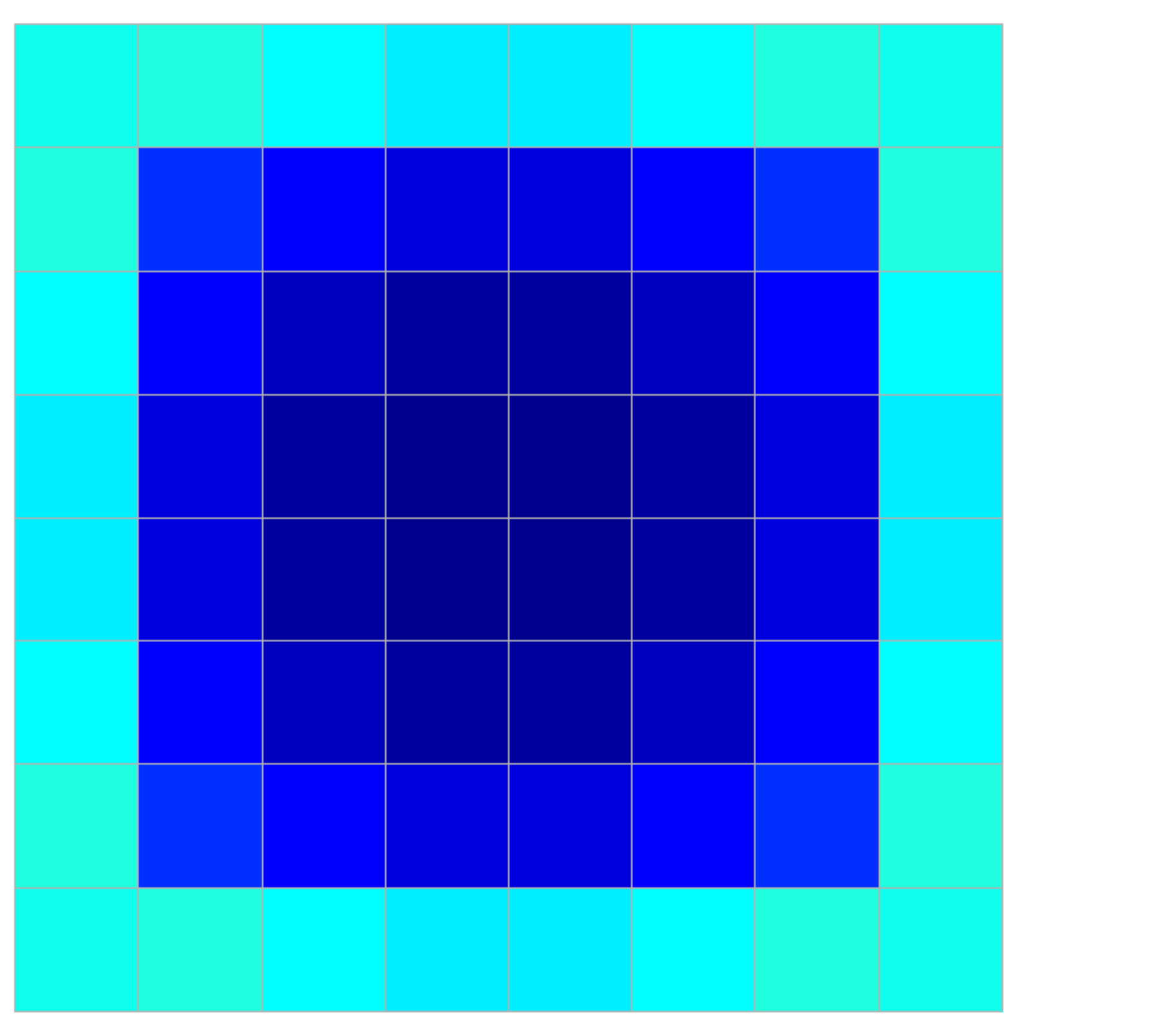}}
    \subfigure[Standard Nitsche's coupling]{\includegraphics[width=0.32\textwidth]{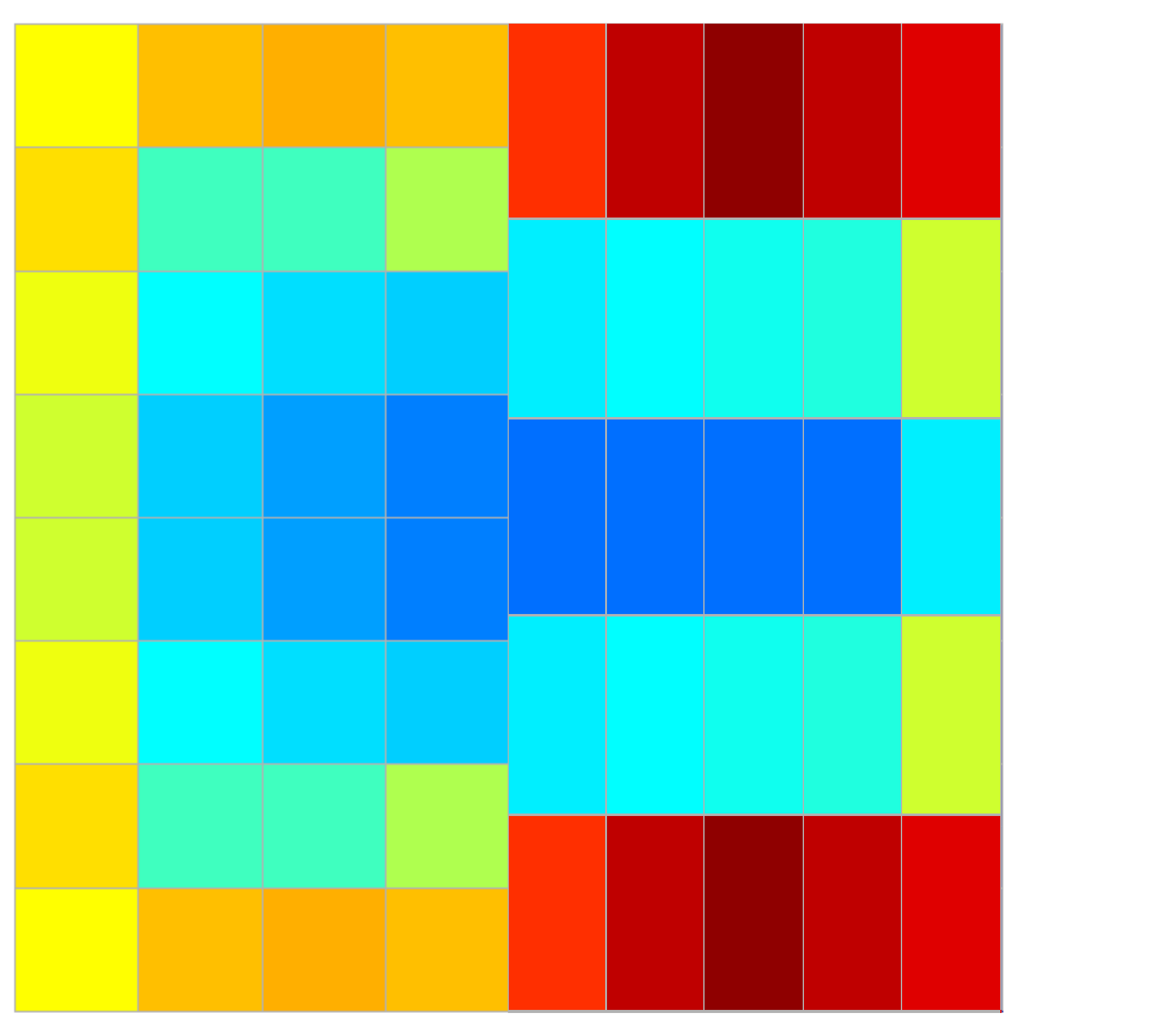}}
    \subfigure[Skew-symmetric Nitsche's coupling]{\includegraphics[width=0.32\textwidth]{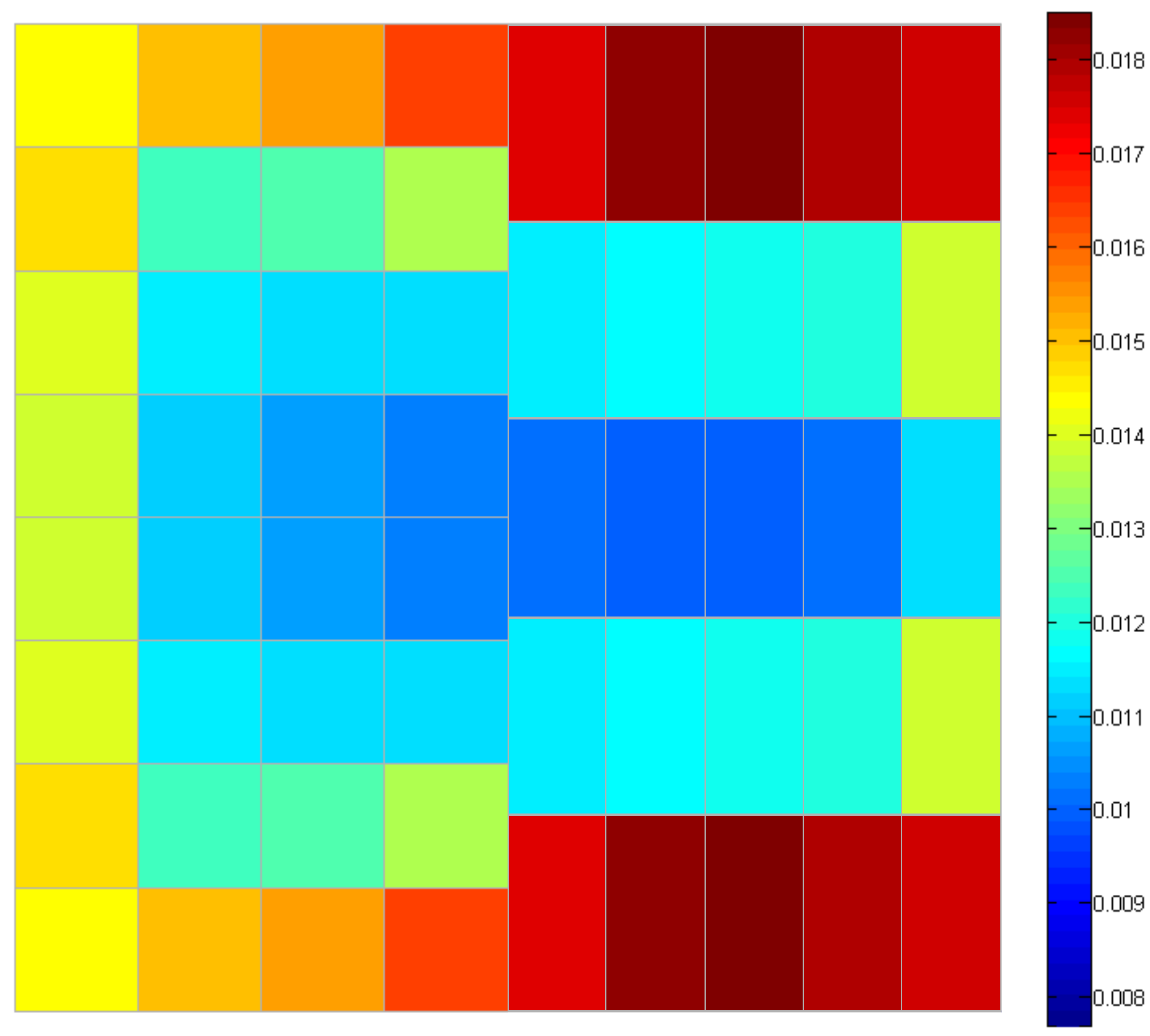}}
	\caption{Patch coupling of two plates: element-wise relative errors of displacement field in $L^2$ norm $||\bu^h - \bu^{\text{ref}}||_{L^2(\Omega^e)} / ||\bu^{\text{ref}}||_{L^2(\Omega^e)}$. For conforming IGA the mesh is $8\times 8$, for two patch coupling the meshes are $4\times 8$ and $5\times 5$.}
    \label{couple_influence_err_L2}
\end{figure}

The element-wise relative errors in $L^2$ norm for the displacement field are plotted in \fref{couple_influence_err_L2} %for one special case:
for specific choices of the meshes:
for conforming IGA we use a $8\times 8$ mesh as before, whereas for Nitsche's formulations a 
$4\times 8$ mesh in the left patch and a $5\times 5$ mesh in the right patch.
Generally the errors due to Nitsche's patch coupling, though acceptable, are larger than conforming IGA,
and the results of the standard and skew-symmetric Nitsche's coupling are comparable.

\subsubsection{Patch coupling of an annular plate}

\begin{figure}[htbp]
\centering
\def\svgwidth{0.75\columnwidth}
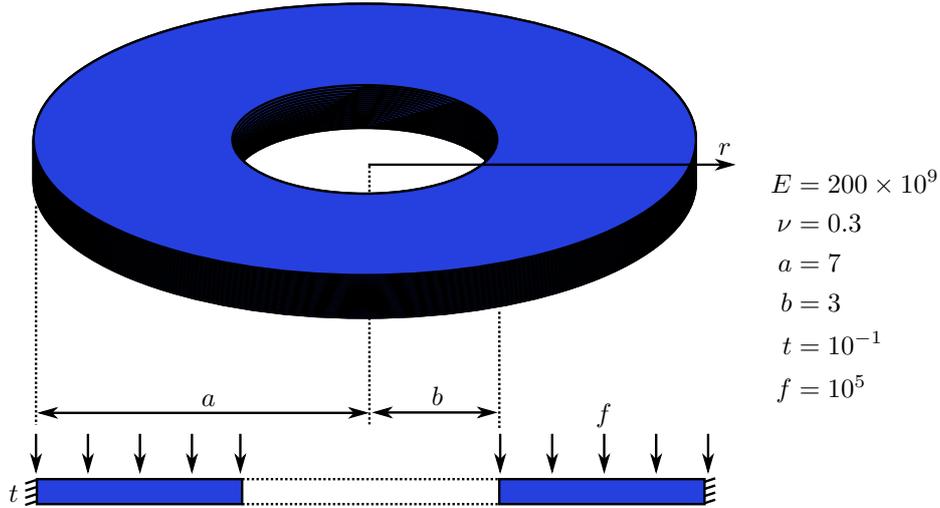
\caption{An annular plate under uniformly distributed load $f$, the outer edge is fixed.}
\label{Fig_anular_mod}
\end{figure}

\begin{figure}[htbp]
\centering
\def\svgwidth{0.75\columnwidth}
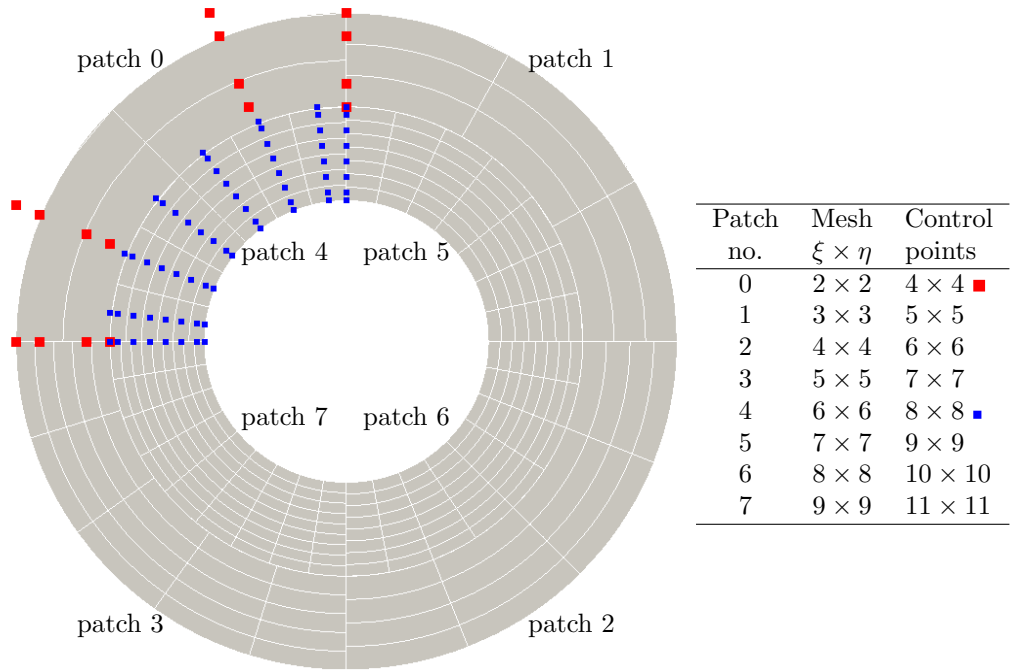
\caption{The annular plate is divided into 8 non-matching patches with 4 curved interfaces and 8 straight interfaces.}
\label{Fig_annular_mesh}
\end{figure}

In this subsection we show that for several curved interfaces that needed to be glued,
Nitsche's formulation is accurate regarding the smoothness and errors of the coupled physical field.
\fref{Fig_anular_mod} shows an annular plate subjected to a uniformly distributed load $f$,
with thickness $t=10^{-1}$, Young's modulus $E=200\times 10^9$ and Poisson's ratio $\nu=0.3$,
The outer edge of the plate is fixed and inner edge is free.
The model is divided into 8 patches with different meshes,
then it leads to 4 curved interfaces and 8 straight interfaces as shown in \fref{Fig_annular_mesh},
thus this problem still corresponds to setting described in \ref{sub:interface}.
We make use of an IGA approximation with bi-quadratic degenerated Reissner-Mindlin elements \cite{adam2015improved}.
The stabilization parameter $\gamma_0$ is still computed by solving the generalized eigenvalue problem \eqref{eigenvalue}.
The analytical solution for the transverse deflection can be found in \cite{young2002roark},
\iffalse
\begin{equation}
w(r)=w_b+\theta_brF_1-\frac{pr^4}{D}G_{11}, \quad r\in [b,a]
\end{equation}
where
\begin{equation}
\begin{split}
w_b &= -\frac{pa^4}{D}\left( \frac{C_1L_{14}}{C_4}-L_{11} \right), \\
\theta_b &= \frac{pa^3}{D}\frac{L_{14}}{C_4}, 
\end{split}
\end{equation}
and
\begin{equation}
\begin{split}
C_1 &= \frac{1+\nu}{2}\frac{b}{a}\ln\frac{a}{b}+\frac{1-\nu}{4}\left( \frac{a}{b}-\frac{b}{a}\right), \\
C_4 &= \frac{1}{2} \left[ (1+\nu)\frac{b}{a}+(1-\nu)\frac{a}{b} \right], \\
L_{11} &= \frac{1}{64}\Bigg\{ 1+4\left(\frac{b}{a}\right)^2-5\left(\frac{b}{a}\right)^4-4\left(\frac{b}{a}\right)^2 \left[2+\left(\frac{b}{a}\right)^2 \right] \ln\frac{a}{b} \Bigg\},   \\
L_{14} &= \frac{1}{16}\left[ 1-\left(\frac{b}{a}\right)^4 -4\left(\frac{b}{a}\right)^2 \ln\frac{a}{b} \right],  \\
F_1 &= \frac{1+\nu}{2}\frac{b}{r} \ln\frac{r}{b} + \frac{1-\nu}{4}\left( \frac{r}{b}-\frac{b}{r} \right), \\
G_{11} &= \frac{1}{64}\Bigg\{ 1+4\left(\frac{b}{r}\right)^2-5\left(\frac{b}{r}\right)^4-4\left(\frac{b}{r}\right)^2 \left[2+\left(\frac{b}{r}\right)^2 \right] \ln\frac{r}{b} \Bigg\},
\end{split}
\end{equation}
and $D=\frac{Eh^3}{12(1-\nu^2)}$ is the flexural rigidity.
\fi
and the largest deflection is $w^{\text{ref}}_{r=b}=-0.10409$.

\begin{figure}[htbp]
	\centering
	\subfigure[Deflection field $\bu^h$]{\includegraphics[width=0.45\textwidth]{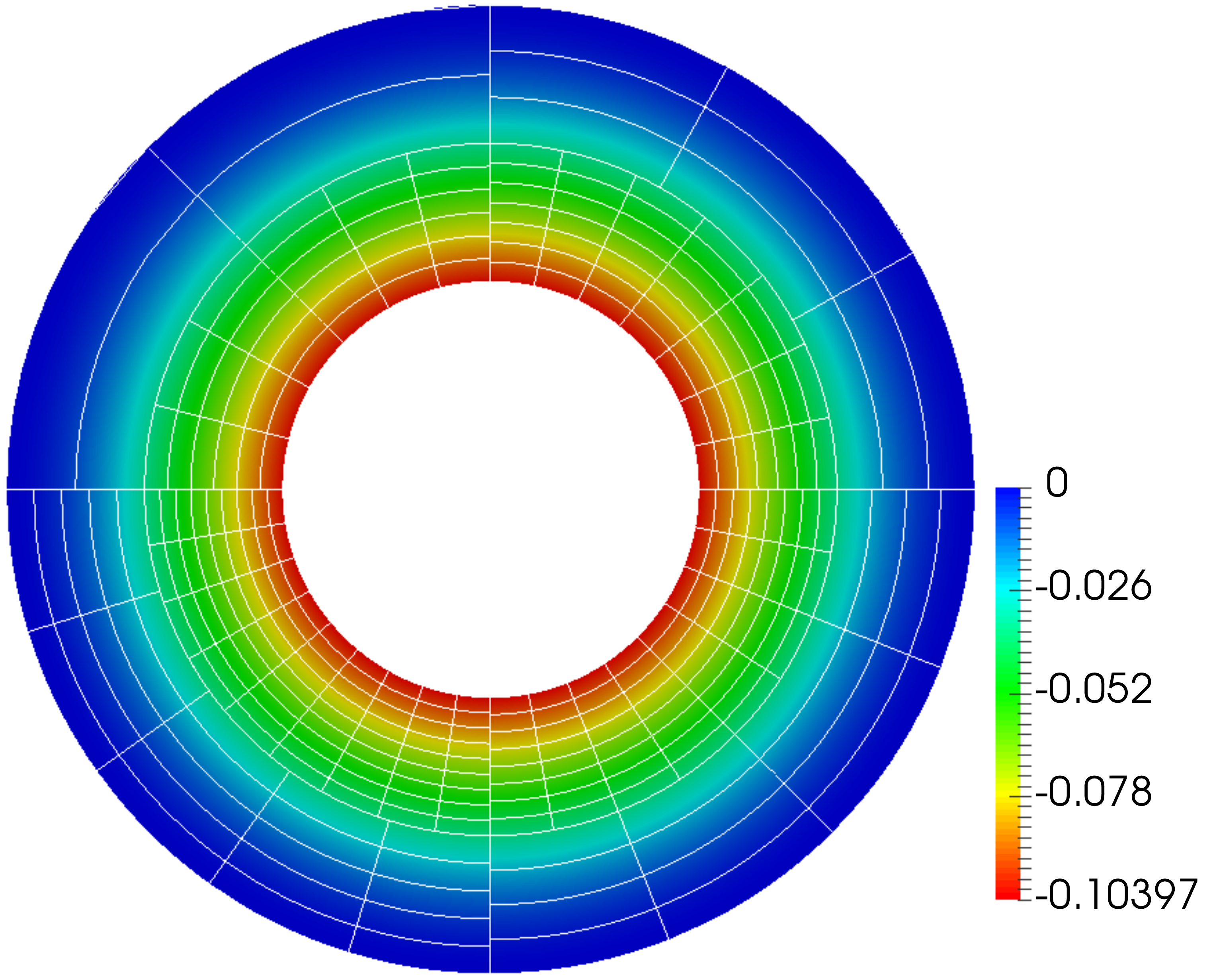}}
    \subfigure[Absolute error $| \bu^h-\bu^{\text{ref}} |$]{\includegraphics[width=0.45\textwidth]{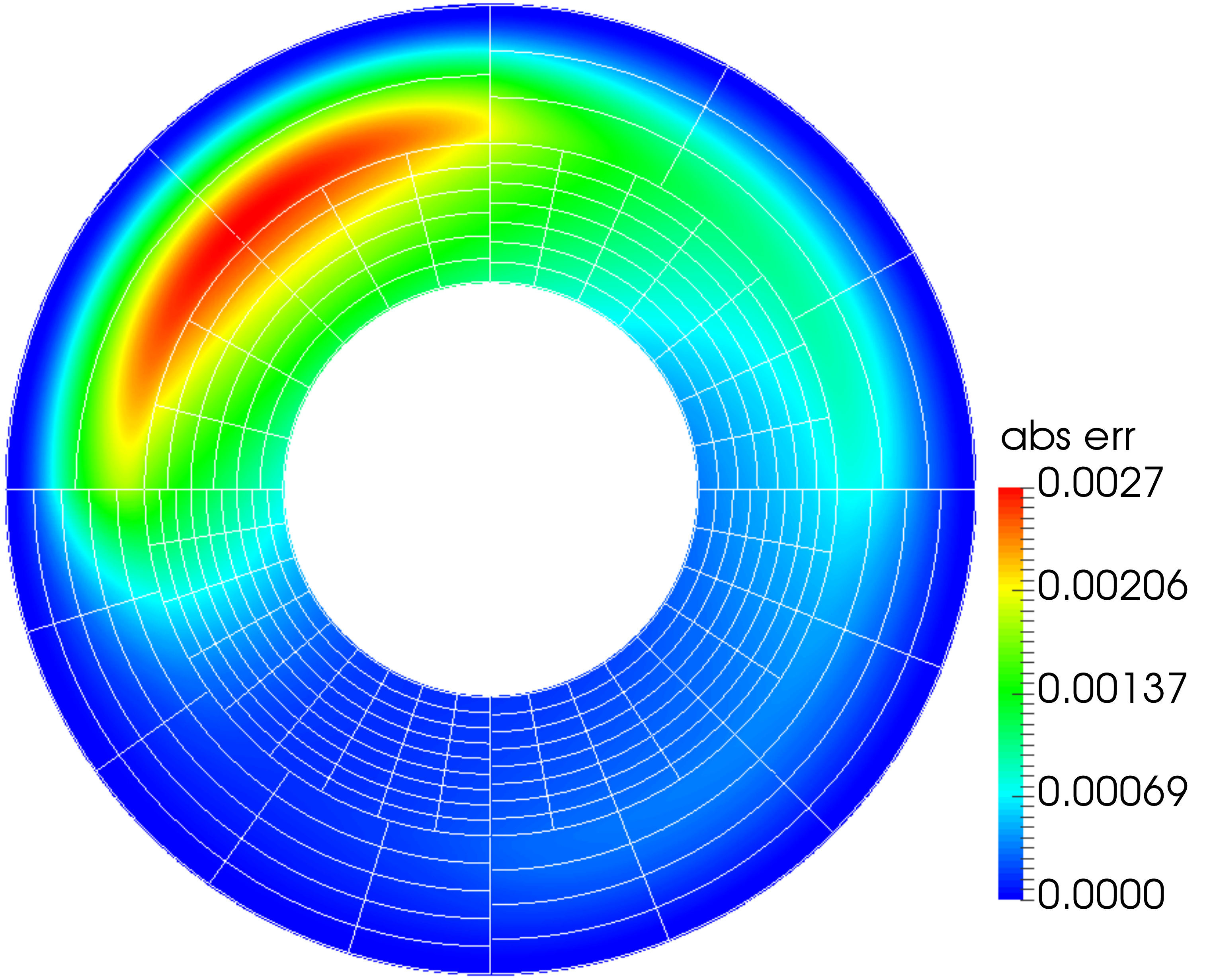}}
	\caption{Results of the annular plate, using standard Nitsche's formulation and the mesh given in \fref{Fig_annular_mesh}.}
	\label{Fig_annular_field_error_stand}
\end{figure}

\begin{figure}[htbp]
	\centering
	\subfigure[Deflection field $\bu^h$]{\includegraphics[width=0.45\textwidth]{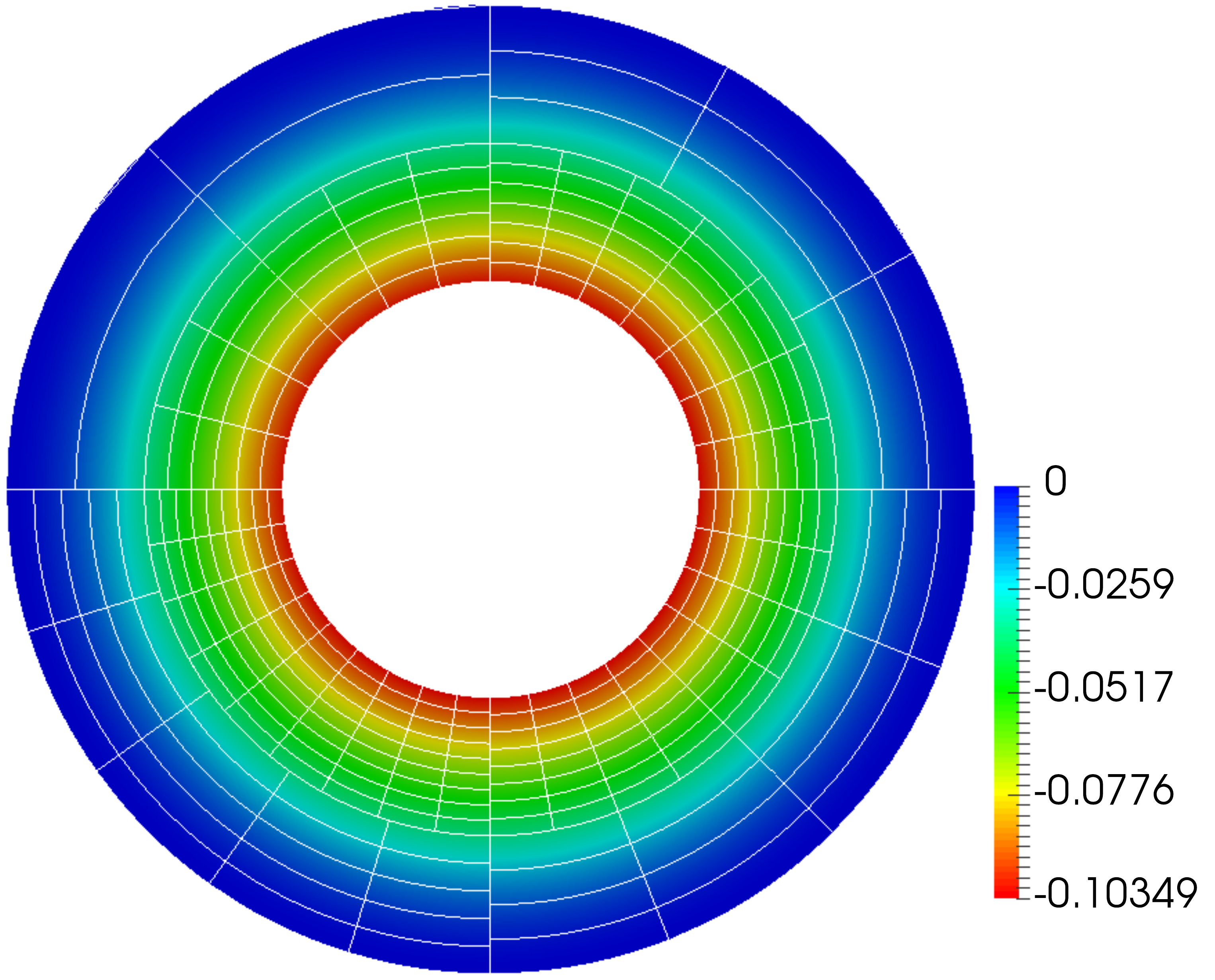}}
    \subfigure[Absolute error $| \bu^h-\bu^{\text{ref}} |$]{\includegraphics[width=0.45\textwidth]{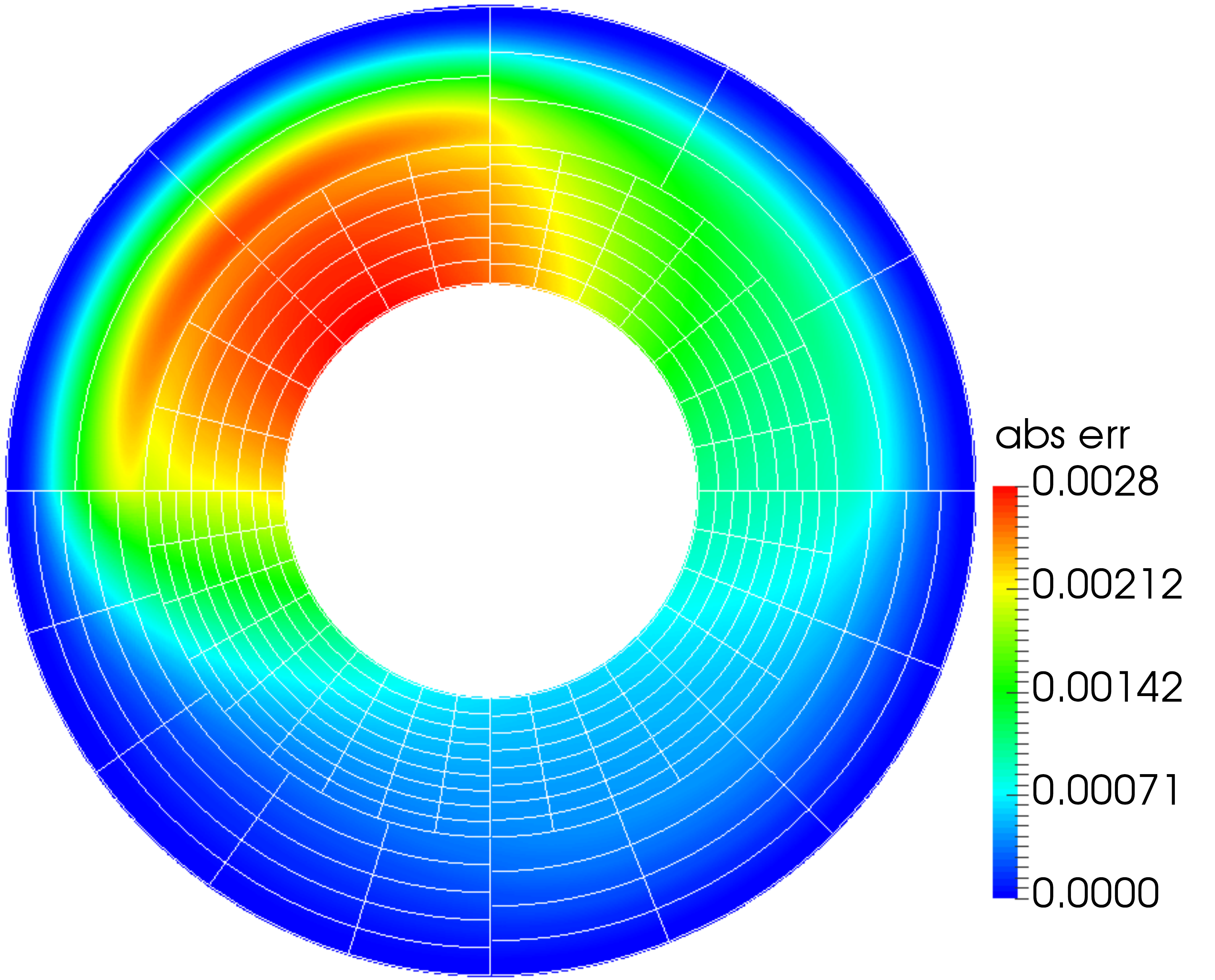}}
	\caption{Results of the annular plate, using skew-symmetric Nitsche's formulation and the mesh given in \fref{Fig_annular_mesh}.}
	\label{Fig_annular_field_error_skew}
\end{figure}

The results obtained with standard and the skew-symmetric Nitsche's formulations are plotted in \fref{Fig_annular_field_error_stand} and \fref{Fig_annular_field_error_skew}, respectively.
Although the model is discretized with different meshes,
the deflection field is quite smooth,
implying that Nitsche's method is effective to glue curved patches with non-conforming meshes.
Visible errors are noticed at patch 0 and patch 4,
which makes sense because the mesh for patch 0 is relatively coarse.
The errors of the largest deflection for the skew-symmetric and the standard Nitsche's formulations are $-0.58\%$ and $-0.12\%$, respectively.

\newpage

\subsubsection{Patch coupling effects: modal analysis}

\begin{figure}[htbp]
\centering
\def\svgwidth{1\columnwidth}
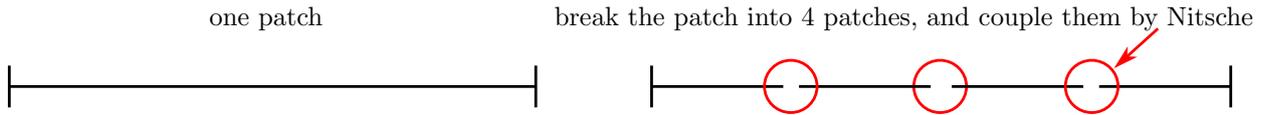
\caption{Rod model. On the left, one patch rod model is adopted for comparison. On the right, the rod is artificially broken into four patches and the additional three interfaces are coupled by Nitsche's method.}
\label{fig_rod}
\end{figure}

To study whether additional effects are introduced in modal analysis by Nitsche's coupling,
the longitudinal vibration of a rod \cite{cottrell2006isogeometric} is considered in \fref{fig_rod}.
On the right side of the figure,
the rod model is broken into four identical patches and they are coupled by Nitsche's method.
We adapt the framework presented in \ref{sub:interface} to this simpler model and to a vibration setting.
%The governing equations of this problem are
So we find $u : (0,1) \rightarrow \mathbb{R}$ and $\omega^2 > 0$ that solve:
\[
\begin{split}
u_{,xx}+\omega^2 u = 0, \textrm{ on } (0,1),\\
u(0)=u(1)=0,
\end{split}
\]
%where $x\in [0,1]$ denotes the rod's length, 
and we know the exact natural frequencies are
\[
\omega_n = n \pi, \quad n=1,\ldots,N,
\]
with $N$ being the total number of DoFs.

\begin{figure}[htbp]
\centering
	\input{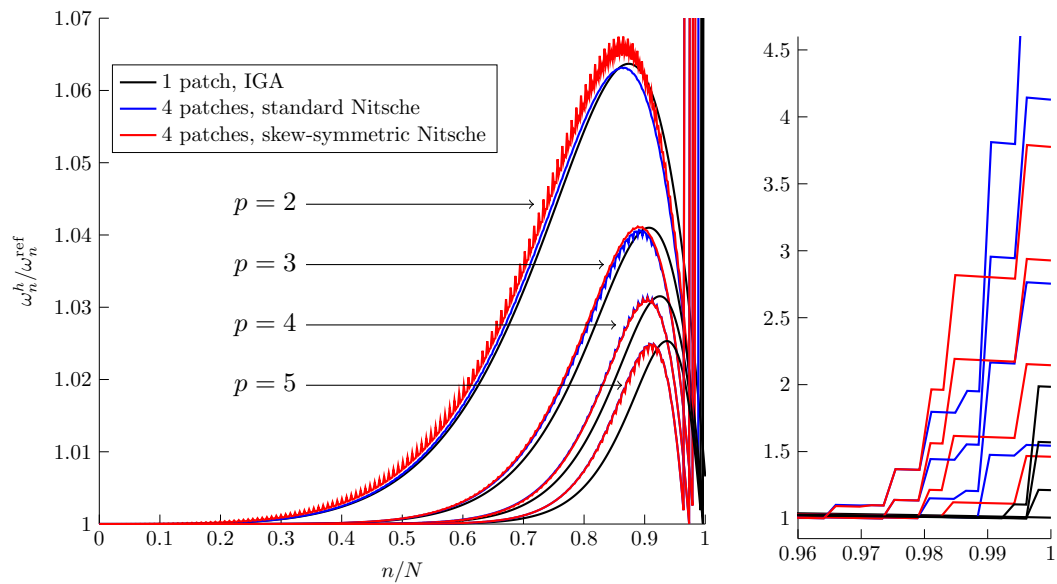}
	\caption{Normalized discrete spectra. As shown in the right figure, the sudden jump of the frequencies identify the ``outliers''. More details are shown in \fref{Fig_vibration_outlier}.}
	\label{Fig_vibration_w}
\end{figure}

\begin{figure}[htbp]
\centering
\begin{tikzpicture}

\tikzset{every picture/.append style={scale=1}}
\setlength\figureheight{12cm} 
\setlength\figurewidth{16cm}

\begin{axis}[%
width=0.951\figurewidth,
height=\figureheight,
at={(0\figurewidth,0\figureheight)},
scale only axis,
unbounded coords=jump,
separate axis lines,
every outer x axis line/.append style={black},
every x tick label/.append style={font=\color{black}},
every x tick/.append style={black},
x tick label style={
    /pgf/number format/.cd,
    fixed,
    fixed zerofill,
    precision=3
  },
xmin=0.965,
xmax=1,
xtick={0.965, 0.97, 0.975, 0.98, 0.985, 0.99, 0.995, 1},
xlabel={$n/N$},
every outer y axis line/.append style={black},
every y tick label/.append style={font=\color{black}},
every y tick/.append style={black},
ymin=1,
ymax=5.7,
ylabel={$\omega_n^h / \omega_n^{\text{ref}}$},
axis background/.style={fill=white},
legend style={at={(0.15,0.5)}, anchor=south west, legend cell align=left, align=left, fill=white}
]

\addplot [ only marks, yellow, mark=10-pointed star, mark size=4.5pt]
  table[row sep=crcr]{%
0.95 0.8\\
};
\addlegendentry{IGA $p=2$}

\addplot [only marks, blue, mark=pentagon, mark size=4.5pt]
  table[row sep=crcr]{%
0.99033	1.475699904\\
0.99226	1.472828893\\
0.9942	1.469969031\\
0.99613	1.549820197\\
0.99807	1.546822479\\
1	1.543836335\\
};
\addlegendentry{standard Nitsche $p=2$}

\addplot [only marks, red, mark=pentagon*, mark size=3pt]
  table[row sep=crcr]{%
0.98453	1.121402661\\
0.98646	1.119208136\\
0.98839	1.117022182\\
0.99033	1.114844751\\
0.99226	1.112675792\\
0.9942	1.110515257\\
0.99613	1.467120254\\
0.99807	1.464282497\\
1	1.461455697\\
};
\addlegendentry{skew-symmetric Nitsche $p=2$}

\addplot [only marks, cyan, mark=x, mark size=4.5pt]
  table[row sep=crcr]{%
0.998046875	1.214424944\\
1	1.212057644\\
};
\addlegendentry{IGA $p=3$}

\addplot [only marks, magenta, mark=square, mark size=4.5pt]
  table[row sep=crcr]{%
0.980806142	1.151228242\\
0.982725528	1.148984133\\
0.984644914	1.146748755\\
0.986564299	1.207350624\\
0.988483685	1.205010797\\
0.990403071	2.165842475\\
0.992322457	2.161661312\\
0.994241843	2.157496261\\
0.996161228	2.764693743\\
0.998080614	2.759387229\\
1	2.754101046\\
};
\addlegendentry{standard Nitsche $p=3$}

\addplot [only marks, green, mark=square*, mark size=3pt]
  table[row sep=crcr]{%
0.980806142	1.214424944\\
0.982725528	1.212057644\\
0.984644914	1.617262346\\
0.986564299	1.614122031\\
0.988483685	1.610993887\\
0.990403071	1.607877845\\
0.992322457	1.604773834\\
0.994241843	1.601681784\\
0.996161228	2.15334723\\
0.998080614	2.149214126\\
1	2.145096857\\
};
\addlegendentry{skew-symmetric Nitsche $p=3$}

\addplot [only marks, olive, mark=star, mark size=4.5pt]
  table[row sep=crcr]{%
0.998050682	1.572173858\\
1	1.569115154\\
};
\addlegendentry{IGA $p=4$}

\addplot [only marks, violet, mark=pentagon, mark size=4.5pt]
  table[row sep=crcr]{%
0.975238095	1.139158034\\
0.977142857	1.136941773\\
0.979047619	1.134734119\\
0.980952381	1.444315093\\
0.982857143	1.441521447\\
0.984761905	1.438738587\\
0.986666667	1.553998438\\
0.988571429	1.551009979\\
0.99047619	2.955723622\\
0.992380952	2.950061317\\
0.994285714	2.944420664\\
0.996190476	4.144183847\\
0.998095238	4.136290163\\
1	4.128426494\\
};
\addlegendentry{standard Nitsche $p=4$}

\addplot [only marks, orange, mark=pentagon*, mark size=3pt]
  table[row sep=crcr]{%
0.975238095	1.139158034\\
0.977142857	1.136941773\\
0.979047619	1.134734119\\
0.980952381	1.563033312\\
0.982857143	1.560010037\\
0.984761905	2.192488613\\
0.986666667	2.188264165\\
0.988571429	2.184055965\\
0.99047619	2.179863919\\
0.992380952	2.175687934\\
0.994285714	2.171527919\\
0.996190476	2.938801541\\
0.998095238	2.933203823\\
1	2.92762739\\
};
\addlegendentry{skew-symmetric Nitsche $p=4$}

\addplot [only marks, black, mark=o, mark size=4.5pt]
  table[row sep=crcr]{%
0.994163424	1.106054507\\
0.996108949	1.103898455\\
0.998054475	1.98761555\\
1	1.983756102\\
};
\addlegendentry{IGA $p=5$}

\addplot [only marks, purple, mark=triangle, mark size=4.5pt]
  table[row sep=crcr]{%
0.965973535	1.102030761\\
0.967863894	1.099882553\\
0.969754253	1.097742704\\
0.971644612	1.099611471\\
0.973534972	1.097480441\\
0.975425331	1.368961422\\
0.97731569	1.366318639\\
0.979206049	1.363686041\\
0.981096408	1.797070195\\
0.982986767	1.793620924\\
0.984877127	1.790184869\\
0.986767486	1.95341184\\
0.988657845	1.949683955\\
0.990548204	3.811406485\\
0.992438563	3.804160465\\
0.994328922	3.796941944\\
0.996219282	5.64573189\\
0.998109641	5.635059429\\
1	5.624427242\\
};
\addlegendentry{standard Nitsche $p=5$}

\addplot [only marks, mark=triangle*, teal, mark size=3pt]
  table[row sep=crcr]{%
0.965973535	1.093057807\\
0.967863894	1.09092709\\
0.969754253	1.088804664\\
0.971644612	1.099611471\\
0.973534972	1.097480441\\
0.975425331	1.368961422\\
0.97731569	1.366318639\\
0.979206049	1.363686041\\
0.981096408	1.964681524\\
0.982986767	1.960910542\\
0.984877127	2.818093255\\
0.986767486	2.812704932\\
0.988657845	2.807337174\\
0.990548204	2.801989865\\
0.992438563	2.796662888\\
0.994328922	2.791356128\\
0.996219282	3.789750766\\
0.998109641	3.782586776\\
1	3.77544982\\
};
\addlegendentry{skew-symmetric Nitsche $p=5$}

\end{axis}
\end{tikzpicture}
\begin{tikzpicture}
\begin{axis}[
ylabel={Number of outliers},
xlabel={order $p$},
xmin=1.5,
xmax=5.25,
xtick={2,3,4,5},
ymin=0,
ymax=20,
stack plots=y,
/tikz/ybar,
/pgf/bar width={0.3},
legend style={at={(0,1.3)}, anchor=south west, legend cell align=left, align=left, fill=white}
]

\legend{IGA,standard Nitsche,skew-symmetric Nitsche}

\addplot [line width=1.0pt, mark=o, black, mark size=3pt]
coordinates {(2,0) (3,2) (4,2) (5,4) };

\addplot [line width=1.0pt, mark=x, blue, mark size=4pt]
coordinates {(2,6) (3,9) (4,12) (5,15) };

\addplot [line width=1.0pt, mark=square, red, mark size=3pt]
coordinates {(2,3) (3,0) (4,0) (5,0) };

\end{axis}

\node[black] at (0.3,0.21) {0}; %red
\node[blue] at (0.3,1.4) {6}; %red
\node[red] at (0.3,2.1) {9}; %red

\node[black] at (1.7,0.5) {2}; %red
\node[blue] at (1.7,2.35) {11}; %red
\node[red] at (1.7,2.7) {11}; %red

\node[black] at (3.2,0.5) {2}; %red
\node[blue] at (3.15,3) {14}; %red
\node[red] at (3.15,3.35) {14}; %red

\node[black] at (4.7,0.9) {4}; %red
\node[blue] at (4.6,4.05) {19}; %red
\node[red] at (4.6,4.4) {19}; %red

\end{tikzpicture}
	\caption{Normalized outlier frequencies. The model coupled by Nitsche's leads to a larger number of ``outliers''. The corresponding outlier eigenmodes of $p=2$ are plotted in \fref{Fig_vibration_eigenmode_stand} for the standard Nitsche's formulation and in \fref{Fig_vibration_eigenmode_skew} for the skew-symmetric Nitsche's formulation.}
	\label{Fig_vibration_outlier}
\end{figure}
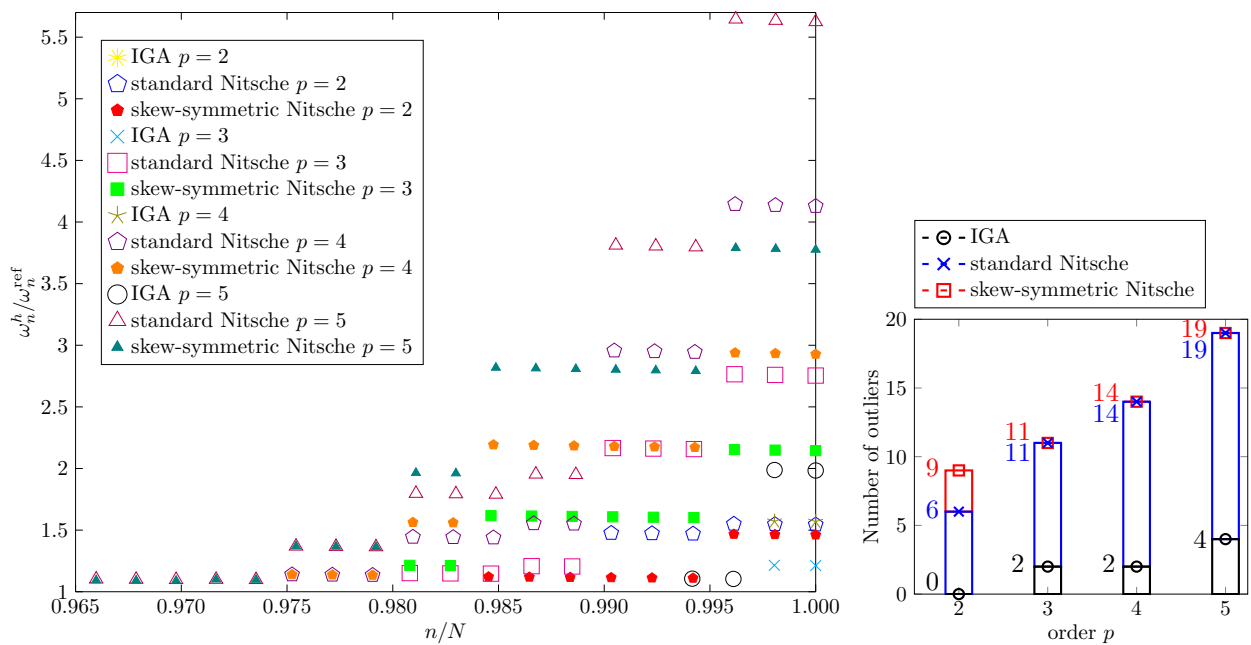

%The discrete spectra is independent of the total number of DoFs $N$ if it is normalized by $N$.
%To get a relatively smooth spectra the mesh should be very well refined,
%here 518 control points are used.
The discrete spectra is normalized by $N$,
and the normalized discrete spectra is given in \fref{Fig_vibration_w},
showing that Nitsche's coupled model is almost as accurate as conforming IGA for lower frequencies,
specifically $n/N<0.2$ for $p=2$ and $n/N<0.5$ for $p>2$. %higher orders.
For higher frequencies,
the skew-symmetric Nitsche's formulation leads to oscillations when $p=2$,
and the solutions from both standard and skew-symmetric Nitsche's formulations are oscillatory when $p=3,4,5$.

\begin{figure}[htbp]
	\centering
	\input{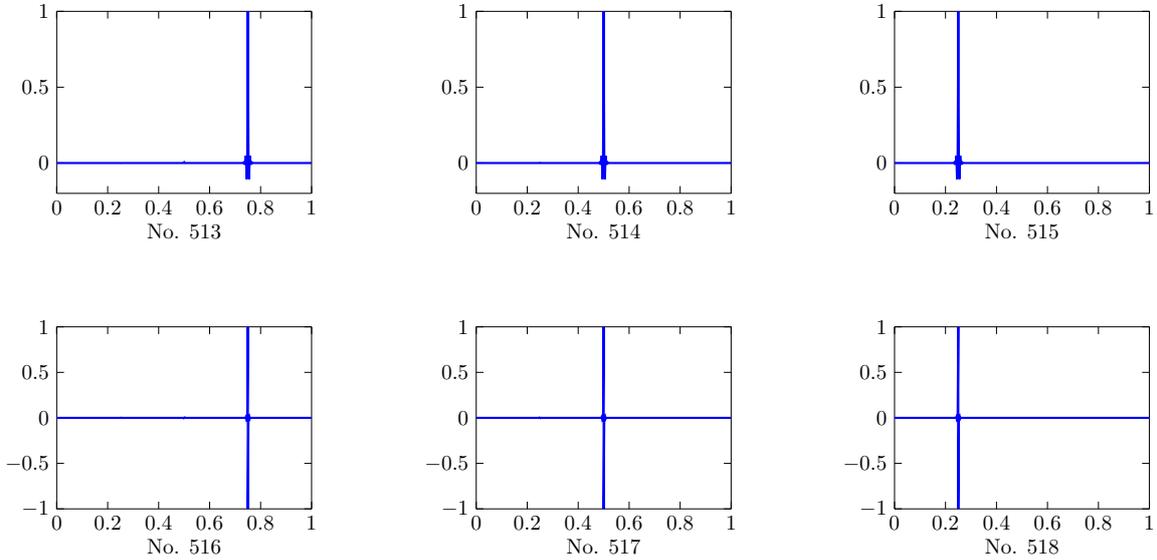}
	\caption{``Outlier'' (last 6 out of 518) eigenmodes obtained by the standard Nitsche's formulation with $p=2$. The longitudinal deformations are plotted along the vertical axis. The eigenmodes that correspond to the ``outlier'' are highly localized at the coupled interfaces.}
	\label{Fig_vibration_eigenmode_stand}
\end{figure}

\begin{figure}[htbp]
	\centering
	\input{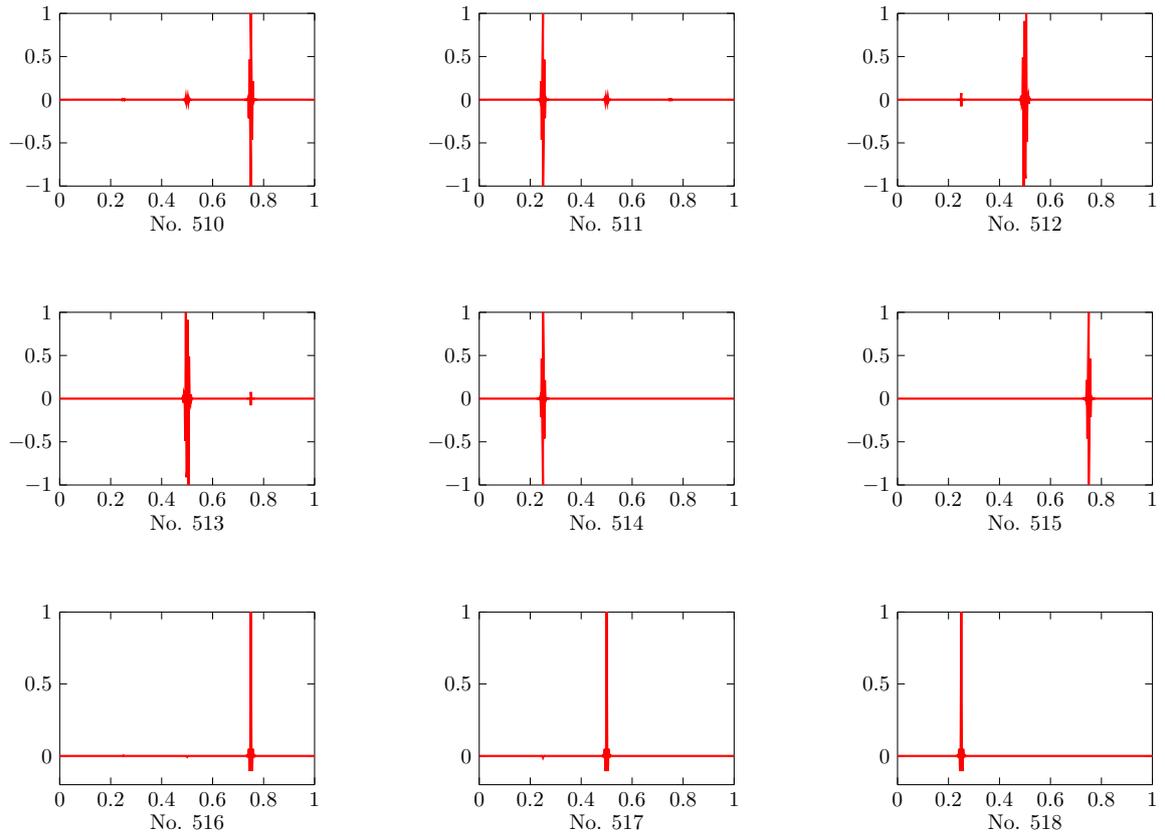}
	\caption{``Outlier'' (last 9 out of 518) eigenmodes obtained by the skew-symmetric Nitsche's formulation with $p=2$. The longitudinal deformations are plotted along the vertical axis. The eigenmodes that correspond to the ``outlier'' are highly localized at the coupled interfaces.}
	\label{Fig_vibration_eigenmode_skew}
\end{figure}

At the very end of the spectra,
the sudden jumps of the frequencies are known as the ``outliers'' \cite{cottrell2006isogeometric,cazzani2016analytical,horger2017improved},
as shown in the enlarged figure on the right side of \fref{Fig_vibration_w}.
These ``outlier'' frequencies are drawn in \fref{Fig_vibration_outlier},
and the number of outliers is also counted.
The ``outlier'' frequencies are captured near $n/N=1$ by the conforming IGA,
and the model coupled by Nitsche's formulation leads to a larger number of ``outliers''.
The corresponding outlier eigenmodes of $p=2$ are plotted in \fref{Fig_vibration_eigenmode_stand} for the standard Nitsche's formulation
and in \fref{Fig_vibration_eigenmode_skew} for the skew-symmetric Nitsche's formulation.
To achieve a better demonstration of the eigenmodes,
we plot the longitudinal deformation along the vertical axis.
%The ``outliers'' obtained with the skew-symmetric Nitsche's formulation are three times more numerous than for the standard Nitsche's formulation when $p=2$.
For the skew-symmetric formulation,
there exist 3 pairs of symmetrical eigenmodes,
i.e. no. 510 and no. 511, no. 512 and no. 513, no. 514 and no. 515.
These results imply that the action of coupling the interfaces by Nitsche's method increases the number of ``outliers'',
and the eigenmodes correspond to these ``outliers'' are highly localized at these coupled interfaces.

\subsection{Frictionless contact}
\label{sub:numcontact}

In the following examples only NURBS basis functions of order $p=q=2$ are employed,
%because they %are already capable of representing a wider variety of geometries.
%moreover to our experience using basis functions of order three or four may lead to numerically unsteadiness (see Remark 4).
since this approximation order is sufficient for the wide majority of contact problems (see Remark \ref{rem-regularity}).
Note that for contact problems the skew-symmetric Nitsche's formulation can not be parameter-free anymore (see Remark \ref{rem-contact-pfree}).
%$\gamma_0$ is used in both standard and skew-symmetric Nitsche's formulations.
%and the value of $\gamma_0^{\textrm{ref}}$ which comes from the generalized eigenvalue problem \eqref{eigenvalue} is adopted  as a reference.

%\red{It is arguable: the interest of skew-symmetric Nitsche is that you may take a lower value of $\gamma_0$ and it should lead to the same accuracy, but with fewer Newton iterations. You may try it \dots (?)}

\subsubsection{Hertz contact}

\begin{figure}[htbp]
\centering
\def\svgwidth{0.7\columnwidth}
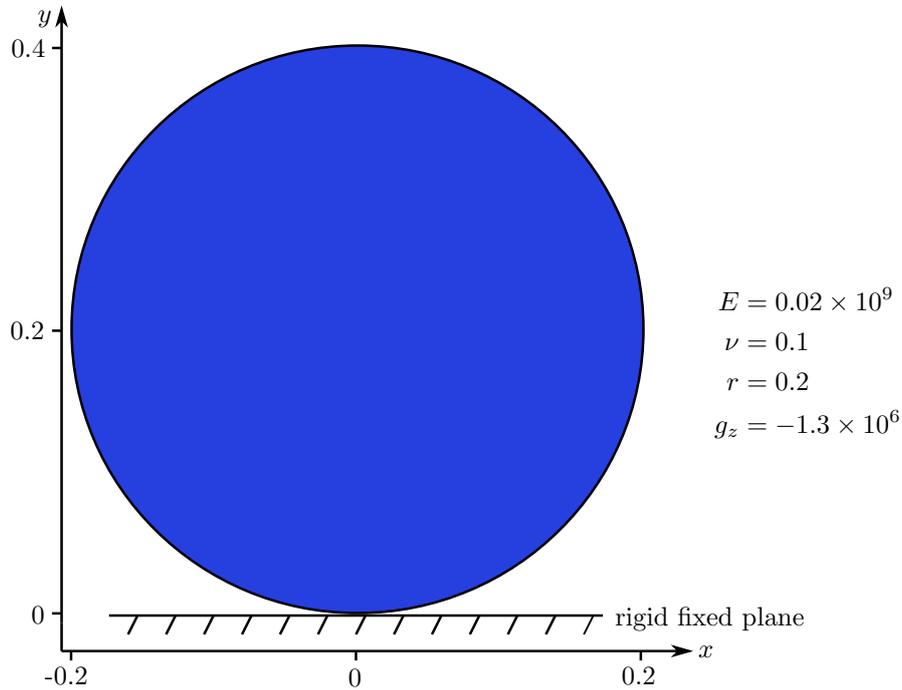
\caption{Hertz contact example, the contact between a linear elastic disc under a vertical gravity force and a rigid fixed plane. Values for the Young modulus $E$, Poisson ratio $\nu$, radius $r$ and gravity load $g_z$ are provided.}
\label{fig_hertz}
\end{figure}

In this section we show that the proposed contact formulation described in Section \ref{sub:contact} is able to predict the contact pressure distribution versus contact width to some degree of accuracy,
and the skew-symmetric formulation is robust w.r.t. the choice of the stabilization parameter $\gamma_0$.
The Hertz contact assumes an elastic frictionless contact without adhesive forces between two cylinders with the same height, radii and elasticity moduli.
To simplify, one cylinder (master body) is fixed, its material is set to be rigid, and its contact surface is flat.
More precisely we study a contact problem as shown in \fref{fig_hertz},
which analytical solution was provided by Hertz \cite{zavarise2006guide}.
Thus we consider the setting \ref{subsub:biased} and a Signorini-type problem (with a rigid support).
The boundary conditions are the following:
the bottom of the disc is specified as the potential contact boundary,
and the whole disc is subjected to a vertical gravity force $g_z$.
To prevent rigid body motions we fix the horizontal displacement at several control points along axis $x=0$.
%The geometry is modeled using NURBS basis functions of order $p=q=2$,
%and we only h-refine the elements hereafter.

\begin{figure}[htbp]
	\centering
	\includegraphics[width=0.6\textwidth]{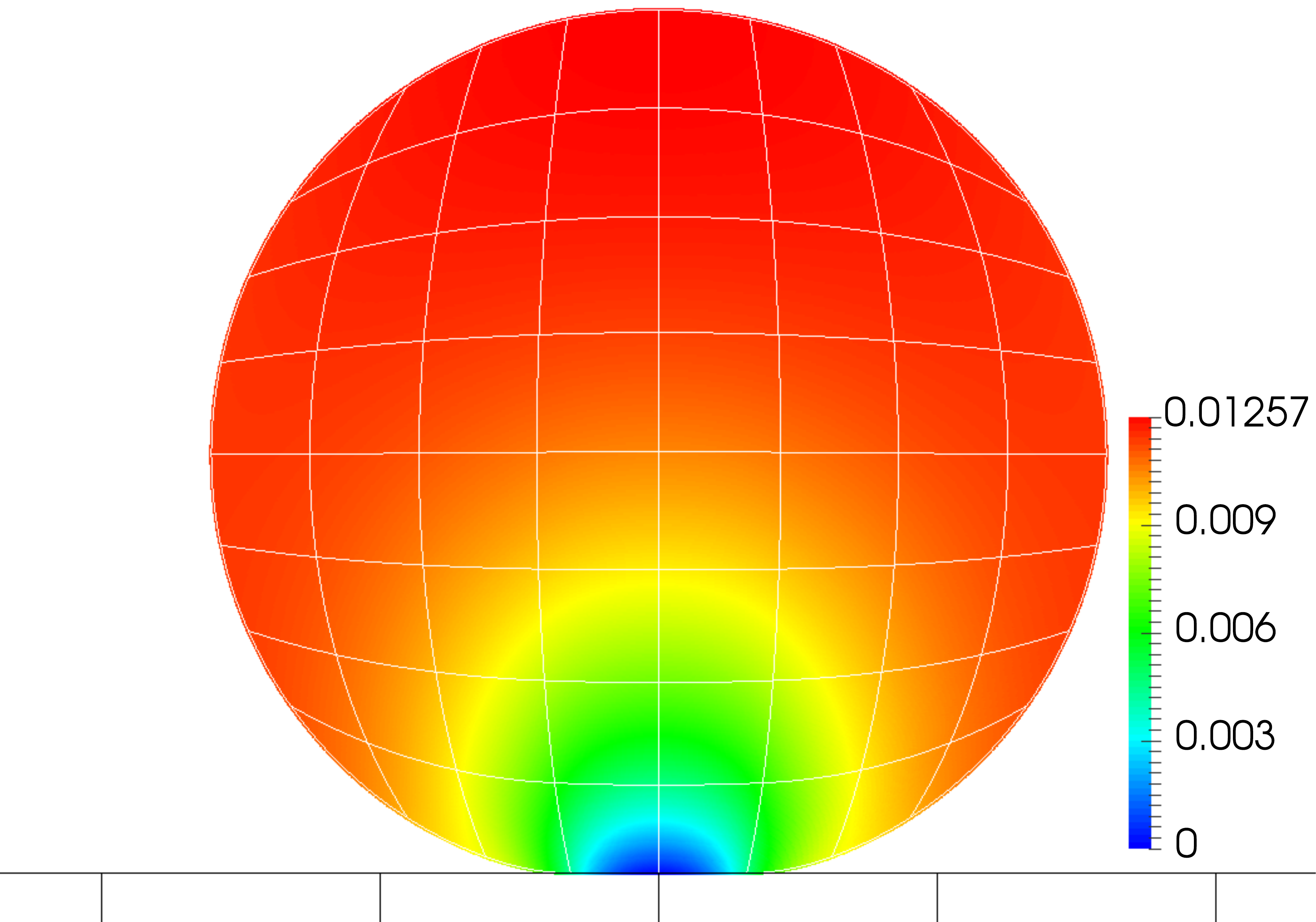}
	\caption{Contour plot of displacement magnitude field obtained by the skew-symmetric Nitsche's method with $8\times 8$ elements.} % and $\gamma = \gamma_0$.}
    \label{hertz_case}
\end{figure}

To begin with we draw the results for one case in \fref{hertz_case} as an illustration,
using the skew-symmetric Nitsche's contact formulation with $8\times 8$ elements. % and $\gamma=\gamma_0$.
The displacement magnitude field shows that the rigid fixed plane has successfully prevented the disc from dropping down,
and the contact surface of the elastic disc adjusts itself to match the rigid fixed plane,
resulting in a straight contact surface.

\begin{figure}[htbp]
\centering
\newcommand{\blackline}{\raisebox{2pt}{\tikz{\draw[black,solid,line width = 1.5pt](0,0) -- (5mm,0);}}}
\newcommand{\blueo}{\raisebox{0.5pt}{\tikz{\draw[blue,solid,line width = 1pt](0,0) -- (5mm,0);\node[draw,scale=0.3,circle,blue,fill=blue](){};}}}
\newcommand{\redo}{\raisebox{0.5pt}{\tikz{\draw[red,solid,line width = 1pt](0,0) -- (5mm,0);\node[draw,scale=0.3,circle,red,fill=red](){};}}}
	\input{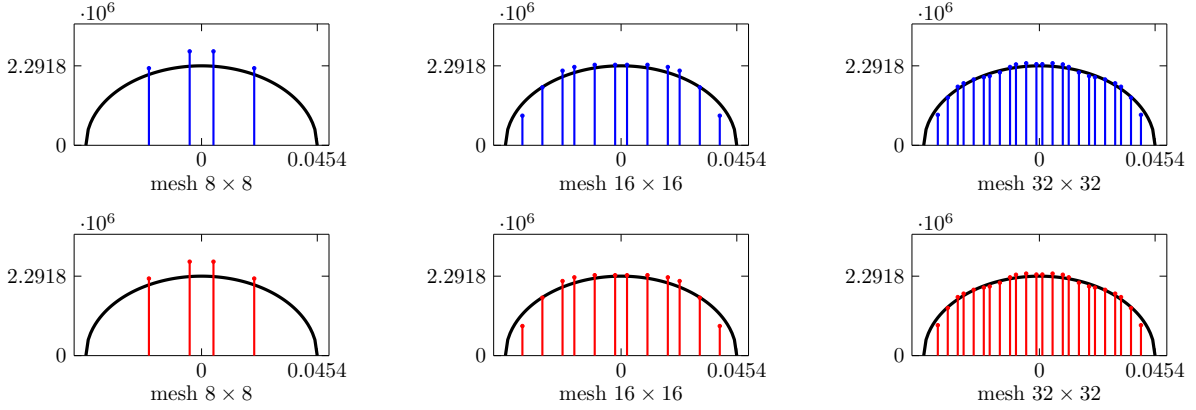}
	\caption{Pressure distribution for Hertz contact. Horizontal axis: contact surface. Vertical axis: contact pressure $p$.
    \protect\\ \protect\hspace{\textwidth}
    \protect\blackline \, analytical solution, \protect\blueo \, the standard Nitsche's method, \protect\redo \, the skew-symmetric Nitsche's method.
    \protect\\ \protect\hspace{\textwidth}
    The contact stresses are calculated at quadrature points.}
    \label{hertz_pressure_distribution_fig}
\end{figure}

The pressure distribution with respect to the contact width is plotted and compared with the analytical solution in \fref{hertz_pressure_distribution_fig},
according to the fact that the stress on the contact surface reaches balance with the contact pressure.
Both standard and skew-symmetric contact formulations are employed. % with $\gamma=\gamma_0$.
The contact stresses in blue and red dots are calculated at quadrature points.
Nitsche's method can properly predict the pressure distribution with respect to the contact width as the mesh is refined.

\begin{table}[htbp]
\centering
\caption{Hertz contact: number of semi-smooth Newton iterations for various $\gamma_0$.} 
\label{hertz_contact_iters}

\medskip

\begin{tabular}{c c ccccc}
 \hline
Method  & $\gamma_0 = $ & Mesh $4\times 4$ & Mesh $8\times 8$ & Mesh $16\times 16$ & Mesh $32\times 32$ & Mesh $64\times 64$ 
\\ 
\noalign{\smallskip}
\hline 
\noalign{\smallskip}
Standard&  $\gamma_0^{\textrm{ref}}$ &  8 & 13 & 21 & 41 & 52 \\ 
\noalign{\smallskip}
Nitsche&  $\gamma_0^{\textrm{ref}} / 10000$ &  6 & 9 & 11 & $>100$ & $>100$ \\
	\noalign{\smallskip}
       &  $\gamma_0^{\textrm{ref}} / 100000$ &  7 & 22 & 9 & 52 & 45 \\
 \noalign{\smallskip}
 \hline
 \noalign{\smallskip}
Skew-symmetric&  $\gamma_0^{\textrm{ref}}$ &  8 & 13 & 21 & 42 & 52 \\ 
\noalign{\smallskip}
Nitsche&  $\gamma_0^{\textrm{ref}} / 10000$ &  6 & 7 & 10 & 9 & 11 \\
	\noalign{\smallskip}
       &  $\gamma_0^{\textrm{ref}} / 100000$ &  7 & 7 & 9 & 9 & 10 \\
 \noalign{\smallskip}
    \hline
\end{tabular}
\end{table}

The number of semi-smooth Newton iterations for various values of the stabilization parameter $\gamma_0$ are displayed in \Tref{hertz_contact_iters}.
We define $\gamma_0^{\textrm{ref}} := 2 \lambda^{h,\mathrm{MAX}}$, and then study the influence of $\gamma_0$ by choosing $\gamma_0 =\gamma_0^{\textrm{ref}}$, $\gamma_0 =\gamma_0^{\textrm{ref}}/10000$, and $\gamma_0 =\gamma_0^{\textrm{ref}}/100000$.
For $\gamma_0 = \gamma_0^{\textrm{ref}}$ the standard Nitsche's formulation and the skew-symmetric one behave similarly,
and the number of iterations increases as the mesh is refined, since the problem becomes stiffer.
For $\gamma_0 =\gamma_0^{\textrm{ref}} /10000$ and $\gamma_0 =\gamma_0^{\textrm{ref}} /100000$
the skew-symmetric formulation remains remarkably robust and converges faster than for $\gamma_0 = \gamma_0^{\textrm{ref}}$. Conversely, convergence is harder to achieve with the standard formulation, especially for finer meshes, because these small values of $\gamma_0$ can not ensure well-posedness anymore. These results are coherent with the behavior observed for FEM in \cite{mlika2017unbiased}.

\begin{figure}[htbp]
\centering
	\subfigure[$\gamma_0 =\gamma_0^{\textrm{ref}}$]{% This file was created by matlab2tikz.
\begin{tikzpicture}

\begin{axis}[%
width=1\figurewidth,
height=1\figureheight,
at={(0\figurewidth,0\figureheight)},
scale only axis,
separate axis lines,
every outer x axis line/.append style={black},
every x tick label/.append style={font=\color{black}},
every x tick/.append style={black},
xmode=log,
xmin=3,
xmax=100,
xminorticks=true,
xlabel={Number of elements per side on slave patch},
every outer y axis line/.append style={black},
every y tick label/.append style={font=\color{black}},
every y tick/.append style={black},
ymode=log,
ymin=1e-2,
ymax=1,
yminorticks=true,
ylabel={$||\bu^h - \bu^{\text{ref}}||_{E(\Omega)} / ||\bu^{\text{ref}}||_{E(\Omega)}$},
axis background/.style={fill=white},
xmajorgrids,
xminorgrids,
ymajorgrids,
yminorgrids,
legend style={at={(0.52,0.79)}, anchor=south west, legend cell align=left, align=left, draw=black}
]

\addplot [color=blue, line width=1.5pt,mark=o]
  table[row sep=crcr]{%
4	4.43054E-01\\
8   2.31449E-01\\
16  8.70145E-02\\
32  4.14345E-02\\
64  2.05459E-02\\
};
\addlegendentry{standard Nitsche}

\addplot [color=red, line width=1.5pt,mark=o]
  table[row sep=crcr]{%
4   4.42835E-01\\
8   2.31473E-01\\
16  8.70143E-02\\
32  4.14346E-02\\
64  2.05463E-02\\
};
\addlegendentry{skew-symmetric Nitsche}

%\draw (20,0.02) -- (5,0.08); %1
%\draw (20,0.02) -- (5,0.32); %2
%\draw (20,0.02) -- (5,1.28); %3

\end{axis}

\node [color=red] at (3.55,3.75) {1.41};

\end{tikzpicture}%
}
	\subfigure[$\gamma_0=\gamma_0^{\textrm{ref}} / 100000$]{% This file was created by matlab2tikz.
\begin{tikzpicture}

\begin{axis}[%
width=1\figurewidth,
height=1\figureheight,
at={(0\figurewidth,0\figureheight)},
scale only axis,
separate axis lines,
every outer x axis line/.append style={black},
every x tick label/.append style={font=\color{black}},
every x tick/.append style={black},
xmode=log,
xmin=3,
xmax=100,
xminorticks=true,
xlabel={Number of elements per side on slave patch},
every outer y axis line/.append style={black},
every y tick label/.append style={font=\color{black}},
every y tick/.append style={black},
ymode=log,
ymin=1e-2,
ymax=1,
yminorticks=true,
ylabel={$||\bu^h - \bu^{\text{ref}}||_{E(\Omega)} / ||\bu^{\text{ref}}||_{E(\Omega)}$},
axis background/.style={fill=white},
xmajorgrids,
xminorgrids,
ymajorgrids,
yminorgrids,
legend style={at={(0.52,0.79)}, anchor=south west, legend cell align=left, align=left, draw=black}
]

\addplot [color=blue, line width=1.5pt,mark=o]
  table[row sep=crcr]{%
4	0.834199\\
8   0.210714\\
16  0.0999895\\
32  0.0429138\\
64  5.46463E-02\\
};
\addlegendentry{standard Nitsche}

\addplot [color=red, line width=1.5pt,mark=o]
  table[row sep=crcr]{%
4   0.486018\\
8   0.202482\\
16  0.0882894\\
32  0.041249\\
64  0.0201706\\
};
\addlegendentry{skew-symmetric Nitsche}

%\draw (20,0.02) -- (5,0.08); %1
%\draw (20,0.02) -- (5,0.32); %2
%\draw (20,0.02) -- (5,1.28); %3

\end{axis}

\node [color=red] at (2.53,3.50) {1.20};
\node [color=blue] at (3.40,3.90) {1.08};

\end{tikzpicture}%
} 
	\caption{Hertz contact: relative errors of displacement field in energy norm.}
	\label{hertz_contact_converge}
\end{figure}
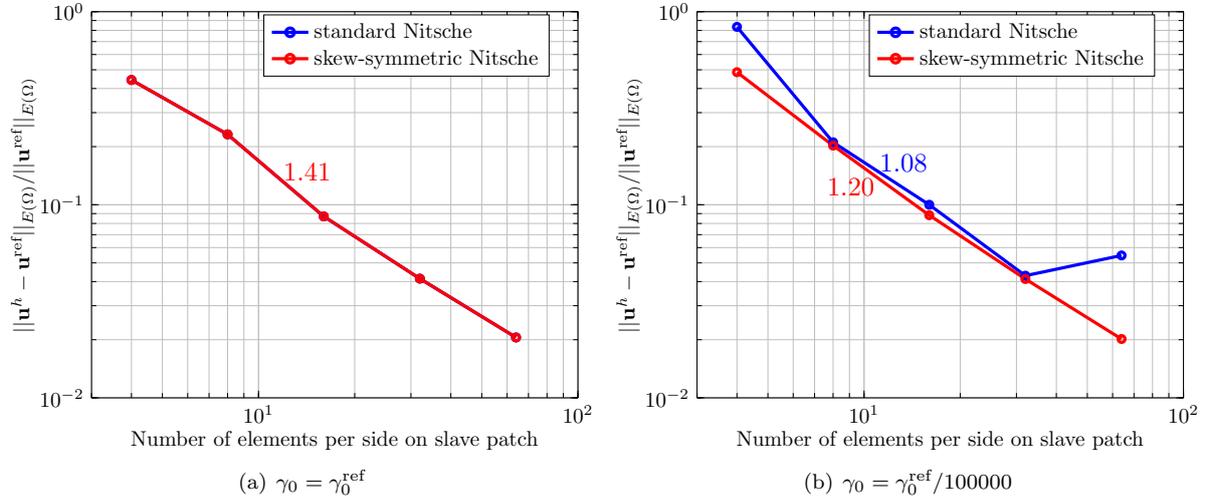

The convergence performance is studied, compared to a reference solution using $128 \times 128$ elements.
For $\gamma_0^{\textrm{ref}}$, the convergence curves for standard and skew-symmetric Nitsche's formulations are similar and we recover a rate of 1.41, close to what is expected for such a problem and a discretization of order 2 (see Remark \ref{rem-regularity} and \cite{chouly2015symmetric}).
%. Note that a "superconvergence" phenomenon is observed, with better rates than predicted by the theory. We have no explanation for its origin.
For $\gamma_0^{\textrm{ref}} /100000$,
as displayed in \fref{hertz_contact_converge},
%the convergence curves for standard Nitsche method with $\gamma_0=E$ are not monotonous and convergence is lost for the finest mesh:
the convergence curves for standard Nitsche's method are perturbed and the convergence rate is lower.
This is in agreement with the theory that standard Nitsche's formulation requires $\gamma_0$ large enough to ensure well-posedness and optimal convergence, see \ref{sub:Schur}.
Conversely,
the convergence performance for skew-symmetric Nitsche's method is not affected by the choice of $\gamma_0$.
All this is in agreement with the theory and previous observations for FEM discretization,
see, e.g., \cite{chouly2015symmetric}.

\subsubsection{Contact between two blocks}

\begin{figure}[htbp]
\centering
\def\svgwidth{0.6\columnwidth}
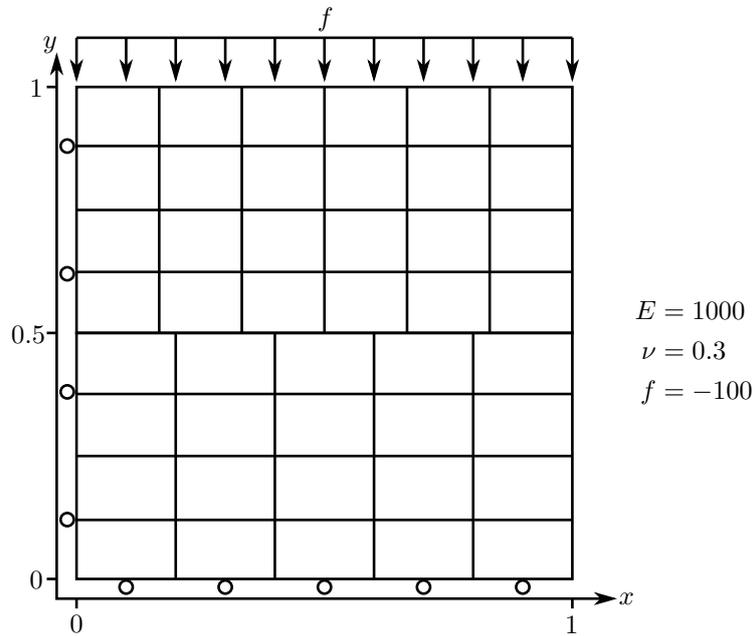
\caption{Two blocks with non-matching meshes. The upper block is subjected to a uniform pressure load $f$ on its top. Both blocks have the same material properties. On the left and at the bottom of the structure, we impose a sliding condition, while on the right we impose no traction. Remark that the problem setup is not symmetric.}
%due to the boundary conditions along the left sides of the two blocks.}
\label{fig_contact_patchtest}
\end{figure}

The goal of this section is to test whether the proposed formulation can properly impose the contact conditions through non-matching elements,
and preliminarily investigate the performance of both biased and unbiased variants of the standard and skew-symmetric contact formulations.
The corresponding framework is described in \ref{subsub:biased}.
The classical contact patch test proposed by Taylor \cite{taylor1991patch} is used to examine the contact algorithm and investigate whether it transfers the constant contact pressure through the contact surface even for non-matching meshes \cite{de2015isogeometric}.
In order to enforce a uniform pressure on the top surface easily,
we adopt a modified patch test \cite{crisfield2000re,el2001stability} as shown in \fref{fig_contact_patchtest}, where the material properties and the boundary conditions are provided.
%Due to the linear strain assumption, we do not set any initial gap between the contact surfaces in order to avoid rigid body motion.
The resulting test setup is similar as in \cite{wriggers2008formulation}.
Because of the boundary conditions along both sides of the two blocks,
this problem setup is not symmetric.
The reference solution of this (plane stress) problem is given by
\begin{equation}
\begin{split}
u_x(x,y) &= 0.03x, \quad u_y(x,y) = -0.1 y, \\
\sigma_{xx}(x,y) &= 0,\; \sigma_{yy} (x,y) = -100,\; \sigma_{xy}(x,y) = 0,
\end{split}
\end{equation}
where $u_x$ and $u_y$ are the components of the displacement $\bu$ and where $\sigma_{xx}$, $\sigma_{xy}$ and $\sigma_{yy}$ are the components of $\boldsymbol{\sigma}(\bu)$.

\begin{figure}[htbp]
	\centering
	\subfigure[Relative error for $u_y$]{\includegraphics[width=0.4\textwidth]{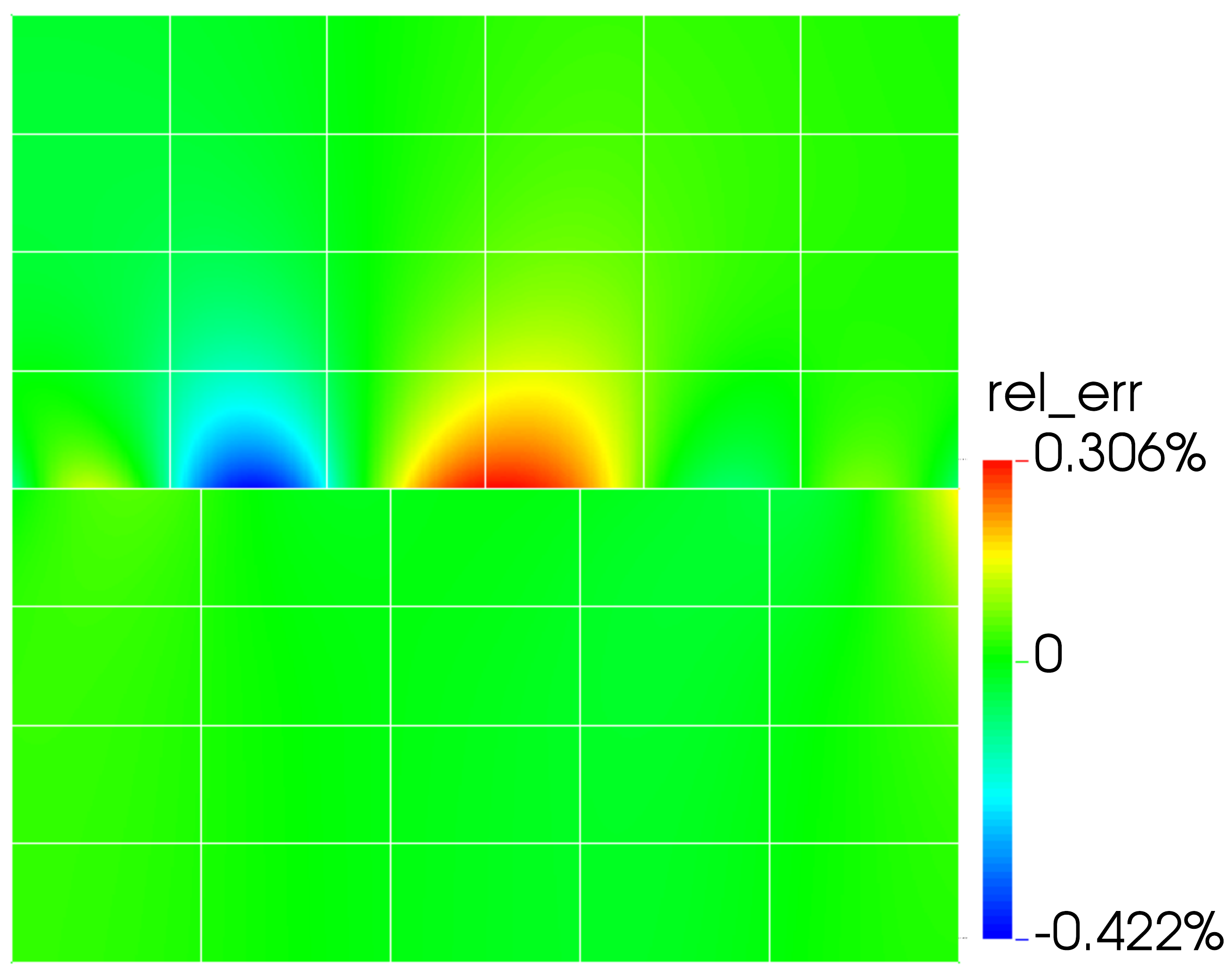}}
    \subfigure[Relative error for $\sigma_{yy}$]{\includegraphics[width=0.4\textwidth]{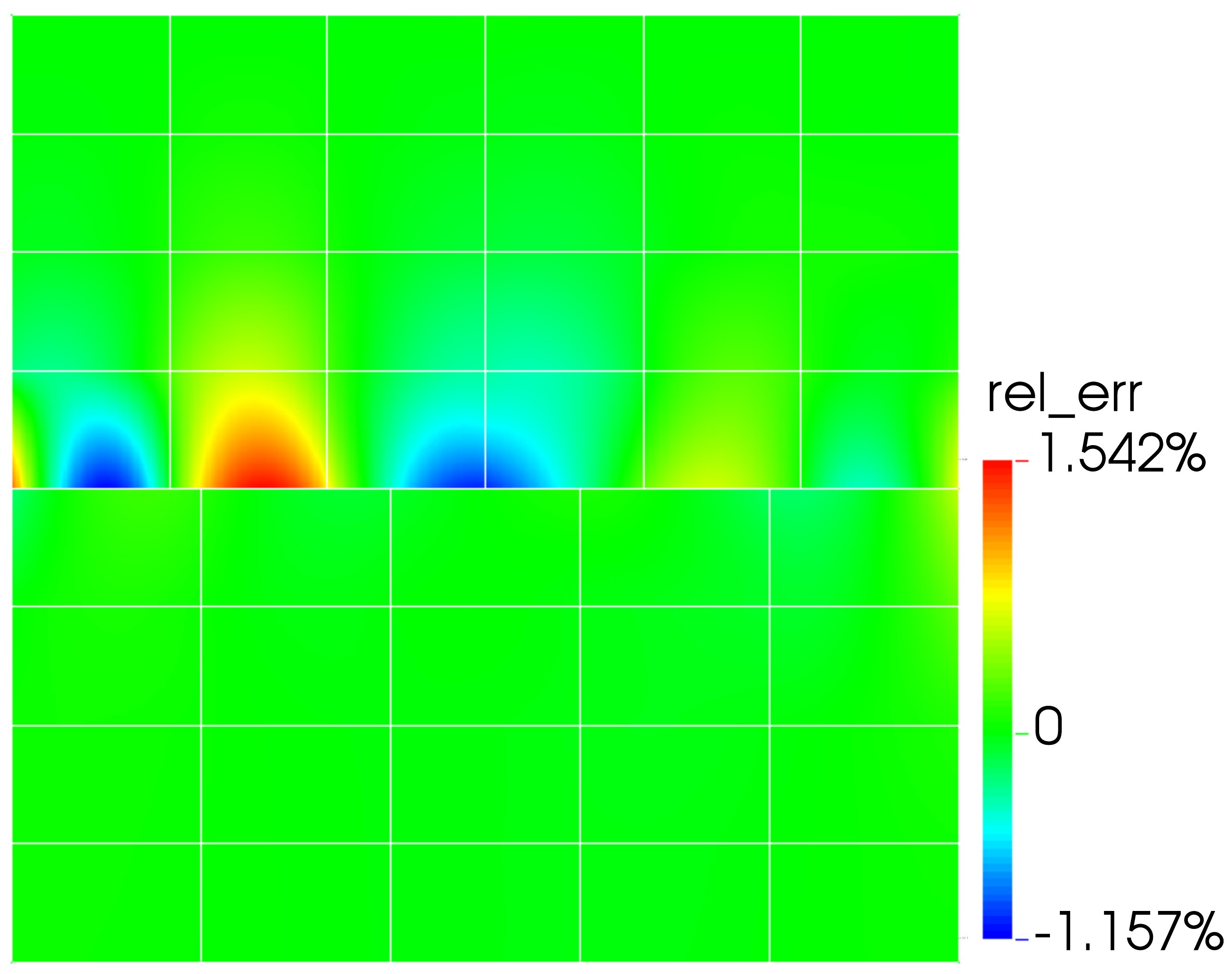}}
	\caption{Contact between two blocks: relative errors for the displacement field $u_y$ and for the stress field $\sigma_{yy}$. The biased skew-symmetric Nitsche's method is employed with the mesh given in \fref{fig_contact_patchtest}.}
	\label{Fig_contact_test_error_biased_skew}
\end{figure}

\begin{figure}[htbp]
	\centering
	\subfigure[Relative error for $u_y$]{\includegraphics[width=0.4\textwidth]{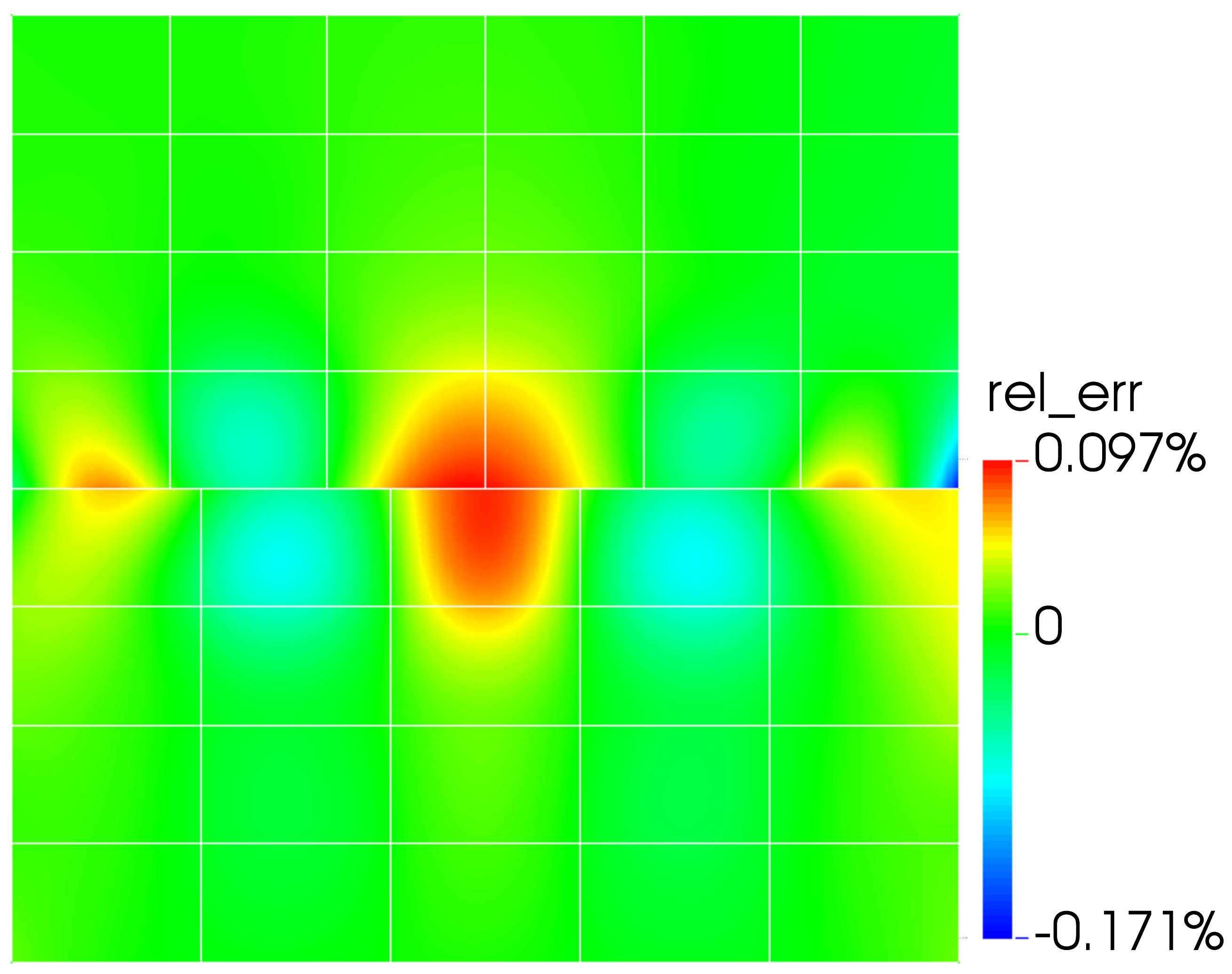}}
    \subfigure[Relative error for $\sigma_{yy}$]{\includegraphics[width=0.4\textwidth]{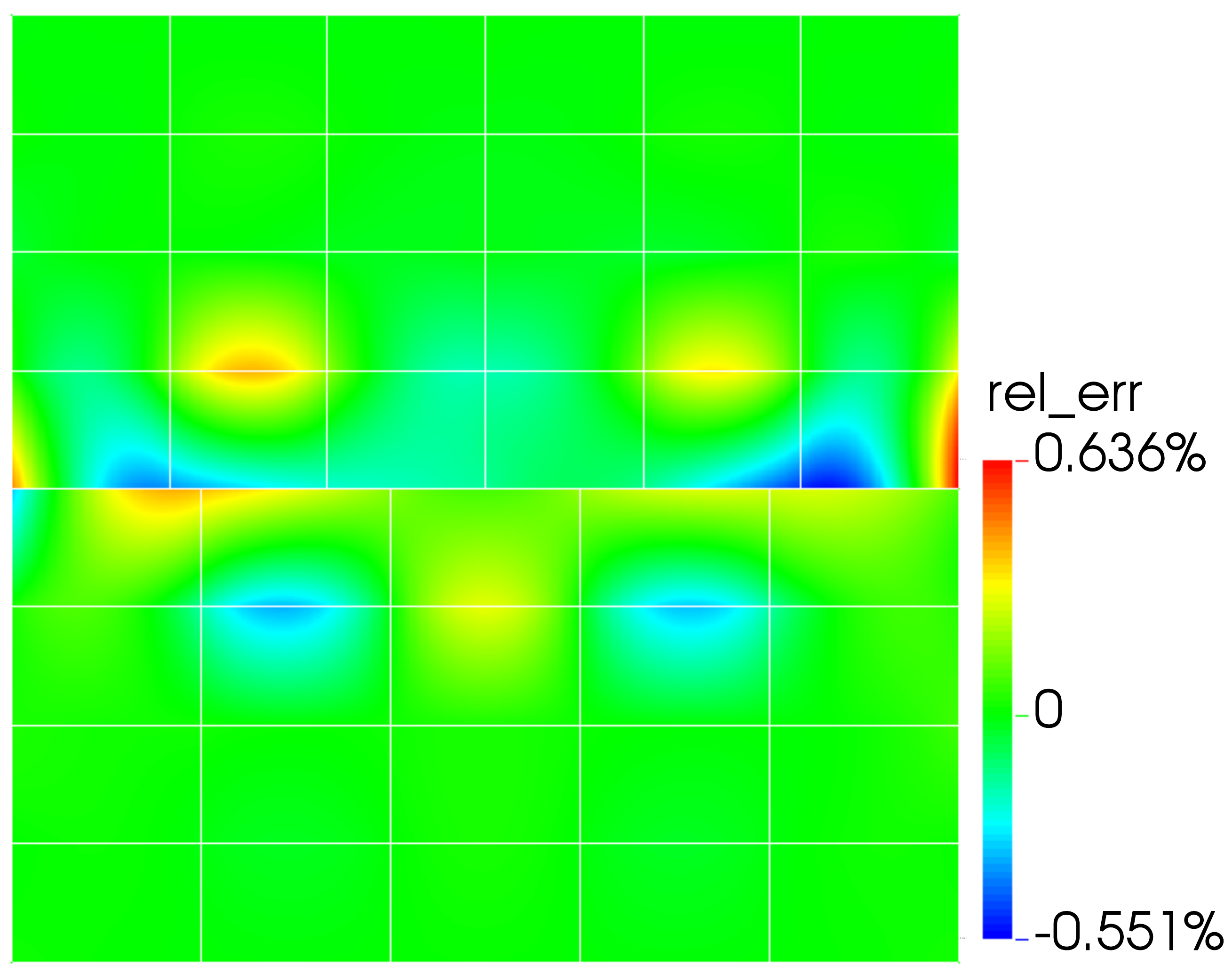}}
	\caption{Contact between two blocks: relative errors for the displacement field $u_y$ and for the stress field $\sigma_{yy}$. The unbiased skew-symmetric Nitsche's method is employed with the mesh given in \fref{fig_contact_patchtest}.}
	\label{Fig_contact_test_error_unbiased_skew}
\end{figure}

We use and compare both the biased Nitsche's contact formulation,
which corresponds to subsection \ref{subsub:biased},
and the unbiased one, which corresponds to subsection \ref{subsub:unbiased}.
Bi-quadratic basis functions are employed,
and only $\gamma_0 = 2 \lambda^{h,\mathrm{MAX}}$ is used.
\fref{Fig_contact_test_error_biased_skew} and \fref{Fig_contact_test_error_unbiased_skew} demonstrate the distribution of the relative errors for biased and unbiased skew-symmetric formulation, respectively.
For the biased formulation visible errors appear on the slave patch (the upper patch),
while for the unbiased formulation they appear on both patches.
We list the relative errors for $u_y$ and $\sigma_{yy}$ in \Tref{taylor_contact_tab},
showing that for both biased and unbiased versions the accuracy is comparable,
which supports the conclusion in \cite{mlika2017unbiased}.

\begin{table}[htbp]
\centering
\caption{Contact between two blocks: ranges of relative errors for $u_y$ and $\sigma_{yy}$.} 
\label{taylor_contact_tab}

\medskip

\begin{tabular}{c c cc}
 \hline
 Formulation type &  Nitsche's formulation type & Relative errors for $u_y$ & Relative errors for $\sigma_{yy}$ \\ \hline \noalign{\smallskip}
Biased &  Standard &  -0.416\% $\sim$ 0.311\% & -1.141\% $\sim$ 1.518\%  \\
Biased &  Skew-symmetric &  -0.422\% $\sim$ 0.306\% & -1.157\% $\sim$ 1.542\%  \\
 \hline \noalign{\smallskip}
Unbiased &  Standard &  -0.279\% $\sim$ 0.384\% & -2.816\% $\sim$ 3.083\%  \\
Unbiased &  Skew-symmetric &  -0.171\% $\sim$ 0.097\% & -0.551\% $\sim$ 0.636\%  \\
    \hline
\end{tabular}
\end{table}

\subsubsection{Self-contact of a 3D clip}

\begin{figure}[htbp]
\centering
\def\svgwidth{1\columnwidth}
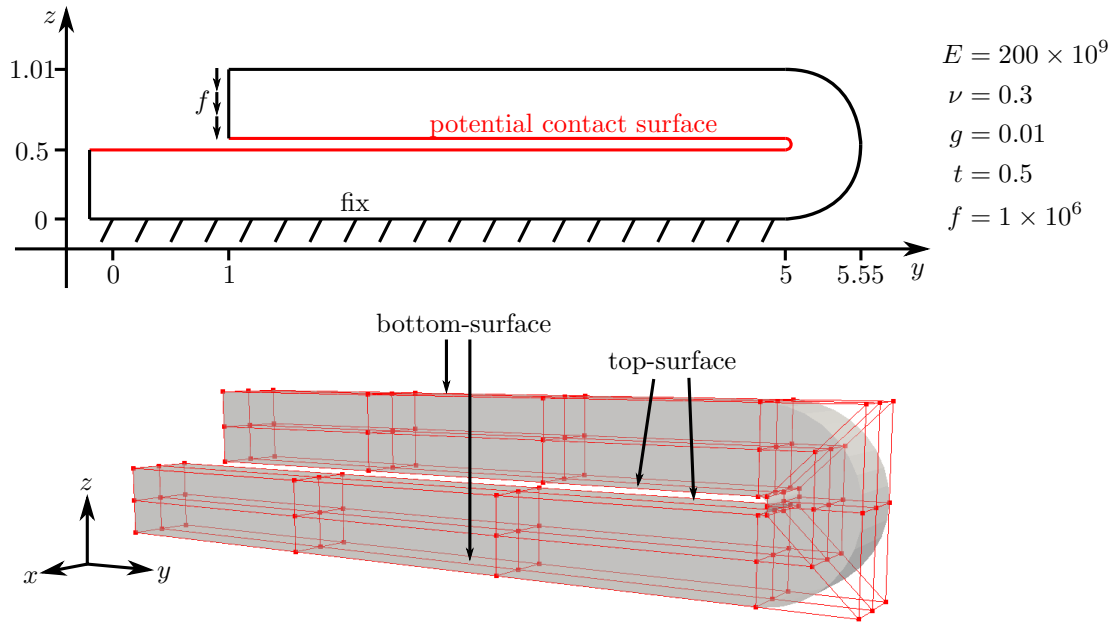
\caption{A 3D clip model with thickness $t$ and a small gap $g$. Its initial control points and control net are also shown. The values of the physical parameters are provided.}
\label{fig_clip_model}
\end{figure}

Here we present a clip model to illustrate the skew-symmetric Nitsche's formulation for self-contact in 3D, see \ref{subsub:unbiased}.
The model in shown in \fref{fig_clip_model}.
It is subjected to a uniform surface traction $f$ on its end.
The surface colored in red is defined as the potential contact surface,
the gap between the contact surface is $g=0.01$.
The control mesh in 3D is also shown below,
from which it is noticed that the contact surface is actually the top-surface itself,
thus this is a top-surface to top-surface self-contact problem.
Once again we adopt quadratic NURBS basis functions for the model and for the displacement field as well,
and use 3 quadrature points for each element boundary along the contact surface.
The value of $\gamma_0 = 2 \lambda^{h,\mathrm{MAX}}$ comes from the generalized eigenvalue problem (Eq. \eqref{eigenvalue}).

\begin{figure}[htbp]
	\centering
	\subfigure[Skew-symmetric Nitsche's self-contact formulation: 1st iteration]{\includegraphics[width=0.7\textwidth]{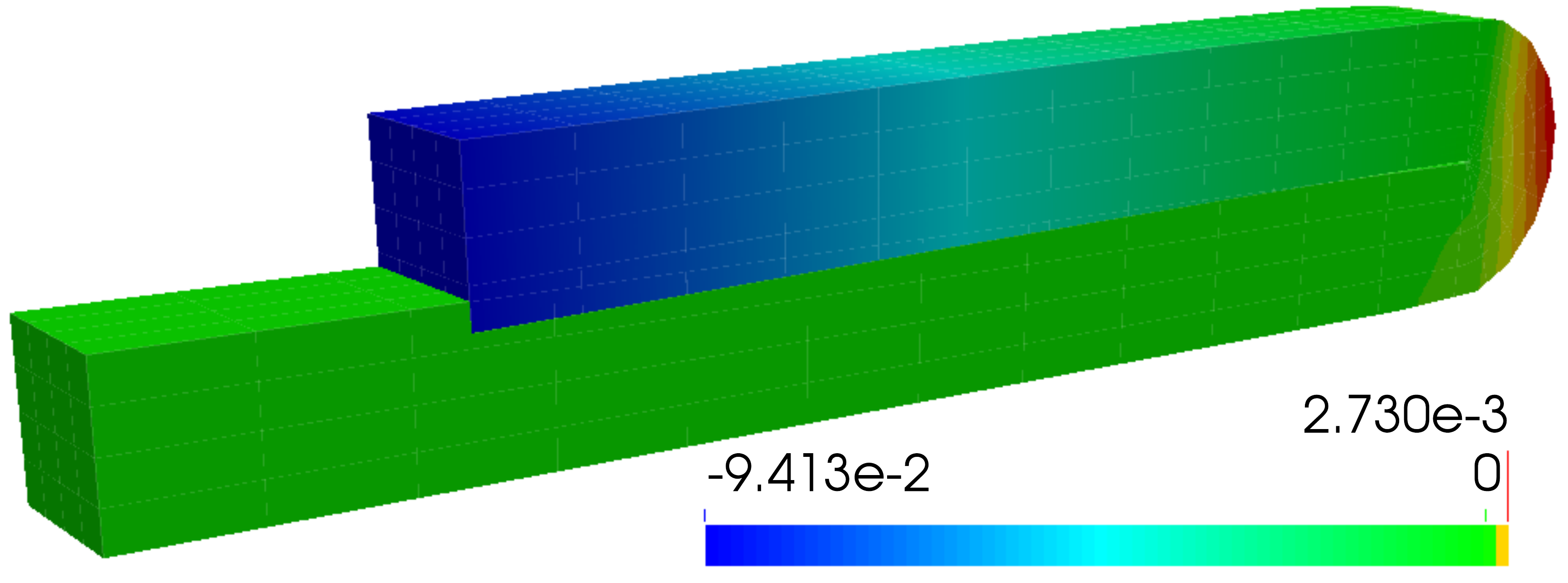}}
    \subfigure[Skew-symmetric Nitsche's self-contact formulation: converged solution]{\includegraphics[width=0.7\textwidth]{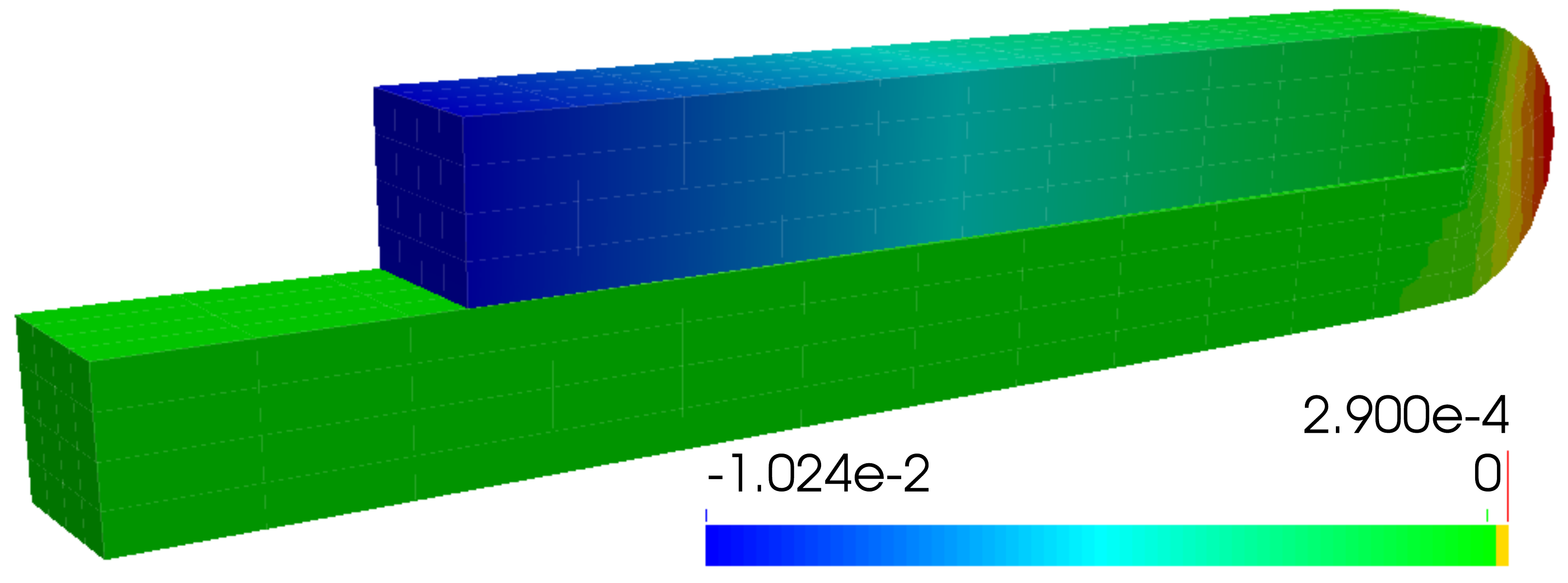}}
    \subfigure[\emph{ABAQUS} reference solution. $u_z^{\textrm{MAX}}=2.991\times 10^{-4}$, $u_z^{\textrm{MIN}}=-1.022\times 10^{-2}$.]{\includegraphics[width=0.8\textwidth]{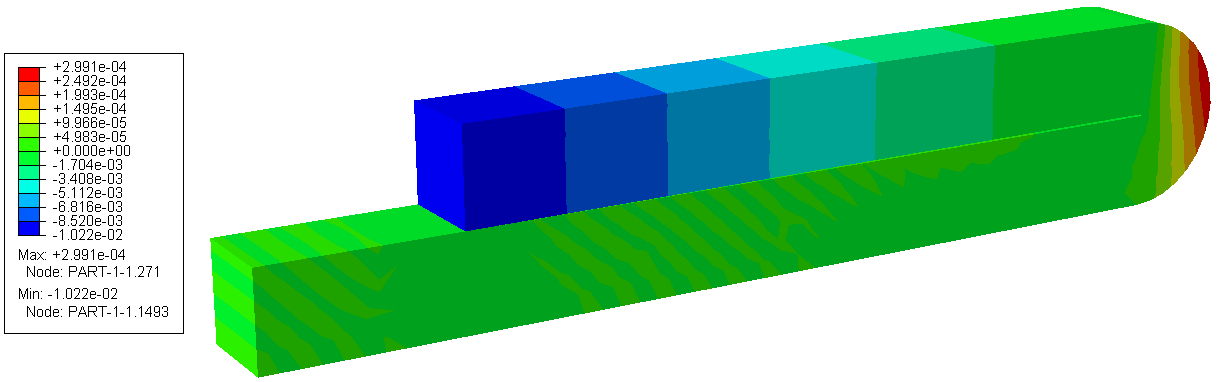}}
	\caption{Contour plot of vertical displacement field $u_z$.}
    \label{3d_selfcontact_uz}
\end{figure}

The contour plot of the vertical displacement field $u_z$ is shown in \fref{3d_selfcontact_uz}.
When the iterative solving starts,
the contact surface is penetrated,
the largest vertical displacement is $u_z=-0.09413$ while the contact gap is $g=0.01$,
see \fref{3d_selfcontact_uz} (a)
(this is what the deformation would be if no contact were taken into account).
When semi-smooth Newton procedure has converged,
the largest vertical displacement becomes $u_z=-0.01024$,
see \fref{3d_selfcontact_uz} (b).
In \fref{3d_selfcontact_uz} (c) we adopt the results obtained by \emph{ABAQUS} using a sufficient number of 3D solid elements:
the vertical displacement distribution obtained by Nitsche's contact formulation is in good agreement with the solution provided by \emph{ABAQUS}.
The relative errors of the maximum positive and negative vertical displacements are $0.196\%$ and $-3.042\%$ respectively.

\section{Conclusions}\label{section_conclusion}

We presented a systematic way to derive Nitsche's formulations for different kind of boundary and interface conditions,
and studied this technique in the context of isogeometric analysis (IGA) discretization.
We recover different variants of Nitsche's method, for different values of the \emph{Nitsche parameter} $\theta$,
and then focused on the skew-symmetric variant, namely $\theta=-1$. This variant
is appealing because it does not need a stabilization term  for linear boundary/interface conditions,
and is robust w.r.t. the stabilization parameter for non-linear boundary/interface conditions.
Several numerical studies were performed to %study and
illustrate the behavior of Nitsche's method,
especially the skew-symmetric variant.  %Nitsche formulation.
From the numerical results we can state the observations below:%se numerical studies we made the following observations:

$\bullet$ The skew-symmetric formulation is effective to impose Dirichlet displacement boundary conditions in small strain elasticity %for second order problems
as well as the symmetric rotational boundary conditions for Kirchhoff-Love plates. %fourth order problems.
The skew-symmetric formulation is parameter-free in this context and achieves good accuracy:
for the circular patch test (\fref{Fig_patch_test_energy}) and the Kirchhoff plate (\fref{kirchoff_result}) we observe the predicted optimal convergence rates in the energy norm.

$\bullet$ For patch coupling in statics,
the skew-symmetric Nitsche's formulation is still parameter-free.
Condition numbers for the global stiffness matrix are far better than for standard Nitsche,
and only slightly above the conforming setting.
They are also almost independent of the mesh size, and basis functions orders:
see \Tref{couple_influence_condition_number} in Section \ref{pach_coupling_statics}.

$\bullet$ For patch coupling in modal analysis,
Nitsche's formulation increases the number of ``outlier'' frequencies.
The reason is believed to be that Nitsche's formulation introduces additional highly localized eigenmodes,
and the positions of these newly added eigenmodes just locate at the coupled interfaces.

$\bullet$ For contact problems in linear elasticity,
the skew-symmetric Nitsche's formulation behaves more robustly than the standard Nische formulation regarding the value of the stabilization parameter.
Nitsche's method can properly impose the contact conditions,
and predict the pressure distribution with respect to the contact width. Moreover, it allows an unbiased variant that can be more appealing for self-body and multi-body contact.

\section*{Acknowledgments}

The authors thank the two anonymous referees for their constructive comments that helped to improve the paper. They thank also Patrick Le Tallec for inspiring discussions and comments that helped to improve the presentation.
Qingyuan Hu is funded by China Scholarship Council.
St\'ephane Bordas thanks the financial support of the European Research Council Starting Independent Research Grant (ERC Stg grant agreement No. 279578).
%entitled Towards real time multiscale simulation of cutting in non-linear materials with applications to surgical simulation and computer guided surgery.
St\'ephane Bordas is also grateful for the support of the Fonds National de la Recherche Luxembourg FNRS-FNR grant INTER/FNRS/15/11019432/EnLightenIt/Bordas.
Franz Chouly thanks R\'egion Bourgogne Franche-Comt\'e for funding (``Convention R\'egion 2015C-4991. Mod\`eles math\'ematiques et m\'ethodes num\'eriques pour l'\'elasticit\'e non-lin\'eaire''), as well as the Centre National de la Recherche Scientifique (``Convention 232789 DEFI InFIniTI 2017 - Projet MEFASIM'').
%Please add your funding here.

%\section*{References}

%\bibliographystyle{model1-num-names} % this is the one we used for the 1st version

% https://www.sharelatex.com/learn/Bibtex_bibliography_styles
%\bibliographystyle{abbrv}
%\bibliographystyle{plain}
%\bibliographystyle{alpha}
%\bibliographystyle{apalike}
\bibliographystyle{acm}

\bibliography{Nitsche}

\end{document}